\pdfoutput=1
\documentclass[english]{jnsao}
\usepackage[utf8]{inputenc}

\usepackage{amsmath}
\usepackage{amsfonts}

\usepackage{enumitem}
\setlist[enumerate,1]{label={(\alph*)}}
\setlist[enumerate,2]{label={(\roman*)}}

\usepackage[ruled,vlined,linesnumbered,algosection,algo2e]{algorithm2e} 

\usepackage[nameinlink,capitalize]{cleveref}

\theoremstyle{definition}
\newtheorem{assum}[theorem]{Assumption}
\crefname{assum}{Assumption}{Assumptions}

\newcommand{\keywordsName}{Keywords}
\newcommand{\keywordsAnd}{~$\cdot$~}
\newcommand{\subclassName}{MSC (2020)}
\newcommand{\subclassAnd}{~$\cdot$~}
\newenvironment{keywords}{%
	\vspace*{0.2cm}%
	\small\noindent%
	{\normalfont\sectfont\nobreak\keywordsName\quad}%
}{}
\newenvironment{subclass}{%
	\vspace*{0.2cm}%
	\small\noindent%
	{\normalfont\sectfont\nobreak\subclassName\quad}%
}{}

\crefname{line}{Step}{Steps} 
\SetKw{KwOr}{or}

\newcommand{\amsmscLink}[1]{\href{http://www.ams.org/mathscinet/msc/msc2020.html?t=#1}{#1}}

\newcommand{\N}{\mathbb{N}}
\newcommand{\R}{\mathbb{R}}
\newcommand{\Rinf}{\overline{\R}} 			
\DeclareMathOperator*{\argmin}{\arg\min}
\DeclareMathOperator{\dom}{dom}
\DeclareMathOperator{\proj}{\Pi}
\DeclareMathOperator{\prox}{prox}
\DeclareMathOperator{\indicator}{\delta}
\DeclareMathOperator*{\minimize}{minimize}
\DeclareMathOperator{\stt}{subject~to}
\newcommand{\func}[3]{#1 \colon #2 \to #3}
\newcommand{\ffunc}[3]{#1 \colon #2 \rightrightarrows #3}
\newcommand{\innprod}[2]{\langle #1, #2 \rangle}
\newcommand{\norm}[1]{\Vert #1 \Vert}
\newcommand{\normalcone}{N}
\newcommand{\limnormalcone}{N}

\DeclareMathOperator\dist{dist}
\DeclareMathOperator\conv{conv}
\DeclareMathOperator\epi{epi}
\DeclareMathOperator\gph{gph}
\DeclareMathOperator\closedball{\mathbb{B}}

\renewcommand{\d}{\mathrm{d}}
\newcommand{\regsub}{\widehat{\partial}}

\newcommand{\conj}{\ast}
\newcommand{\downto}{\downarrow}
\newcommand{\downtoeq}{\to}
\newcommand{\nto}{\nrightarrow}
\newenvironment{mybox}{}{}

\DeclareMathAlphabet{\mathpzc}{OT1}{pzc}{m}{it}
\newcommand\oo{\mathpzc{o}}

\numberwithin{equation}{section}

\newcommand{\XX}{\R^n}
\newcommand{\YY}{\R^m}
\newcommand{\toattentive}[1]{\:
	\overset{%
		\text{\raisebox{1.2ex}{\smash{\scalebox{0.8}{$#1$}}}}%
	}{%
		\text{\raisebox{0.2ex}{\smash{$\to$}}}%
	}
\:}
\newcommand{\LL}{\mathcal{L}}
\newcommand{\LLslack}{\LL^{\mathrm{S}}}
\newcommand{\Ybounded}{Y}
\newcommand{\FBO}{\mathbf{O}}

\newcommand{\TheAuthorADM}{Alberto De~Marchi}
\newcommand{\TheEmailADM}{alberto.demarchi@unibw.de}
\newcommand{\TheAffiliationADM}{%
	University of the Bundeswehr Munich,
	Department of Aerospace Engineering,
	Institute of Applied Mathematics and Scientific Computing,
	85577 Neubiberg, Germany%
}
\newcommand{\TheAuthorPM}{Patrick Mehlitz}
\newcommand{\TheEmailPM}{mehlitz@uni-marburg.de}
\newcommand{\TheAffiliationPM}{%
	Philipps-Universität Marburg,
	Department of Mathematics and Computer Science,
	35032 Marburg, Germany%
}
\newcommand{\TheTitle}{Local properties and augmented Lagrangians in fully nonconvex composite optimization}
\newcommand{\TheShortTitle}{Local properties in fully nonconvex composite optimization}
\newcommand{\TheKeywords}{%
	Augmented Lagrangian framework\keywordsAnd%
	Composite nonconvex optimization\keywordsAnd%
	Error bounds\keywordsAnd%
	Local convergence properties\keywordsAnd%
	Second-order variational analysis%
}
\newcommand{\TheAMSsubj}{%
	\amsmscLink{49J52}\subclassAnd
	\amsmscLink{49J53}\subclassAnd
	\amsmscLink{65K10}\subclassAnd
	\amsmscLink{90C30}\subclassAnd
	\amsmscLink{90C33}
}

	\author{%
		\TheAuthorADM\thanks{%
			\TheAffiliationADM.
		\email{\TheEmailADM},
		\orcid{0000-0002-3545-6898}.
		}\and%
		\TheAuthorPM\thanks{%
			\TheAffiliationPM.
		\email{\TheEmailPM},
		\orcid{0000-0002-9355-850X}.
		}
	}%

\title{\TheTitle}
\shorttitle{\TheShortTitle}
\shortauthor{De Marchi and Mehlitz}
\date{\ISOToday}

\manuscriptlicense{CC-BY 4.0}

\manuscriptsubmitted{2023-09-06}
\manuscriptaccepted{2024-05-14}
\manuscriptvolume{5}
\manuscriptnumber{12235}
\manuscriptyear{2024}
\manuscriptdoi{10.46298/jnsao-2024-12235}

\manuscripteprinttype{arXiv}
\manuscripteprint{2309.01980}

\begin{document}
	
\maketitle

\begin{center}
	``\emph{In fact the great watershed in optimization isn't between linearity and nonlinearity, but between convexity and nonconvexity}.''\\
	\hfill--- R. T. Rockafellar \cite{rockafellar1993lagrange}
\end{center}

\medskip

\begin{abstract}
	A broad class of optimization problems can be cast in composite form,
	that is,
	considering the minimization
	of the composition of a lower semicontinuous function
	with a differentiable mapping.
	This paper investigates the versatile template of composite optimization 
	without any convexity assumptions.
	First- and second-order optimality conditions are discussed.
	We highlight the difficulties that stem from the lack of convexity
	when dealing with necessary conditions in a Lagrangian framework
	and when considering error bounds.
	Building upon these characterizations,
	a local convergence analysis is delineated for a recently developed augmented Lagrangian method,
	deriving rates of convergence in the fully nonconvex setting.
	
	\begin{keywords}
		\TheKeywords
	\end{keywords}
	
	\begin{subclass}
		\TheAMSsubj
	\end{subclass}
\end{abstract}

\section{Introduction}\label{sec:introduction}

In this paper, we are concerned with finite-dimensional optimization problems of the form
\begin{equation}
	\minimize_{x}
	{}\quad{}
	\varphi(x) \coloneqq f(x) + g( c(x) ) ,
	\tag{P}\label{eq:P}
\end{equation}
where $x\in\XX$ is the decision variable, $\func{f}{\XX}{\R}$ and $\func{c}{\XX}{\YY}$ are smooth mappings, 
and $\func{g}{\YY}{\Rinf \coloneqq \R \cup \{\infty\}}$ is merely proper and lower semicontinuous.
The data functions $f$ and $g$ are allowed to be nonconvex mappings, the nonsmooth cost $g$ is not necessarily continuous nor real-valued, 
and the mapping $c$ is potentially nonlinear.
Such setting of nonsmooth nonconvex composite optimization was named
\emph{generalized} \cite{rockafellar2022convergence} (or \emph{extended} \cite{rockafellar2000extended}) nonlinear programming by Rockafellar,
and it is well known that the model problem \eqref{eq:P} covers numerous applications from 
signal processing, sparse optimization, compressed sensing, machine learning, and disjunctive programming.
Let us mention that the situation where $g$ is \emph{convex} has been intensively studied in
the literature, see e.g.\ \cite{rockafellar2000extended,rockafellar2022augmented,rockafellar2022convergence}.
We do \emph{not} make such an assumption and consider the far more challenging general setup.

Our motivation behind the potential nonconvexity of $g$
is driven by applications from sparse or low-rank optimization.
Although the (convex) $\ell_1$- and nuclear norm are known to
promote sparse or low-rank behavior,
solutions are often not sparse \emph{enough} in certain settings.
In order to overcome this issue, one can rely on the $\ell_0$-quasi-norm
or the matrix rank as ``regularizers'', which are discontinuous functions.
Intermediate choices for $g$ like the $\ell_q$-quasi-norm or the $q$-Schatten-quasi-norm for $q\in(0,1)$
are nonconvex but uniformly continuous, and have turned out to work well in numerical practice.
Another driving force behind this work comes from disjunctive programming,
in particular from the observation that constraints can be naturally formulated
in a function-in-set format whereby sets are nonconvex yet simple (to project onto).
The template \eqref{eq:P} lends itself to capture this scenario,
taking full advantage of $g$ as the indicator of a nonconvex set,
comprising structures typical for, e.g., complementarity, switching, and vanishing constraints,
see the classical monographs \cite{LuoPangRalph1996,OutrataKocvaraZowe1998}
and the more recent contributions \cite{Gfrerer2014,mehlitz2020}.

The possibility to include constraints in the model problem \eqref{eq:P}
becomes apparent with a direct reformulation.
Introducing an auxiliary variable $z \in \YY$, \eqref{eq:P} can be equivalently rewritten in the form
\begin{equation}
	\minimize_{x,z}
	{}\quad{}
	f(x) + g(z)
	{}\qquad{}
	\stt
	{}\quad{}
	c(x) - z = 0 ,
	\tag{P$_{\text{S}}$}\label{eq:ALMz:P}
\end{equation}
which has a separable objective function, without nontrivial compositions,
and explicitly includes some equality constraints.
An analogous template has been studied in \cite{demarchi2023constrained},
demonstrating its modeling versatility,
mostly enabled by accepting potentially nonconvex $g$.

A fundamental technique for solving constrained optimization problems is the
augmented Lagrangian (AL) framework,
which can effortlessly handle nonsmooth objectives, see e.g.\ 
\cite{bolte2018nonconvex,chen2017augmented,demarchi2024implicit,demarchi2023constrained,DhingraKhongJavanovic2019,HangSarabi2021,hang2023convergence,
rockafellar2022augmented,rockafellar2022convergence,sabach2019lagrangian}
for some recent contributions,
and \cite{bertsekas1996constrained,birgin2014practical,conn1991globally} for some fundamental
literature which addresses the setting of standard nonlinear programming.
Particularly, Rockafellar extended the approach in \cite{rockafellar2022augmented,rockafellar2022convergence}
to the broad setting of \eqref{eq:P} with $g$ convex,
relying on some local duality to build a connection with the
proximal point algorithm (PPA), see \cite{rockafellar1973dual,rockafellar1976augmented}.
Embracing the fully nonconvex setting,
we are interested here in investigating AL methods for
generalized nonlinear programming without any convexity assumption.
Although the shifted-penalty approach underpinning the seminal
``method of multipliers'' still applies in our setting,
it appears more difficult to leverage the perspective of PPA.
Moreover, the nonconvexity of $g$ leads to a lack of regularity,
as its proximal mapping is potentially set-valued.
Here, we seek a better understanding of the variational properties
of \eqref{eq:P} and the convergence guarantees of AL methods for this class of problems.
Building upon the global characterization in \cite{demarchi2024implicit},
we will focus on second-order optimality conditions and local analysis,
portraying a convergence theory for the fully nonconvex setting including rates-of-convergence results.

The following blanket assumptions are considered throughout, without further mention.
Technical definitions are given in \cref{sec:preliminaries}.
\begin{mybox}
	\begin{assum}\label{ass:P}
		The following hold in \eqref{eq:P}:
		\begin{enumerate}[label=(\roman*)]
			\item\label{ass:f}%
			$\func{f}{\XX}{\R}$ and $\func{c}{\XX}{\YY}$ are twice continuously differentiable;
			\item\label{ass:g}%
			\(\func{g}{\YY}{\Rinf}\) is proper, lower semicontinuous, and prox-bounded;
			\item\label{ass:phi}%
			\(\inf\varphi\in\R\).
		\end{enumerate}
	\end{assum}
\end{mybox}

The prox-boundedness assumption on $g$ in \cref{ass:P}\ref{ass:g} is included to ensure that,
for some suitable parameters,
the proximal mapping of $g$ is well-defined,
and so is the overall numerical scheme considered in this work.
However, such stipulation is not necessary.
As suggested in \cite{rockafellar2022convergence},
a ``trust region'' can be specified to localize the proximal minimization and to support its attainment.
While avoiding the artificial unboundedness potentially introduced
by relaxing the composition constraint,
this localization would affect some algorithmic and global aspects,
but not the local behavior and properties we are interested in here.

\subsection{Related Work and Contributions}\label{sec:relatedWorkContributions}

Since its inception \cite{hestenes1969multiplier,powell1969method,rockafellar1973dual},
the AL framework has been extensively investigated and developed
\cite{birgin2014practical,conn1991globally},
also in infinite dimensions \cite{kanzow2018augmented}.
It was soon recognized that, in the convex setting,
the method of multipliers can be associated to
the PPA applied to a dual problem, see \cite{rockafellar1976augmented}.
Following this pattern,
local convexity enabled by some second-order optimality conditions
allowed to reconcile the AL scheme with an application of the PPA,
and thereby establishing convergence, beyond the convex setting
\cite{bertsekas1996constrained,rockafellar2022augmented,rockafellar2022convergence}.
However, when it comes to \emph{local} convergence properties,
available results remain confined to the case with $g$ convex.

\paragraph{Contributions}

With this work,
we extend recent results by Rockafellar from
\cite{rockafellar2022augmented,rockafellar2022convergence},
where composite optimization problems with convex $g$ are considered,
to the more general setting.
In particular, we study the implicit AL method from \cite{demarchi2024implicit}
and characterize its local convergence behavior.
Particularly, under suitable conditions, we show convergence of the full sequence with
linear or superlinear rate in \cref{thm:rate_convergence}.
To proceed, we make use of problem-tailored second-order conditions 
which have been developed recently in \cite{BenkoMehlitz2023}.
Moreover, the Lagrange multiplier has to be locally unique,
see \cref{lem:uniqueness_of_multiplier_CQ}.

Sparsity-promoting terms and nonconvex constraint sets have turned out
to work well in the AL framework---at least from a global
perspective, e.g.\ in \cite{demarchi2023constrained,JiaKanzowMehlitzWachsmuth2023}.
We are also interested in local properties now,
with a focus on the numerical method proposed in \cite{demarchi2024implicit},
which favorably avoids the use of slack variables,
see \cite{benko2021implicit} for a recent study.

Local fast convergence of an AL method in composite optimization
has been considered from the viewpoint of
variational analysis in the recent paper
\cite{HangSarabi2021} in the context where $g$ is a continuous, piecewise
quadratic, convex function. This allows for a unified analysis as the
standard second-order sufficient condition already gives the necessary error bound condition
(due to the result from \cref{lem:subderivative_vs_graphical_derivative} and
the analysis in \cite{MohammadiMordukhovichSarabi2022,rockafellar1998variational}).
A recent analysis of the local convergence of AL methods for \eqref{eq:P} is that in \cite{hang2023convergence},
restricted to convex $g$,
which considers second-order sufficient conditions 
and establishes Q-linear convergence of the primal-dual sequence,
without assuming any constraint qualification (CQ).

Our analysis also took inspiration from 
\cite{BoergensKanzowSteck2019,KanzowSteck2018,steck2018dissertation}
where, among other things, the local analysis of
AL methods for (smooth) optimization
problems with geometric constraints of type $c(x)\in K$ for some
closed, convex set $K$ is considered in a Banach space setting.
We, at least roughly, follow the arguments in \cite{steck2018dissertation}
and (apart from the fact that we are working in a fully finite-dimensional setting)
generalize the findings therein to nonsmooth composite problems.

Let us point out that desirable local convergence properties of AL methods
in standard nonlinear programming can be guaranteed with no more than a
second-order sufficient condition, 
i.e., no additional CQ is necessary, see \cite{FernandezSolodov2012},
and the second-order sufficient condition 
can even be replaced by a weaker noncriticality assumption on
the involved multipliers, as shown later in \cite{IzmailovKurennoySolodov2015}.
One reason for this behavior is the inherent (convex) polyhedrality of the involved
constraint set, see \cref{ex:NLP} below, which also gives (convex) polyhedrality of the associated set
of Lagrange multipliers. 
The fact that polyhedrality comes along with certain stability properties (in the sense of error bounds)
is well known from the seminal papers \cite{Robinson1981,WalkupWets1969}.
In the general, nonpolyhedral situation, such a result is not likely
to hold, see \cite{KanzowSteck2018}, and an additional CQ might be necessary.
Exemplary, this has been visualized in the papers \cite{HangMordukhovichSarabi2022,KanzowSteck2019}
where AL methods in second-order cone programming have been investigated.
In order to obtain convergence rates, the authors do not only postulate the validity of a
second-order condition, but make use of an addition assumption.
In \cite{KanzowSteck2019}, the authors exploit the strict Robinson condition 
(which guarantees uniqueness of the underlying Lagrange multiplier) while in
\cite{HangMordukhovichSarabi2022}, a certain multiplier mapping is assumed to be calm
while, at the point of interest, uniqueness of the Lagrange multiplier is also needed.

\paragraph{Roadmap}

The remainder of the paper is organized as follows.
\cref{sec:preliminaries}
provides some preliminary results from variational analysis 
and generalized differentiation.
\cref{sec:compositeOptimization} is dedicated to the
investigation of first-order necessary and second-order sufficient 
optimality conditions in nonconvex composite optimization.
Furthermore, we comment on a reasonable choice for an AL function
and investigate error bounds for a system of necessary optimality conditions
associated with \eqref{eq:P}.
In \cref{sec:ALM}, we first introduce the AL method of our interest before presenting
some global convergence results which complement the analysis provided in \cite{demarchi2024implicit}.
Then, local convergence results are presented. 
We start by clarifying the existence and convergence of minimizers for the associated
AL subproblems before focusing on the derivation of rates-of-convergence results.
\cref{sec:examples} illustrates our findings by means of two
exemplary settings: sparsity-promoting nonlinear optimization and 
complementarity-constrained optimization.
The paper closes with some concluding remarks in \cref{sec:conclusions}.

\section{Preliminaries}\label{sec:preliminaries}

This section provides
notation, preliminaries, and known facts based on 
\cite{mordukhovich2018,rockafellar1998variational}, 
with some additional basic results.

With $\R$ and $\Rinf \coloneqq \R \cup \{\infty\}$ we indicate the real and extended-real line, respectively.
The set of natural numbers is denoted by $\N$.
We equip the appearing Euclidean spaces 
possessing the standard Euclidean inner product $\innprod{\cdot}{\cdot}$
with the associated Euclidean norm $\norm{\cdot}$.
In product spaces, we make use of the associated maximum norm.
With $\closedball_r(x)$ we indicate the closed ball centered at $x\in\R^n$ with radius $r>0$.
Given a set $A\subseteq\R^n$, we use 
$x+A\coloneqq A+x\coloneqq\{a+x\in\XX\,|\,a\in A\}$ for brevity.
The notation $\{a^k\}_{k\in K}$ represents a sequence indexed by elements of the set $K\subseteq\N$,
and we write $\{a^k\}_{k\in K} \subseteq A$ to indicate that
$a^k \in A$ for all indices $k\in K$.
Whenever clear from context, we may simply write $\{a^k\}$ to indicate $\{a^k\}_{k\in\N}$.
Notation $a^k\to_K x$ ($a^k\to x$) is used to express convergence of $\{a^k\}_{k\in K}$ (of $\{a^k\}$) to $x$.
If $n = 1$, we use $\{a_k\}_{k\in K}$ and $\{a_k\}$ to emphasize that we are dealing with sequences of scalars.
We will adopt the \emph{little-o} notation for asymptotics:
given sequences $\{a_k\}$ and $\{\varepsilon_k\}\subset(0,\infty)$,
we write $a_k \in \oo(\varepsilon_k)$ to indicate that $\lim_{k\to\infty} {|a_k|}/{\varepsilon_k} = 0$.

A function $\func{f}{\XX \times \YY}{\Rinf}$ with values $f(x,z)$ is 
\emph{level-bounded in $x$ locally uniformly in $z$} 
if for each $\alpha \in \R$ and $\bar{z} \in \YY$ there exists $\varepsilon > 0$ 
such that the set $\{ (x,z)\in\XX\times\YY \,\vert\, f(x,z) \leq \alpha, \norm{z - \bar{z}} \leq \varepsilon \}$ is bounded.
The \emph{effective domain} of a function $\func{h}{\YY}{\Rinf}$ is denoted by $\dom h \coloneqq \{ z\in\YY \,\vert\, h(z) < \infty \}$.
The set $\epi h\coloneqq \{(z,\alpha)\in\YY\times\R\,|\,\alpha\geq h(z)\}$ 
is called the \emph{epigraph} of $h$.
We say that $h$ is \emph{proper} if $\dom h \neq \emptyset$
and \emph{lower semicontinuous} if $h(\bar{z}) \leq \liminf_{z\to\bar{z}} h(z)$ for all $\bar{z} \in \YY$.
Note that $h$ is lower semicontinuous if and only if $\epi h$ is closed.
Given a point $\bar{z}\in\dom h$, 
we may avoid to assume $h$ continuous and instead appeal to 
$h$-\emph{attentive} convergence of a sequence $\{ z^k \}$, 
denoted as $z^k \toattentive{h} \bar{z}$ and given by $z^k \to \bar{z}$ 
with $h(z^k) \to h(\bar{z})$.
For some real number $\lambda\geq 1$,
we refer to $h$ as \emph{positively homogeneous of degree $\lambda$}
if $h(\alpha y)=\alpha^\lambda h(y)$ holds for all
$y\in\YY$ and real numbers $\alpha>0$.
The \emph{conjugate function} $\func{h^\conj}{\YY}{\Rinf}$ 
associated with (proper and lower semicontinuous) $h$ is defined by
\begin{equation*}
	h^\conj(y)
	\coloneqq
	\sup\limits_z\{ \innprod{y}{z} - h(z) \}
	.
\end{equation*}
We note that $h^\conj$ is a convex function by definition since it is a supremum of affine functions.

For a proper and lower semicontinuous function $\func{h}{\YY}{\Rinf}$, 
a point $\bar z\in \YY$ is called \emph{feasible} if $\bar z \in \dom h$.
A feasible point $\bar z \in \YY$ is said to be \emph{locally optimal}, or called a \emph{local minimizer}, 
if there exists $r > 0$ such that $h(\bar z) \leq h(z)$ holds for all feasible $z \in \closedball_r(\bar z)$.
Additionally, if this inequality holds for all feasible $z \in \YY$, then $\bar z$ is said to be \emph{(globally) optimal}.

We use the notation $\ffunc{\Gamma}{\XX}{\YY}$ to indicate a point-to-set function 
$\func{\Gamma}{\XX}{2^{\YY}}$. 
The set $\gph\Gamma\coloneqq\{(x,y)\in\XX\times\YY\,|\,y\in\Gamma(x)\}$ is called the \emph{graph} of $\Gamma$.
The set-valued mapping $\ffunc{\Gamma^{-1}}{\YY}{\XX}$ given by
$\gph\Gamma^{-1}\coloneqq\{(y,x)\in\YY\times\XX\,|\,(x,y)\in\gph\Gamma\}$
is referred to as the \emph{inverse} of $\Gamma$.
The set $\ker\Gamma\coloneqq \{x\in\XX\,|\,0\in\Gamma(x)\}$ is the \emph{kernel} of $\Gamma$.
Recall that $\Gamma$ is said to be a \emph{polyhedral} mapping if $\gph\Gamma$ can be represented
as the union of finitely many convex polyhedral sets.

\subsection{Proximal mappings}

Let $\func{h}{\YY}{\Rinf}$ be proper and lower semicontinuous.
Given a parameter value $\mu>0$, the \emph{proximal} mapping 
$\ffunc{\prox_{\mu h}}{\YY}{\YY}$ is defined by
\begin{equation}
	\prox_{\mu h}(z)
	{}\coloneqq{} 
	\argmin_{z^\prime} \left\{ h(z^\prime) + \frac{1}{2\mu}\norm{z^\prime-z}^2 \right\} .
	\label{eq:proxmapping}
\end{equation}
We say that $h$ is \emph{prox-bounded} if  
$h + \norm{\cdot}^2 / (2\mu)$ is bounded below on $\YY$ 
for some $\mu > 0$, see \cite[Def.\ 1.23]{rockafellar1998variational}.
The supremum of all such $\mu$ is the threshold $\mu_h$ of prox-boundedness for $h$.
In particular, if $h$ is bounded below by an affine function, then $\mu_h = \infty$.
For any $\mu \in (0,\mu_h)$, the proximal mapping $\prox_{\mu h}$ 
is locally bounded as well as nonempty- and compact-valued, 
see \cite[Thm~1.25]{rockafellar1998variational}.
The value function of the minimization problem defining the proximal mapping is the 
\emph{Moreau envelope} with parameter $\mu \in (0,\mu_h)$, denoted $\func{h^\mu}{\YY}{\R}$, namely
\begin{equation*}
	h^\mu(z) \coloneqq \inf_{z^\prime} \left\{ h(z^\prime) + \frac{1}{2\mu}\norm{z^\prime-z}^2 \right\} .
\end{equation*}

The \emph{projection} mapping $\ffunc{\proj_\Omega}{\YY}{\YY}$ 
and the \emph{distance} function $\func{\dist_\Omega}{\YY}{\R}$ 
of a nonempty set $\Omega\subseteq\YY$ are defined by
\begin{align*}
	\proj_\Omega(z)
	{}\coloneqq{}&
	\argmin_{z^\prime \in \Omega} \norm{z^\prime - z},
	&
	\dist(z, \Omega)
	{}\coloneqq{}&
	\inf_{z^\prime \in \Omega} \norm{z^\prime - z}.
\end{align*}
The former is a set-valued mapping whenever $\Omega$ is nonconvex, 
whereas the latter is always single-valued.

The following technical lemmas are used later on.

\begin{mybox}
	\begin{lemma}\label{lem:dist_to_domain}
		Let $\func{h}{\YY}{\overline\R}$ be proper, lower semicontinuous, and prox-bounded.
		Let $\dom h$ be closed and fix $\bar z\in\YY$.
		Then
		\[
			\lim\limits_{z^\prime\to\bar z,\,\mu\downto 0}\inf\limits_{z\in\dom h}
				\left\{\mu h(z) + \frac12\norm{z-z^\prime}^2\right\}
			=
			\frac12\dist^2(\bar z,\dom h).
		\]
	\end{lemma}
\end{mybox}
\begin{proof}
	As $\dom h$ is nonempty and closed, we find $\tilde z\in\dom h$ 
	such that $\dist(\bar z,\dom h)=\norm{\tilde z-\bar z}$.
	For every $z^\prime\in\YY$ and $\mu>0$, this gives
	\[
		\inf\limits_{z\in\dom h}\left\{\mu h(z)+\frac12\norm{z-z^\prime}^2\right\}
		\leq
		\mu h(\tilde z)+\frac12\norm{\tilde z-z^\prime}^2,
	\]
	and taking the upper limit, we find
	\begin{equation}\label{eq:upper_semicontinuity_relation}
		\limsup\limits_{z^\prime\to\bar z,\,\mu\downto 0}
		\inf\limits_{z\in\dom h}\left\{\mu h(z)+\frac12\norm{z-z^\prime}^2\right\}
		\leq
		\frac12\norm{\tilde z-\bar z}^2
		=
		\frac12\dist^2(\bar z,\dom h).
	\end{equation}
	
	Next, we define $\func{\psi}{\YY\times[0,\infty)}{\R\cup\{-\infty\}}$ by means of
	\[
		\forall z^\prime\in\YY,\,\forall \mu\in[0,\infty)\colon\quad
		\psi(z^\prime,\mu)
		\coloneqq
		\inf\limits_{z\in\dom h}\left\{\mu h(z)+\frac12\norm{z-z^\prime}^2\right\}.
	\]
	As $h$ is prox-bounded, $\psi$ takes finite values for all sufficiently small $\mu$, 
	and these finite values are attained, see \cite[Thm~1.25]{rockafellar1998variational}.
	Suppose that there are sequences $\{z^k\},\{\bar z^k\}\subseteq\YY$ 
	and $\{\mu_k\}\subseteq[0,\infty)$ with $z^k\to\bar z$ and $\mu_k\to 0$, 
	$\psi(z^k,\mu_k)=\mu_k h(\bar z^k)+\tfrac12\norm{\bar z^k-z^k}^2$, 
	and $\norm{\bar z^k}\to\infty$.
	On the one hand, boundedness of $\{z^k\}$ and $\{\mu_k\}$ 
	gives the existence of a constant $C>0$ such that
	\begin{equation}\label{eq:optimal_value_bounded_from_above}
		\forall k\in\N\colon\quad
		\psi(z^k,\mu_k)\leq \mu_k h(\tilde z)+\frac12\norm{\tilde z-z^k}^2\leq C.
	\end{equation}
	On the other hand, the prox-boundedness of $h$ implies 
	that $h$ is minorized by a quadratic function,
	see \cite[Ex.\ 1.24]{rockafellar1998variational}.
	Hence, there are constants $c_1,c_2,c_3>0$ such that, for sufficiently large $k\in\N$,
	\begin{align*}
		\psi(z^k,\mu_k)
		&\geq
		-\mu_kc_1\norm{\bar z^k}^2-\mu_kc_2\norm{\bar z^k}-\mu_k c_3+\frac12\norm{\bar z^k-z^k}^2\\
		&\geq
		\left(\frac12-\mu_kc_1\right)\norm{\bar z^k}^2-(\mu_kc_2+\norm{z^k})\norm{\bar z^k}-\mu_kc_3.
	\end{align*}
	Boundedness of $\{z^k\}$ and $\mu_k\to 0$ thus yield $\psi(z^k,\mu_k)\to\infty$ since $\norm{\bar z^k}\to\infty$.
	This, however, is a contradiction to \eqref{eq:optimal_value_bounded_from_above}.
	Hence, we can choose a compact set $\mathcal C\subseteq\YY$ such that, for each $\mu\geq 0$ small enough and each
	$z^\prime$ sufficiently close to $\bar z$, we have
	\[
		\psi(z^\prime,\mu)=\inf\limits_{z\in\mathcal C\cap\dom h}
			\left\{\mu h(z)+\frac12\norm{z-z^\prime}^2\right\}.
	\]
	Thus, due to the lower semicontinuity of $h$, we can apply \cite[Thm~4.2.1(1)]{BankGuddatKlatteKummerTammer1983}
	in order to obtain
	\[
		\liminf\limits_{z^\prime\to\bar z,\,\mu\downto 0}\inf\limits_{z\in\dom h}\left\{\mu h(z)+\frac12\norm{z-z^\prime}^2\right\}
		\geq
		\psi(\bar z,0)
		=
		\frac12\dist^2(\bar z,\dom h).
	\]
	Together with \eqref{eq:upper_semicontinuity_relation}, the assertion follows.
\end{proof}

\begin{mybox}
	\begin{lemma}\label{lem:norm_to_zero}
		Let $\func{h}{\YY}{\overline\R}$ be proper, lower semicontinuous, and prox-bounded.
		Fix $\bar z\in\YY$ as well as sequences $\{\mu_k\}\subseteq(0,\infty)$, $\{z^k\}\subseteq\dom h$, and $\{\tilde z^k\}\subseteq\YY$
		such that $\mu_k\downto 0$, $\{\tilde z^k\}$ is bounded, and 
		\begin{equation}\label{eq:some_upper_limit}
			\limsup\limits_{k\to\infty}\left(\mu_k h(z^k)+\frac12\norm{z^k-\tilde z^k}^2\right)
			\leq 0.
		\end{equation}
		Then $\mu_kh(z^k)\to 0$ and $\norm{z^k-\tilde z^k}\to 0$.
	\end{lemma}
\end{mybox}
\begin{proof}
	As in the proof of \cref{lem:dist_to_domain}, we use \cite[Ex.\ 1.24]{rockafellar1998variational} to find 
	constants $c_1,c_2,c_3>0$ such that $h(z)\geq -c_1\norm{z}^2-c_2\norm{z}-c_3$ holds for all $z\in\YY$.
	Thus, we have
	\[
		\limsup\limits_{k\to\infty}
		\left(
			\left(\frac12-\mu_kc_1\right)\norm{z^k}^2-(\mu_kc_2+\norm{\tilde z^k})\norm{z^k}-\mu_kc_3
		\right)
		\leq 0
	\]
	from \eqref{eq:some_upper_limit}.
	By $\mu_k\downto 0$ and boundedness of $\{\tilde z^k\}$, 
	this implies that $\{z^k\}$ is bounded as well.
	Hence, $\{h(z^k)\}$ is bounded below, which gives $\liminf_{k\to\infty}\mu_k h(z^k)\geq 0$.
	Now, \eqref{eq:some_upper_limit} yields the claim.
\end{proof}

\subsection{Variational analysis and generalized differentiation}\label{sec:VA}

\subsubsection*{Tangent and normal cones}

We start by repeating the definition of some cones 
which are well known in variational analysis.
Therefore, we fix some closed set $\Omega\subseteq\YY$ and $\bar z\in\Omega$.
We refer to
\begin{align*}
	T_\Omega(\bar z)
	&\coloneqq
	\left\{
		v\in\YY\,\middle|\,
		\begin{aligned}
			\exists\{t_k\}\subset(0,\infty),\,\exists\{v^k\}\subseteq\YY\colon\,
			t_k\downto 0,\,v^k\to v,\,\bar z+t_kv^k\in\Omega\,\forall k\in\N
		\end{aligned}
	\right\}
\end{align*} 
as the \emph{tangent cone} to $\Omega$ at $\bar z$,
and we point out that it is always a closed cone.
Furthermore, we make use of
\begin{align*}
	\widehat N_\Omega(\bar z)
	&\coloneqq
	\{
		v\in\YY\,|\,
		\forall z\in\Omega\colon\,\innprod{v}{z-\bar z}\leq\oo(\norm{z-\bar z})
	\},\\
	N_\Omega(\bar z)
	&\coloneqq
	\{
		v\in\YY\,|\,
		\begin{aligned}
			\exists\{z^k\}\subseteq\Omega,\,\exists\{v^k\}\subseteq\YY\colon\,
			z^k\to\bar z,\,v^k\to v,\,v^k\in\widehat N_\Omega(z^k)\,\forall k\in\N
		\end{aligned}
	\}
\end{align*}
which are called \emph{regular} (or Fr{\'e}chet) and \emph{limiting} (or Mordukhovich) 
\emph{normal cone} to $\Omega$ at $\bar z$. 
Both of these cones are closed, 
and $\widehat N_\Omega(\bar z)$ is, additionally, convex. 
For a convex set $\Omega$, we have
\[
	\widehat N_\Omega(\bar z)=N_\Omega(\bar z)
	=
	\{v\in\YY\,|\,\forall z\in\Omega\colon\,\innprod{v}{z-\bar z}\leq 0\}.
\]
We would like to point out the polar relation
\begin{equation}\label{eq:polarization_rules}
	\widehat N_\Omega(\bar z)
	=
	T_\Omega(\bar z)^\circ,
\end{equation}
Here, we made use of $A^\circ\coloneqq \{v\in\YY\,|\,\forall z\in A\colon\,\innprod{v}{z}\leq 0\}$, 
the \emph{polar cone} of $A\subseteq\YY$.

\subsubsection*{Subdifferentials and stationarity}

For a lower semicontinuous function $\func{h}{\YY}{\overline\R}$ and $\bar z\in\dom h$,
\begin{align*}
	\widehat{\partial} h(\bar{z})
	&\coloneqq
	\{
		v\in\YY\,|\,
		(v,-1)\in\widehat N_{\epi h}(\bar z,h(\bar z))
	\},\\
	\partial h(\bar z)
	&\coloneqq
	\{
		v\in\YY\,|\,
		(v,-1)\in N_{\epi h}(\bar z,h(\bar z))
	\},\\
	\partial^\infty h(\bar z)
	&\coloneqq
	\{
		v\in\YY\,|\,
		(v,0)\in N_{\epi h}(\bar z,h(\bar z))
	\}
\end{align*}
are referred to as the the \emph{regular} (or Fr{\'e}chet), \emph{limiting} (or Mordukhovich),
and \emph{singular} (or horizon) \emph{subdifferential} of $h$ at $\bar z$.
Whenever $h$ is Lipschitz continuous around $\bar z$, then $\partial^\infty h(\bar z)=\{0\}$.
Let us mention that, among others, 
the subdifferential operators $\widehat\partial$, $\partial$, and $\partial^\infty$ 
are compatible
with respect to smooth additions. 
Indeed, for each continuously differentiable function
$\func{h_0}{\YY}{\R}$, it holds
\[
	\widehat\partial (h_0+h)(\bar z)
	=
	\nabla h_0(\bar z)+\widehat\partial h(\bar z),
	\quad
	\partial(h_0+h)(\bar z)
	=
	\nabla h_0(\bar z)+\partial h(\bar z),
	\quad
	\partial^\infty(h_0+h)(\bar z)
	=
	\partial^\infty h(\bar z).
\]
Whenever $h\coloneqq \indicator_\Omega$, where $\indicator_\Omega$ is the 
(proper and lower semicontinuous)
\emph{indicator function} of 
the nonempty, closed set $\Omega\subseteq\YY$,
vanishing on $\Omega$ and being $\infty$ otherwise,
we have $\dom\indicator_\Omega=\Omega$,
and
\[
	\widehat\partial\indicator_\Omega(\bar z)
	=
	\widehat N_\Omega(\bar z),
	\quad
	\partial\indicator_\Omega(\bar z)
	=
	\partial^\infty\indicator_\Omega(\bar z)
	=
	N_\Omega(\bar z)
\]
for $\bar z\in\Omega$.
The proximal mapping of $\indicator_\Omega$ is the projection $\proj_\Omega$, 
so that $\proj_\Omega$ is locally bounded.

\begin{mybox}
	\begin{lemma}\label{lem:some_helpful_property}
		Let $\func{h}{\YY}{\overline\R}$ be proper, lower semicontinuous,
		and positively homogeneous of degree $2$.
		Then, for each $\bar z\in\dom h$ and $v\in\partial h(\bar z)$,
		$h(\bar z)=\tfrac12 \innprod{v}{\bar z}$ is valid.
	\end{lemma}
\end{mybox}
\begin{proof}
	Let us note that the assertion is trivially true whenever $\bar z=0$ holds,
	so let us assume that $\bar z\neq 0$.
	First, suppose that $v\in\widehat\partial h(\bar z)$.
	Due to \cite[Thm~1.26]{mordukhovich2018}, this yields
	\[
		\liminf\limits_{z\to\bar z,\,z\neq\bar z}
			\frac{h(z)-h(\bar z)-\innprod{v}{z-\bar z}}{\norm{z-\bar z}}
		\geq 0.
	\]
	Considering $z\coloneqq(1\pm t)\bar z$ for $t\downto 0$ in the above limit,
	and exploiting positive homogeneity of degree $2$ of $h$, we find
	$(\pm 2h(\bar z)\mp \innprod{v}{\bar z})/\norm{\bar z}\geq 0$
	which yields $h(\bar z)=\tfrac12 \innprod{v}{\bar z}$.
	To obtain the lemma's assertion in the more general case where
	$v\in\partial h(\bar z)$, we combine the above findings with
	\cite[Thm~1.28]{mordukhovich2018}.
\end{proof}

A point $\bar z\in\dom h$ is said to be \emph{M-stationary} whenever
$0 \in \partial h(\bar z)$ is valid, and this constitutes 
a necessary condition for the local minimality of $\bar z$ for $h$
also known as \emph{Fermat's rule},
see \cite[Thm~10.1]{rockafellar1998variational}.
It should be noted that $0\in\widehat\partial h(\bar z)$ serves as a
(potentially sharper) necessary optimality condition as well.
Given some tolerance $\varepsilon \geq 0$, 
an approximate M-stationarity concept for the minimization of $h$ refers to 
$\dist( 0, \partial h(\bar z) ) \leq \varepsilon$ which we refer to as
\emph{$\varepsilon$-M-stationarity}.
By closedness of $\partial h(\bar z)$, $\varepsilon$-M-stationarity with 
$\varepsilon = 0$ recovers the notion of M-stationarity.

Below, we introduce a stationarity concept
that will be use later to qualify the iterates of our 
implicit AL algorithm,
see \cite[Sec.\ 4.2]{demarchi2024implicit}.
Therefore, let us consider a parametric optimization problem 
with an objective $\func{p}{\XX}{\Rinf}$ and an oracle $\ffunc{\FBO}{\XX}{\YY}$ given by
\begin{equation}
	p(x) \coloneqq \inf_{z\in\YY} P(x,z),
	\qquad
	\FBO(x) \coloneqq \argmin_{z\in\YY} P(x,z)
	\label{eq:fbo}
\end{equation}
for a proper, lower semicontinuous function $\func{P}{\XX \times \YY}{\Rinf}$.
Recalling that the notion of uniform level-boundedness 
corresponds to a parametric extension of level-boundedness, 
see \cite[Def.\ 1.16]{rockafellar1998variational}, 
we suppose that $P$ is level-bounded in $z$ (second argument) 
locally uniformly in $x$ (first argument).
Then, from \cite[Thm~10.13]{rockafellar1998variational} we have 
for every $\bar{x}\in\dom p$ the inclusion
\begin{equation}
	\partial p(\bar{x}) 
	{}\subseteq{} 
	\Upsilon(\bar{x})
	{}\coloneqq{} 
	\bigcup_{\bar{z} \in \FBO(\bar{x})} \left\{ \xi\in\XX \,\middle\vert\, (\xi,0) \in \partial P(\bar{x},\bar{z}) \right\}.
	\label{eq:Upsilon}
\end{equation}
In the setting \eqref{eq:fbo}, because of the parametric nature of $p$, 
the subdifferential mapping $\ffunc{\partial p}{\XX}{\XX}$ is not a simple object in general,
making M-stationarity difficult to check.
Therefore, for the minimization of $p$, one can resort to the concept of \emph{$\Upsilon$-stationarity}, coined in \cite[Def.\ 4.1]{demarchi2024implicit}.

\begin{definition}[$\Upsilon$-stationarity]\label{def:UpsilonStationarity}
	Let $\varepsilon \geq 0$ be fixed and let $\func{P}{\XX \times \YY}{\Rinf}$ 
	be chosen as specified above.
	Define $\func{p}{\XX}{\Rinf}$ and $\ffunc{\Upsilon}{\XX}{\XX}$ 
	as in \eqref{eq:fbo} and \eqref{eq:Upsilon}, respectively.
	Then, relatively to the minimization of $p$, 
	a point $\bar x \in \dom p$ is called \emph{$\varepsilon$-$\Upsilon$-stationary} 
	if $\dist(0, \Upsilon(\bar x)) \leq \varepsilon$.
	In the case $\varepsilon = 0$, such a point $\bar x$ 
	is said to be \emph{$\Upsilon$-stationary}.
\end{definition}

Notice that the inclusion $\bar x \in \dom p$ is implicitly required to have the set $\Upsilon(\bar x)$ nonempty.
Furthermore, in the exact case, namely $\varepsilon = 0$, $\Upsilon$-stationarity of $\bar{x}$ coincides with $0 \in \Upsilon(\bar x)$, 
by closedness of $\Upsilon(\bar x)$ \cite[Thm~10.13]{rockafellar1998variational}.
We shall point out that $\Upsilon$-stationarity provides an intermediate qualification
between (the stronger) M-stationarity for $p$ and
(the weaker) M-stationarity for $P$,
see \cite[Prop.\ 4.1 and 4.2]{demarchi2024implicit} for details.
Finally, it appears from \eqref{eq:Upsilon} that
an $\varepsilon$-$\Upsilon$-stationary point $\bar x\in\XX$ can be
associated to a (possibly nonunique) \emph{certificate} $\bar z \in \FBO(\bar x)$ that satisfies
\begin{equation*}
	\dist(0, \Upsilon(\bar x))
	{}\leq{}
	\min_{\xi\in\XX} \left\{ \norm{\xi} \,\middle\vert\, (\xi,0) \in \partial P(\bar x,\bar z) \right\}
	\leq
	\varepsilon
	.
	\label{eq:UpsilonCertificate}
\end{equation*}
Given such upper bound, the \emph{pair} $(\bar x, \bar z)$ certificates the $\varepsilon$-$\Upsilon$-stationarity of $\bar x$ for $p$.

\subsubsection*{Generalized derivatives of set-valued mapping}

Let us fix a set-valued mapping $\ffunc{\Gamma}{\XX}{\YY}$ and some point
$(\bar x,\bar z)\in\gph\Gamma$.
We refer to the set-valued mappings 
$\ffunc{D\Gamma(\bar x,\bar z)}{\XX}{\YY}$ and
$\ffunc{D^*\Gamma(\bar x,\bar z)}{\YY}{\XX}$, given by
\[
	\begin{aligned}
	D\Gamma(\bar x,\bar z)(u)
	&\coloneqq \{v\in\YY\,|\,(u,v)\in T_{\gph\Gamma}(\bar x,\bar z)\},&\\
	D^*\Gamma(\bar x,\bar z)(y^*)
	&\coloneqq \{x^*\in\XX\,|\,(x^*,-y^*)\in N_{\gph\Gamma}(\bar x,\bar z)\},&
	\end{aligned}
\]
as the \emph{graphical derivative} and the \emph{limiting coderivative} of $\Gamma$ at $(\bar x,\bar z)$.
Subsequently, we review some stability properties of set-valued mappings,
see e.g.\ \cite{BenkoMehlitz2022,rockafellar1998variational}.
We say that $\Gamma$ is \emph{metrically regular} at $(\bar x,\bar z)$ whenever there
are neighborhoods $U\subseteq\XX$ of $\bar x$ and $V\subseteq\YY$ of $\bar z$
as well as a constant $\kappa>0$ such that
\[
	\forall x\in U,\,\forall z\in V\colon\quad
	\dist(x,\Gamma^{-1}(z))\leq\kappa\dist(z,\Gamma(x)).
\]
If just
\[
	\forall x\in U\colon\quad
	\dist(x,\Gamma^{-1}(\bar z))\leq\kappa\dist(\bar z,\Gamma(x))
\]
holds, i.e., if $z\coloneqq\bar z$ can be fixed in the estimate required for metric regularity,
then $\Gamma$ is called \emph{metrically subregular} at $(\bar x,\bar z)$.
Furthermore, $\Gamma$ is said to be \emph{strongly metrically subregular} at $(\bar x,\bar z)$,
whenever there exist a neighborhood $U\subseteq\XX$ of $\bar x$ 
and a constant $\kappa>0$ such that
\[
	\forall x\in U\colon\quad
	\norm{x-\bar x}\leq\kappa\dist(\bar z,\Gamma(x))
\]
is valid.
Recall that strong metric subregularity of $\Gamma$ at $(\bar x,\bar z)$ is equivalent to
$\ker D\Gamma(\bar x,\bar z)=\{0\}$
by the so-called \emph{Levy--Rockafellar criterion}, see \cite[Thm~4E.1]{DontchevRockafellar2014}
and \cite[Prop.\ 4.1]{Levy1996}.
Furthermore, $\Gamma$ is metrically regular at $(\bar x,\bar z)$ if and only if
$\ker D^*\Gamma(\bar x,\bar z)=\{0\}$
by the so-called \emph{Mordukhovich criterion}, see \cite[Thm~9.40]{rockafellar1998variational}.

Let $\func{F}{\XX}{\YY}$ be continuously differentiable and let $\Omega\subseteq\YY$ be closed.
We consider the so-called \emph{feasibility mapping} $\ffunc{\Phi}{\XX}{\YY}$ given by
\begin{equation}\label{eq:feasibility_mapping}
	\Phi(x) \coloneqq F(x)-\Omega.
\end{equation}
We fix some point $\bar x\in\XX$ satisfying $F(\bar x)\in \Omega$, i.e., $(\bar x,0)\in\gph\Phi$.
It is well known that $\Phi$ is metrically regular at $(\bar x,0)$ if and only if
\begin{equation}\label{eq:MRCQ}
	\limnormalcone_\Omega(F(\bar x))\cap\ker F^\prime(\bar x)^\top=\{0\}
\end{equation}
is valid, as we have
\[
	D^*\Phi(\bar x,0)(y)
	=
	\begin{cases}
		\{F^\prime(\bar x)^\top y\}	&	y\in \limnormalcone_\Omega(F(\bar x)),\\
		\emptyset				&	\text{otherwise}
	\end{cases}
\]
from the change-of-coordinates formula in 
\cite[Ex.\ 6.7]{rockafellar1998variational} and the representation
\begin{equation}\label{eq:graph_Phi}
	\gph\Phi=\{(x,w)\in\XX\times\YY\,|\,F(x)-w\in\Omega\}.
\end{equation}
The following lemma, 
which is a direct consequence of \cite[Thm~4B.1]{DontchevRockafellar2014}, 
provides a certain openness-type property of the feasibility mapping
from \eqref{eq:feasibility_mapping} around points of its graph where it is metrically regular. 

\begin{mybox}
\begin{lemma}\label{lem:some_stability_of_MR}
	Fix $(\bar x,0)\in\gph\Phi$ where $\Phi$, the mapping given in
	\eqref{eq:feasibility_mapping}, is metrically regular.
	Then there exist $s>0$ and $\varepsilon>0$ such that
	\[
		\closedball_{s}(0)
		\subseteq
		F^\prime(x)\closedball_1(0)
		-
		\bigl(T_\Omega(z)\cap \closedball_1(0)\bigr)
	\]
	holds true for all $x\in \closedball_\varepsilon(\bar x)$
	and all $z\in \Omega\cap \closedball_\varepsilon(F(\bar x))$.
\end{lemma}
\end{mybox}
\begin{proof}
	Since $\Phi$ is assumed to be metrically regular at $(\bar x,0)$,
	\cite[Thm~4B.1]{DontchevRockafellar2014} yields the existence of
	constants $\delta>0$ and $r>0$ such that
	\[
		\forall(x,w)\in\gph\Phi\cap\closedball_\delta(\bar x,0),\,
		\forall v\in\closedball_1(0),\,
		\exists u\in\closedball_r(0)\colon\quad
		v\in D\Phi(x,w)(u).
	\]
	We apply \cite[Ex.\ 6.7]{rockafellar1998variational} once more
	to the representation \eqref{eq:graph_Phi} in order to find
	\[
		D\Phi(x,w)(u)
		=
		F'(x)u-T_\Omega(F(x)-w)
	\]
	for $(x,w)\in\gph\Phi$.
	Hence, we have
	\[
		\forall(x,w)\in\gph\Phi\cap\closedball_\delta(\bar x,0)\colon\quad
		\closedball_1(0)
		\subseteq
		F'(x)\closedball_r(0)-T_\Omega(F(x)-w).
	\]
	By continuous differentiablity of $F$,
	there is a constant $C>0$ such that $\norm{F'(x)}\leq C$ holds for all
	$x\in\closedball_\delta(\bar x)$.
	With the aid of $\kappa:=\max(r,1+Cr)$, this yields
	\[
		\forall(x,w)\in\gph\Phi\cap\closedball_\delta(\bar x,0)\colon\quad
		\closedball_1(0)
		\subseteq
		F'(x)\closedball_\kappa(0)-\bigl(T_\Omega(F(x)-w)\cap\closedball_{\kappa}(0)\bigr).
	\]
	Let us now choose $\varepsilon\in(0,\delta/2)$ so small such that
	$\norm{F(x)-F(\bar x)}\leq\delta/2$ holds for all $x\in\closedball_\varepsilon(\bar x)$.
	Then, for arbitrary $x\in\closedball_\varepsilon(\bar x)$ and
	$z\in\Omega\cap\closedball_\varepsilon(F(\bar x))$, 
	we can set $w:=F(x)-z$ in order to find 
	$(x,w)\in\gph\Phi\cap\closedball_\delta(\bar x,0)$,
	and the above guarantees
	\[
		\closedball_{1/\kappa}(0)
		\subseteq
		F'(x)\closedball_1(0)-\bigl(T_\Omega(z)\cap\closedball_1(0)\bigr).
	\]
	Choosing $s:=1/\kappa$, the assertion follows.
\end{proof}

\subsubsection*{Subderivatives}

Let us fix a lower semicontinuous function $\func{h}{\YY}{\Rinf}$.
For $\bar z\in\dom h$ and $v\in\YY$, the lower limit
\[
	\d h(\bar z)(v)
	\coloneqq
	\liminf\limits_{t\downto 0,\,v^\prime\to v}
	\frac{h(\bar z+tv^\prime)-h(\bar z)}{t}
\]
is called the \emph{subderivative} of $h$ at $\bar z$ in direction $v$,
and the mapping $v\mapsto \d h(\bar z)(v)$, which, by definition, is lower semicontinuous and
positively homogeneous, is referred to as the subderivative of $h$ at $\bar z$.
We note that $\epi\d h(\bar z)=T_{\epi h}(\bar z,h(\bar z))$, see \cite[Thm~8.2(a)]{rockafellar1998variational}.
Furthermore, for $\bar y\in\YY$,
\begin{equation}\label{eq:def_second_subderivative}
	\d^2 h(\bar z,\bar y)(v)
	\coloneqq
	\liminf\limits_{t\downto 0,\,v^\prime\to v}
	\frac{h(\bar z+tv^\prime)-h(\bar z)-t\innprod{\bar y}{v^\prime}}{\tfrac12 t^2}
\end{equation}
is called the \emph{second subderivative} of $h$ at $\bar z$ for $\bar y$ in direction $v$.
The mapping $v\mapsto\d^2h(\bar z,\bar y)(v)$, which, by definition, is lower semicontinuous and
positively homogeneous of degree $2$, is referred to as the second subderivative of $h$ at $\bar z$ for $\bar y$.
The recent study \cite{BenkoMehlitz2023} presents an overview of calculus rules addressing these variational tools.

\begin{mybox}
\begin{lemma}\label{lem:trivial_second_subderivative}
	Let $\func{h}{\YY}{\Rinf}$ be a lower semicontinuous function, and
	fix $\bar z\in\dom h$ and $\bar y\in\YY$.
	Then we have $\d^2h(\bar z,\bar y)(0)\in\{-\infty,0\}$.
\end{lemma}
\end{mybox}
\begin{proof}
	Observe that $\d^2h(\bar z,\bar y)(0)\leq 0$ holds by definition
	of the second subderivative simply by choosing $v^\prime \coloneqq 0$ in \eqref{eq:def_second_subderivative}.
	Positive homogeneity of degree $2$ of the second subderivative guarantees validity of
	$\d^2 h(\bar z,\bar y)(0)=\alpha^2\d^2 h(\bar z,\bar y)(0)$ for each $\alpha>0$, 
	and this is only possible if $\d^2h(\bar z,\bar y)(0)\in\{-\infty,0\}$.
\end{proof}

For brevity of presentation,
we do not formally introduce the notions of
prox-regularity, subdifferential continuity, and twice epi-differentiability,
which will be used in the next lemma.
Instead,
as the precise meaning of these concepts is not exploited in this paper,
we refer the interested reader to 
\cite[Def.\ 13.27, 13.28, and 13.6(b)]{rockafellar1998variational}
for proper definitions.

\begin{mybox}
\begin{lemma}\label{lem:subderivative_vs_graphical_derivative}
	Let $\func{h}{\YY}{\Rinf}$ be a lower semicontinuous function, and
	fix $\bar z\in\dom h$ and $\bar y\in\partial h(\bar z)$.
	Assume that $h$ is prox-regular, subdifferentially continuous, and twice epi-differentiable at $\bar z$ for $\bar y$.
	Then we have
	\[
		\forall v\in\YY\colon\quad
		D(\partial h)(\bar z,\bar y)(v)
		=
		\frac12\partial \d^2 h(\bar z,\bar y)(v),
	\] 
	and
	\[
		\forall v,w\in\YY\colon\quad
		w\in D(\partial h)(\bar z,\bar y)(v)
		\implies
		\d^2 h(\bar z,\bar y)(v)=\innprod{w}{v}.
	\]
\end{lemma}
\end{mybox}
\begin{proof}
	Recalling that $\d^2 h(\bar z,\bar y)$ is positively homogeneous of degree $2$,
	the second property follows with the aid of
	\cref{lem:some_helpful_property} from the first one,
	which is taken from \cite[Thm~13.40]{rockafellar1998variational}.
\end{proof}

\section{Fundamentals of composite optimization}
\label{sec:compositeOptimization}

We now move our attention to \eqref{eq:P} and discuss relevant optimality and stationarity notions.
Furthermore, we investigate local characterizations using second-order tools, regularity concepts, and error bounds.

\subsection{Stationarity concepts and Lagrangian-type functions}
\label{sec:concepts}

Before dealing with optimality conditions,
we consider some Lagrangian terminology and notions
useful for first-order analysis.
Including an auxiliary variable $z\in\YY$,
we can lift \eqref{eq:P} as \eqref{eq:ALMz:P}
involving merely equality constraints but no (nontrivial) compositions.
Introducing a Lagrange multiplier $y\in\YY$ for the constraints,
we define 
a Lagrangian-type function
$\func{\LLslack}{\XX\times\YY\times\YY}{\Rinf}$ 
associated with \eqref{eq:ALMz:P} by means of
\begin{equation}
	\label{eq:ALMz:Lagrangian}
	\LLslack(x,z,y)
	{}\coloneqq{}
	f(x) + g(z) 
	+ 
	\innprod{y}{c(x) - z}.
\end{equation}
Focusing on those terms of $\LLslack$ depending on $x$,
we call the function $\func{\LL}{\XX\times\YY}{\R}$ given by
\begin{equation}
	\label{eq:Lagrangian}
	\LL(x,y)
	\coloneqq
	f(x) + \innprod{y}{c(x)}
\end{equation}
the \emph{Lagrangian} function of \eqref{eq:P}.
Then, acting as a precursor of $\LL$,
we refer to $\LLslack$ as the \emph{pre-Lagrangian} function of \eqref{eq:P}.
These objects are tightly related to the so-called M-stationarity conditions 
of both problems \eqref{eq:ALMz:P} and \eqref{eq:P}, see \cref{def:Mstationary} below.
In fact,
these first-order optimality conditions can be expressed in Lagrangian form as
\begin{align*}
	0
	{}\in{}&
	\partial_x \LLslack(\bar{x}, \bar{z}, \bar{y}) ,&
	0
	{}\in{}&
	\partial_z \LLslack(\bar{x}, \bar{z}, \bar{y}) ,&
	0
	{}\in{}&
	\partial_y\LLslack(\bar{x}, \bar{z}, \bar{y})
\end{align*}
or more explicitly as
\begin{align}\label{eq:ALMz:Mstationary}
	0
	{}={}&
	\nabla_x \LL(\bar{x},\bar{y}) ,&
	\bar{y}
	{}\in{}&
	\partial g( \bar{z} ) ,&
	0
	{}={}&
	c(\bar{x})-\bar{z} .
\end{align}
Equivalently, albeit omitting the auxiliary variable $\bar{z} = c(\bar{x})$, these read
\begin{subequations}\label{eq:Mstationary}
	\begin{align}
		0 {}={}& \nabla_x \LL(\bar{x},\bar{y}) ,
		\label{eq:Mstationary:x}\\
		\bar{y} {}\in{}& \partial g( c(\bar{x}) ) .
		\label{eq:Mstationary:y}
	\end{align}
\end{subequations}
Notice that \eqref{eq:Mstationary:y}
implicitly requires the feasibility of $\bar{x}$ for \eqref{eq:P}, 
namely $c(\bar{x}) \in \dom g$,
for otherwise
the subdifferential $\partial g( c(\bar{x}) )$
is empty.

Interpreting \eqref{eq:P} as an \emph{unconstrained} problem,
first-order necessary optimality conditions
using the notion of M-stationarity
pertain a point $\bar{x}\in\XX$ such that
$0\in\partial\varphi(\bar x)$.
We now aim to rewrite this condition in terms of initial
problem data, i.e., first-order (generalized) derivatives of $f$, $c$, and $g$.
Exploiting compatibility of the limiting subdifferential with respect to smooth additions,
we find $0\in\nabla f(x)+\partial(g\circ c)(\bar x)$. 
It has been recognized, e.g.\ in \cite[Sec.\ 3.2]{IoffeOutrata2008}, 
that metric subregularity of the set-valued mapping
$\ffunc{\Xi}{\R^n\times\R}{\R^m\times\R}$ given by
\begin{equation}\label{eq:svm_composition}
	\Xi(x,\alpha)
	\coloneqq
	(c(x),\alpha)-\epi g
\end{equation}
is enough to guarantee that a subdifferential chain rule can be used to
approximate the limiting subdifferential of $g\circ c$ from above in terms
of the subdifferential of $g$ and the derivative of $c$.
More precisely, if $\Xi$ is metrically subregular at $((\bar x,g(c(\bar x))),(0,0))$,
then $\partial(g\circ c)(\bar x)\subseteq c^\prime(\bar x)^\top\partial g(c(\bar x))$.
From a Lagrangian perspective,
this gives the existence of some Lagrange multiplier 
$\bar{y}\in\YY$ such that the stationarity conditions \eqref{eq:Mstationary} hold.
This stationarity characterization resembles, at least in spirit,
the so called Karush--Kuhn--Tucker (or KKT) conditions in
nonlinear programming, see e.g.\ \cite{bertsekas1996constrained,birgin2014practical}.

These considerations lead to the following definition,
which uses in accordance with \eqref{eq:Mstationary} the Lagrangian function
$\func{\LL}{\XX\times\YY}{\R}$
associated with \eqref{eq:P},
given by \eqref{eq:Lagrangian}.

\begin{definition}[M-stationarity]\label{def:Mstationary}
	Relative to \eqref{eq:P}, a point $\bar x \in \XX$ is called \emph{M-stationary} if there exists a multiplier $\bar{y} \in \YY$ such that \eqref{eq:Mstationary} holds.
	Let $\Lambda(\bar{x})$ denote the set of multipliers $\bar{y}\in\YY$ such that the M-stationarity conditions 
	\eqref{eq:Mstationary} are satisfied by $(\bar{x},\bar{y})$.
\end{definition}

As a reminder of the possible gap highlighted above,
where metric subregularity of $\Xi$ is invoked to formulate
first-order optimality conditions in Lagrangian terms,
the notion given in \cref{def:Mstationary} could be referred to as \emph{KKT}-stationarity, as in \cite{demarchi2024implicit}.
For simplicity, we stick to the nomenclature of \emph{M}-stationarity.

Following the nomenclature in \cite{bolte2018nonconvex,sabach2019lagrangian},
$\LLslack$ would be referred to as the Lagrangian of \eqref{eq:ALMz:P},
as a standalone problem,
and not only as the pre-Lagrangian in view of \eqref{eq:P}.
In fact, the definition of $\LLslack$ in \eqref{eq:ALMz:Lagrangian} complies with the classical concept of Lagrangian function
for equality-constrained optimization problems, such as \eqref{eq:ALMz:P},
and reflects the (nonsmooth, extended real-valued) objective 
$(x,z) \mapsto f(x) + g(z)$ of \eqref{eq:ALMz:P}
and its equality constraints $c(x) - z = 0$.
However, containing (primal) nonsmooth terms, $\LLslack$ is not differentiable.
The object $\LL$ defined in \eqref{eq:Lagrangian}
corresponds to the \emph{ordinary} Lagrangian function of \eqref{eq:P}
as described in \cite{rockafellar2022convergence},
and this is consistent with several other papers which exploit the 
variational analysis approach to composite optimization,
see e.g.\ \cite{BenkoGfrererYeZhangZhou2022,BenkoMehlitz2023,HangMordukhovichSarabi2022,HangSarabi2021,MohammadiMordukhovichSarabi2022,MordukhovichSarabi2018,MordukhovichSarabi2019}
and, particularly, the setting of standard nonlinear programming,
see \cref{ex:NLP} below.

Above, we derived the M-stationarity conditions of \eqref{eq:P} at some feasible point $\bar x\in\XX$
by using the chain rule for the limiting subdifferential which, in general, 
requires a qualification condition like metric subregularity of $\Xi$ from \eqref{eq:svm_composition}
at $((\bar x,g(c(\bar x))),(0,0))$.
Note that \eqref{eq:MRCQ} reduces to
\begin{equation}\label{eq:metric_regularity_CQ}
	\partial^\infty g(c(\bar x))\cap\ker c^\prime(\bar x)^\top
	=
	\{0\}
\end{equation}
when applied to $\Xi$ at the given reference point, and the latter is
equivalent to the mapping $\Xi$ being metrically regular
at $((\bar x,g(c(\bar x))),(0,0))$, which also extends to
a neighborhood of the point of interest.
Thus, \eqref{eq:metric_regularity_CQ} is
sufficient for the subregularity requirement stated earlier.
Clearly, \eqref{eq:metric_regularity_CQ} is valid whenever
$g$ is locally Lipschitz continuous at $c(\bar x)$ 
or if $c^\prime(\bar x)$ has full row rank.
As we know that $0\in\partial\varphi(\bar x)$ provides a
necessary optimality condition for the local optimality of
$\bar x$, the M-stationarity conditions from \cref{def:Mstationary}
do so as well in the presence of a suitable CQ 
like \eqref{eq:metric_regularity_CQ} as outlined above.

\paragraph{Augmented Lagrangian}

We shall introduce \emph{augmented} Lagrangian functions,
which not only offer the basic component for AL methods,
but also closely relate to first-order optimality concepts.
An AL function for \eqref{eq:P}
can be obtained in two steps:
augmenting the pre-Lagrangian $\LLslack$ with a penalty term,
and then marginalizing over the auxiliary variables.

For some penalty parameter $\mu>0$,
the AL function
$\func{\LLslack_\mu}{\XX \times \YY \times \YY}{\Rinf}$
associated to \eqref{eq:ALMz:P}
entails the sum of the pre-Lagrangian $\LLslack$ and a quadratic penalty for the constraint violation,
scaled by $\mu$.
This leads to the definition of $\LLslack_\mu$ as
\begin{equation}\label{eq:ALMz:L}
	\LLslack_\mu(x,z,y)
	{}\coloneqq{}
	\LLslack(x,z,y) + \frac{1}{2 \mu} \norm{c(x) - z}^2 . \nonumber
\end{equation}
Then, since \eqref{eq:ALMz:P} involves the minimization over both original and auxiliary variables,
whereas the latter ones do not appear in the original problem \eqref{eq:P},
we consider the marginalization of $\LLslack$ over $z$,
which yields the
\emph{augmented Lagrangian} function
$\func{\LL_\mu}{\XX \times \YY}{\R}$ associated to \eqref{eq:P}:
\begin{align}\label{eq:L}
	\LL_\mu(x,y)
	{}\coloneqq{}
	\inf_z \LLslack_\mu(x,z,y)
	{}={}&
	f(x)
	+ \inf_z \left\{ g(z) + \frac{1}{2 \mu} \norm{c(x) + \mu y - z}^2 \right\} - \frac{\mu}{2} \norm{y}^2 \\
	{}={}&
	f(x) + g^\mu( c(x) + \mu y ) - \frac{\mu}{2} \norm{y}^2 . \nonumber
\end{align}
Notice that the minimization over $z$ is well-defined only for sufficiently small penalty parameters,
relative to the prox-boundedness threshold of $g$,
in particular $\mu \in (0, \mu_g)$.
Moreover, we highlight that the Moreau envelope $\func{g^\mu}{\YY}{\R}$ of $g$ is real-valued 
and strictly continuous \cite[Ex.\ 10.32]{rockafellar1998variational}, 
but not continuously differentiable in general, as the proximal mapping of $g$ is possibly set-valued.

With the AL tools at hand,
one can readily recover the M-stationarity conditions \eqref{eq:Mstationary} for \eqref{eq:P}.
Through the augmented pre-Lagrangian function $\LLslack_\mu$ of \eqref{eq:P},
the first-order optimality conditions in the form of \eqref{eq:ALMz:Mstationary}
can be written, for any $\mu > 0$, as
$0 \in \partial_{x} \LLslack_\mu(\bar{x}, c(\bar{x}), \bar{y})$, $0\in\partial_z\LLslack_\mu(\bar x,c(\bar x),\bar y)$,
and $0\in\partial_y\LLslack_\mu(\bar x,c(\bar x),\bar y)$,
which hold if and only if $(\bar{x},\bar{y}) \in \XX \times \YY$ satisfies \eqref{eq:Mstationary}.

The following lemma, which is inspired by \cite[Lem.\ 3.1]{BoergensKanzowSteck2019}, will come in handy later on.
\begin{mybox}
	\begin{lemma}\label{lem:upper_bound_AL}
		If $x\in\XX$ is feasible for \eqref{eq:P}, then $\LL_\mu(x,y) \leq \varphi(x)$ for all $\mu > 0$ and $y \in \YY$.
	\end{lemma}
\end{mybox}
\begin{proof}
	By feasibility of $x$, we have $c(x) \in \dom g$.
	Then, for all $\mu > 0$ and $y\in\YY$,
	\begin{equation*}
		g^\mu(c(x)+\mu y )
		{}={}
		\inf_z \left\{ g(z) + \frac{1}{2\mu} \norm{c(x)+\mu y-z}^2 \right\}
		{}\leq{}
		g(c(x)) + \frac{1}{2\mu} \norm{\mu y}^2
		{}={}
		g(c(x)) + \frac{\mu}{2} \norm{y}^2
	\end{equation*}
	by selecting $z = c(x) \in \dom g$.
	Hence,
	$\LL_\mu(x,y)
	=
	f(x) + g^\mu(c(x)+\mu y ) - \tfrac{\mu}{2} \norm{y}^2
	\leq
	f(x) + g(c(x))
	=
	\varphi(x)$,
	concluding the proof.
\end{proof}

In view of the AL subproblems arising in \cref{sec:ALM} below,
the subsequent remark considers the notion of $\Upsilon$-stationarity,
discussed already in \cref{sec:VA},
to the AL function $\LL_\mu$.

\begin{remark}\label{rem:Upsilon_stationarity_AL}
	Motivated by the minimization of the AL function $\LL_\mu(\cdot,\hat y)$
	where $\hat y\in\YY$ is fixed,
	we are interested in pairs $(\bar x,\bar z)\in\XX\times\YY$ which certificate $\varepsilon$-$\Upsilon$-stationarity
	of $\bar x$ for $\LL_\mu(\cdot,\hat y)$,
	for some given $\varepsilon\geq 0$.
	A simple calculation reveals that
	\[
		\Upsilon(\bar x)
		=
		\bigcup\limits_{\bar z\in\prox_{\mu g}(c(\bar x)+\mu \hat y)}
		\left\{
		\nabla_x\LL(\bar x,\bar y) \,\middle|\, \bar y\coloneqq\hat y + \frac{c(\bar x)-\bar z}{\mu} \in\partial g(\bar z)
		\right\}
	\]
	holds in this situation.
	Clearly, $\bar z\in\prox_{\mu g}(c(\bar x)+\mu \hat y)$ always gives $\hat y+(c(\bar x)-\bar z)/\mu\in\partial g(\bar z)$ by Fermat's rule
	(the converse is true for convex $g$),
	so $\varepsilon$-$\Upsilon$-stationarity boils down to the existence of some
	$\bar z\in\prox_{\mu g}(c(\bar x)+\mu \hat y)$
	such that $\norm{\nabla_x\LL(\bar x,\bar y)}\leq\varepsilon$ where $\bar y \coloneqq \hat y+(c(\bar x)-\bar z)/\mu$.
	Note that this implicitly demands $\mu\in(0,\mu_g)$.
	Obviously, for arbitrary $\mu>0$ and any pair $(\bar x,\bar z)$ certificating $\Upsilon$-stationarity (where $\varepsilon\coloneqq 0$) 
	of $\bar x$ for $\LL_\mu(\cdot,\hat y)$ in the above sense such that $\bar z=c(\bar x)$ holds, $\bar x$ is also M-stationary.
	The converse is true whenever $g$ is a convex function, and,
	in this case, the proximal mapping is
	single-valued.
\end{remark}

Some of the concepts addressed in this section are visualized in the following example
in terms of standard nonlinear programming.

\begin{example}\label{ex:NLP}
	Nonlinear programming can be cast in the form \eqref{eq:P} via many reformulations.
	Let us consider the setting
	\begin{equation}
		\tag{NLP}\label{eq:NLP}
		\minimize_x
		{}\quad{}
		f(x)
		{}\qquad{}
		\stt
		{}\quad{}
		c(x) \in C
	\end{equation}
	with $g \equiv \indicator_C$ being the indicator of a nonempty, closed, convex set $C \coloneqq [c_l, c_u]$.
	Allowing entries of $c_l$ and $c_u$ to take infinite values,
	namely $c_l \in (\R\cup\{-\infty\})^m$ and $c_u \in (\R\cup \{\infty\})^m$,
	the model includes equalities, inequalities, and bounds in a compact form,
	and the constraint set $C$ is convex polyhedral.
	The pre-Lagrangian for \eqref{eq:NLP} with auxiliary variable $z\in\YY$ and multiplier $y\in\YY$ reads
	\begin{equation*}
		\LLslack(x,z,y)
		{}={}
		f(x) + \indicator_C(z) + \innprod{y}{c(x) - z} .
	\end{equation*}
	The 
	M-stationarity conditions of \eqref{eq:NLP}
	can be expressed as
	\begin{align*}
		\nabla f(\bar{x}) + c^\prime(\bar{x})^\top \bar{y} {}={}& 0 ,
		&
		\bar{y} {}\in{}& \normalcone_C( c(\bar{x}) ) ,
	\end{align*}
	where the inclusion coincides with the classical complementarity conditions
	and imposes the feasibility condition $c(\bar{x}) \in C$ as well.
	The Lagrangian $\LL$ for \eqref{eq:NLP} is
	$
	\LL(x,y)
	{}={}
	f(x) + \innprod{y}{c(x)}
	$
	and the AL $\LL_\mu$ is given by
	\begin{equation*}
		\LL_\mu(x,y)
		{}={}
		f(x) + \frac{1}{2\mu} \dist_C^2( c(x) + \mu y ) - \frac{\mu}{2} \norm{y}^2
		,
	\end{equation*}
	recovering all classical quantities.
	As $C$ is convex in \eqref{eq:NLP}, the squared distance term in the AL function is continuously differentiable,
	see e.g.\ \cite[Cor.\ 12.30]{BauschkeCombettes2011}.
\end{example}

\begin{remark}\label{rem:generalized_Lagrangian}
	Yet another way to the definition of a Lagrangian-type function in composite optimization
	with \emph{convex} function $g$ has been promoted by Rockafellar in his recent papers
	\cite{rockafellar2022augmented,rockafellar2022convergence} where he introduces
	the so-called \emph{generalized} Lagrangian of \eqref{eq:P}
	by marginalization of the pre-Lagrangian $\LLslack$ given in \eqref{eq:ALMz:Lagrangian}.
	The marginalization step enters here because
	\eqref{eq:ALMz:P} involves the minimization over both original and auxiliary variables,
	whereas the latter ones do not appear in \eqref{eq:P}.
	Marginalization of $\LLslack$ over $z$
	results in the generalized Lagrangian $\func{\ell}{\XX\times\YY}{\R\cup\{-\infty\}}$ given by
	\begin{equation*}
		\ell(x,y)
		{}\coloneqq{}
		\inf_z
		\LLslack(x,z,y) 
		{}={}
		f(x) 
		+ 
		\innprod{y}{c(x)} 
		+ 
		\inf_z \{ g(z) - \innprod{y}{z} \} 
		{}={}
		\LL(x,y) - g^\conj(y) , 
	\end{equation*}
	where $\LL$ is the Lagrangian defined in \eqref{eq:Lagrangian}.
	Clearly, it is $\nabla_x \LL = \nabla_x \ell$.
	One could hope that
	the generalized Lagrangian $\ell$ encapsulates all information
	needed to state the M-stationarity conditions \eqref{eq:Mstationary} for \eqref{eq:P},
	exemplary
	as $0 \in \partial_x\ell(\bar{x},\bar{y}),0 \in \partial_y(-\ell)(\bar{x},\bar{y})$.
	The negative sign appearing for the multipliers relates to the (generalized) saddle-point
	properties of the (generalized) Lagrangian.
	Expanding terms, this gives
	\begin{subequations}\label{eq:stat_via_gen_Lag}
		\begin{align}
		\label{eq:stat_via_gen_Lag_x}
		0
		{}={}&
		\nabla_x \LL(\bar{x},\bar{y}), \\
		\label{eq:stat_via_gen_Lag_y}
		c(\bar{x}) {}\in{}& \partial g^\conj(\bar{y}).
		\end{align}
	\end{subequations}
	If $g$ is convex, \eqref{eq:stat_via_gen_Lag_y} is equivalent to \eqref{eq:Mstationary:y},
	see \cite[Prop.\ 11.3]{rockafellar1998variational},
	so that M-stationarity can be fully characterized via the derivatives of the generalized Lagrangian.
	Whenever $g$ is a nonconvex function, however, this reasoning is no longer possible.
	Under additional assumptions on $g$ (and $g^\conj$), one may apply the convex hull property,
	see e.g.\ \cite[formula (2.7)]{DempeDuttaMordukhovich2007}, and a marginal function rule,
	see e.g.\ \cite[Thm~5.1]{BenkoMehlitz2022} or \cite[Thm~10.13]{rockafellar1998variational}, 
	in order to find
	\[
		\partial g^\conj(y)
		\subseteq
		-\conv\partial(-g^\conj)(y)
		\subseteq
		-\conv\{-z\in\YY\,|\,y\in\partial g(z)\}
		=
		\conv(\partial g)^{-1}(y).
	\]
	Hence, whenever $(\partial g)^{-1}(y)$ is convex, \eqref{eq:stat_via_gen_Lag_y} yields
	\eqref{eq:Mstationary:y} if the aforementioned calculus rules apply.
	Consequently, under additional assumptions, \eqref{eq:stat_via_gen_Lag}
	implies the M-stationarity conditions \eqref{eq:Mstationary} even for nonconvex $g$.
	However,
	the converse implication is likely to fail even in very simple situations,
	as illustrated in the subsequently stated \cref{ex:conjugates_for_nonconvex_g_messy}.
\end{remark}

\begin{example}\label{ex:conjugates_for_nonconvex_g_messy}
	We investigate the model problem
	\begin{equation}\label{eq:sparse_programming}
	\minimize_{x}
	{}\quad{}
	f(x) + \norm{c(x)}_0,
	\end{equation}
	where $g$ plays the role of the $\ell_0$-quasi-norm $\func{\norm{\cdot}_0}{\R^m}{\R}$, which simply counts
	the nonzero entries of the argument vector. 
	Clearly, $\norm{\cdot}_0$ is a merely lower semicontinuous function and is not convex.
	For some point $\bar x\in\R^n$, we will exploit the index sets
	\[
		I^0(\bar x)\coloneqq \{i\in\{1,\ldots,m\}\,|\,c_i(\bar x)=0\},
		\quad
		I^{\pm}(\bar x)\coloneqq \{1,\ldots,m\}\setminus I^0(\bar x).
	\]
	One can easily check that
	\[
		\partial\norm{\cdot}_0(c(\bar x))
		=
		\{y\in\R^m\,|\,\forall i\in I^{\pm}(\bar x)\colon\,y_i=0\}
	\]	
	holds true. Hence, $\bar x$ is M-stationary for \eqref{eq:sparse_programming}
	if and only if there is some $\bar y\in\YY$ such that \eqref{eq:Mstationary:x} is valid and,
	for all $i\in I^\pm(\bar x)$, it is $\bar y_i=0$.
	A simple calculation reveals that $\norm{\cdot}_0^\conj	= \indicator_{\{0\}}$,
	which is why condition \eqref{eq:stat_via_gen_Lag_y} reduces to $\bar y=0$.
	This is a much stronger requirement on the multiplier
	than the one demanded by M-stationarity.
\end{example}

\subsection{Second-order sufficient optimality conditions}
\label{sec:SOSC:P}

In this subsection, we briefly review the second-order sufficient optimality condition
for \eqref{eq:P} which has been derived in \cite[Sec.\ 6]{BenkoMehlitz2023}.

Let us fix a feasible point $\bar x\in\R^n$ of \eqref{eq:P}, and define the \emph{critical cone}
of \eqref{eq:P} at $\bar x$ by means of
\[
	\mathcal C(\bar x)
	\coloneqq
	\{u\in\R^n\,|\,f^\prime(\bar x)u+\d g(c(\bar x))(c^\prime(\bar x)u)\leq 0\}.
\]
Furthermore, for each $u\in\mathcal C(\bar x)$, we make use of the 
\emph{directional multiplier set} $\Lambda(\bar x,u)$ given by
\[
	\Lambda(\bar x,u)
	\coloneqq
	\{y\in\R^m\,|\,
		\nabla_x\LL(\bar x,y)=0,\,
		\d g(c(\bar x))(c'(\bar x)u)=\innprod{y}{c'(\bar x)u},\,
		\d^2g(c(\bar x),y)(c'(\bar x)u)>-\infty
	\}.
\]
Let us mention that this definition
can be stated equivalently in terms of the so-called directional proximal
subdifferential of $g$, see \cite[Sec.\ 3.2]{BenkoMehlitz2023} for details.

\begin{definition}[Second-order sufficient condition]\label{def:SOSC}
	For a feasible point $\bar x\in\R^n$ of \eqref{eq:P}, we say that the
	\emph{Second-Order Sufficient Condition} (SOSC for brevity) is valid, 
	whenever 
	\begin{equation*}
		\forall u\in\mathcal C(\bar x)\setminus\{0\},\,
		\exists y\in\Lambda(\bar x,u)\colon\quad
		\nabla^2_{xx}\LL(\bar x,y)[u,u]
		+
		\d^2g(c(\bar x),y)(c^\prime(\bar x)u)
		>
		0.
	\end{equation*}
\end{definition}

Let us fix a feasible point $\bar x\in\XX$ of \eqref{eq:P} where SOSC is valid. 
To avoid trivial situations, we assume that $\mathcal C(\bar x)$ contains
a non-vanishing direction $u\in\XX$. 
Clearly, SOSC requires that $\Lambda(\bar x,u)$ is nonempty, 
and since we have $\Lambda(\bar x,u)\subseteq\Lambda(\bar x)$ from
\cite[Prop.\ 2.9]{BenkoGfrererYeZhangZhou2022} and \cite[Cor.\ 3.14]{BenkoMehlitz2023},
$\Lambda(\bar x)$ is nonempty as well, i.e., $\bar x$ is M-stationary.
Thus, checking validity of SOSC is only reasonable at M-stationary points.

The following result is a consequence of \cite[Thm~6.1]{BenkoMehlitz2023}
and the associated discussions.
\begin{mybox}
	\begin{proposition}\label{prop:SOSC_P}
		Let $\bar x\in\R^n$ be a feasible point of \eqref{eq:P}
		where SOSC is valid.
		Then there exist constants $\varepsilon>0$ and $\beta>0$ 
		such that
		the second-order growth condition
		\begin{equation}\label{eq:second_order_growth_P}
			\forall x\in\closedball_\varepsilon(\bar x)\colon\quad
			\varphi(x)-\varphi(\bar x)\geq\frac{\beta}{2}\norm{x-\bar x}^2
		\end{equation}
		is valid. 
		Particularly, $\bar x$ is a strict local minimizer of \eqref{eq:P}.
	\end{proposition}
\end{mybox}

Since the composite optimization problem \eqref{eq:P} can
be recast as the constrained problem
\[
	\minimize_{x,\tau}
	{}\quad{}
	f(x)+\tau
	{}\quad{}
	\stt
	{}\quad{}
	(c(x),\tau)\in\epi g,
\]
\cref{prop:SOSC_P} is also a consequence of \cite[Thm~3.3]{BenkoGfrererYeZhangZhou2022}
when taking into account the inequality
\[
	\d^2\indicator_{\epi g}((c(\bar x),g(c(\bar x))),(y,-1))(c'(\bar x)u,\upsilon)
	\geq
	\d^2 g(c(\bar x),y)(c'(\bar x)u)
\]
which holds for each $u\in\XX$ and $\upsilon\in\R$, see \cite[Prop.\ 3.13]{BenkoMehlitz2023}.

The following corollary provides a sufficient condition for SOSC which will be of
interest later on.
\begin{mybox}	
	\begin{corollary}\label{cor:SSOSC}
		Let $\bar x\in\R^n$ be an M-stationary point for \eqref{eq:P}, and fix $\bar y\in\Lambda(\bar x)$.
		Furthermore, let the condition
		\begin{equation}\label{eq:SSOSC}
			\forall u\in\mathcal C(\bar x)\setminus\{0\}\colon\quad
			\nabla^2_{xx}\LL(\bar x,\bar y)[u,u]
			+
			\d^2g(c(\bar x),\bar y)(c^\prime(\bar x)u)
			>
			0
		\end{equation}
		hold.
		Then SOSC and, thus, 
		the second-order growth condition \eqref{eq:second_order_growth_P} are valid.
	\end{corollary}
\end{mybox}
\begin{proof}
	We will show that $\bar y\in\Lambda(\bar x,u)$ holds for each $u\in\mathcal C(\bar x)\setminus\{0\}$.
	Then it is clear that \eqref{eq:SSOSC} implies validity of SOSC, 
	and the final assertion is just a consequence of \cref{prop:SOSC_P}.
	
	Fix $u\in\mathcal C(\bar x)\setminus\{0\}$.
	Then \eqref{eq:SSOSC} obviously implies $\d^2g(c(\bar x),\bar y)(c'(\bar x)u)>-\infty$,
	and \cite[formula (5)]{BenkoMehlitz2023} immediately yields
	$\d g(c(\bar x))(c'(\bar x)u)\geq \innprod{\bar y}{c'(\bar x)u}$.
	By definition of the critical cone and $\nabla_x\LL(\bar x,\bar y)=0$, we also have
	\[
		\d g(c(\bar x))(c'(\bar x)u)
		\leq
		-f'(\bar x)u
		=
		\innprod{\bar y}{c'(\bar x)u}.
	\]
	Thus, $\bar y\in\Lambda(\bar x,u)$ follows.
\end{proof}

\subsection{Error bounds}\label{sec:ErrorBounds}

Here, we aim to establish a connection between the second-order sufficient conditions
from \cref{def:SOSC} and an error bound property.
Relating to stability properties and involving the distance to the primal-dual solution set,
error bounds are an essential ingredient for deriving rates of local convergence for
numerical methods addressing \eqref{eq:P}.
In order to quantify the violation of the M-stationarity conditions from \cref{def:Mstationary}, 
it is (almost) natural to define the residual mapping
\begin{equation}\label{eq:ErrorBound:error}
	\Theta(x,z,y)
	\coloneqq
	\norm{\nabla_x \LL(x,y)}
	+
	\norm{ c(x) - z }
	+
	\dist( y, \partial g(z) ).
\end{equation}
Clearly, the M-stationarity conditions \eqref{eq:Mstationary} for some $(\bar x,\bar y)\in\XX\times\YY$ 
are equivalent to $\Theta(\bar x,c(\bar x),\bar y) = 0$.
We shall see now that, under certain assumptions, $\Theta$ allows us to quantify not only the violation of \eqref{eq:Mstationary}, 
but also the distance to the (primal-dual) solution set.

The proof of the following proposition, 
which provides the foundations of our analysis in this section, 
relates the error bound property of our interest with the strong metric subregularity 
of a certain set-valued mapping.
Moreover, the latter characterization is quantifiable via a condition which
can be stated in terms of initial problem data,
thanks to the Levy--Rockafellar criterion.

\begin{mybox}
	\begin{proposition}\label{lem:strong_ms_yields_error_bound}
		Let $\bar x\in\R^n$ be an M-stationary point of \eqref{eq:P} and pick $\bar y\in\Lambda(\bar x)$.
		Assume that the qualification condition
		\begin{equation}\label{eq:strong_metric_subregularity}
			0=\nabla^2_{xx}\LL(\bar x,\bar y)u+c^\prime(\bar x)^\top \eta
			,\,
			\eta\in D(\partial g)(c(\bar x),\bar y)(c^\prime(\bar x)u)
			\implies
			u=0,\,\eta=0
		\end{equation}
		is valid.
		Then there are a constant $\varrho_{\mathrm{u}}>0$ and
		a neighborhood $U$ of $(\bar x,c(\bar x),\bar y)$ such that,
		for each $(x,z,y)\in U\cap(\R^n\times\dom g\times\R^m)$, we have
		the upper estimate
		\begin{equation}\label{eq:error_bound}
			\norm{x-\bar x}
			+
			\norm{z-c(\bar x)}
			+
			\norm{y-\bar y}
			\leq
			\varrho_{\mathrm{u}}\,\Theta(x,z,y).
		\end{equation}
	\end{proposition}
\end{mybox}
\begin{proof}
	We define a set-valued mapping $\ffunc{G}{\R^n\times\R^m\times\R^m}{\R^n\times\R^m\times\R^m}$ by means of
	\begin{equation*}
		G(x,z,y)
		\coloneqq
		(\nabla_x\LL(x,y),c(x)-z,y)-\{0\}\times\{0\}\times \partial g(z).
	\end{equation*}
	By continuous differentiability of the single-valued part of $G$, one can easily check,
	e.g., by means of the change-of-coordinates formula for tangents from \cite[Ex.\ 6.7]{rockafellar1998variational}, that
	\begin{multline*}
		DG((\bar x,c(\bar x),\bar y),(0,0,0))(u,v,\eta)
		\\
		=
		(\nabla^2_{xx}\LL(\bar x,\bar y)u+c^\prime(\bar x)^\top\eta,c^\prime(\bar x)u-v,\eta)
		-
		\{0\}\times\{0\}\times D(\partial g)(\bar y,c(\bar x))(v)
	\end{multline*}
	holds. Thus,
	\eqref{eq:strong_metric_subregularity} is equivalent to
	$
		\ker DG((\bar x,c(\bar x),\bar y),(0,0,0))=\{(0,0,0)\}
	$.
	By the Levy--Rockafellar criterion, $G$ is strongly metrically subregular at $(\bar x,c(\bar x),\bar y)$,
	and the latter is equivalent to the desired error bound condition.
\end{proof}

Note that the proof of \cref{lem:strong_ms_yields_error_bound} actually shows that \eqref{eq:strong_metric_subregularity}
is \emph{equivalent} to the local validity of the error bound property \eqref{eq:error_bound}.

\begin{mybox}
	\begin{corollary}\label{lem:sufficient_condition_for_error_bound}
		Let $\bar x\in\XX$ be an M-stationary point for \eqref{eq:P}, and fix $\bar y\in\Lambda(\bar x)$.
		If the second-order condition \eqref{eq:SSOSC} is valid, and if we have
		\begin{equation}\label{eq:crucial_assumption_for_error_bound}
			\forall u\in\XX,\,\forall \eta\in D(\partial g)(c(\bar x),\bar y)(c^\prime(\bar x)u)\colon\quad
			\innprod{\eta}{c^\prime(\bar x)u}
			\geq
			\d^2 g(c(\bar x),\bar y)(c^\prime(\bar x)u)
		\end{equation}
		and
		\begin{equation}\label{eq:CQ_uniqueness_of_multiplier}
			D(\partial g)(c(\bar x),\bar y)(0)\cap\ker c^\prime(\bar x)^\top=\{0\},
		\end{equation}
		then there are a constant $\varrho_{\mathrm{u}}>0$ and
		a neighborhood $U$ of $(\bar x,c(\bar x),\bar y)$ such that
		the upper estimate \eqref{eq:error_bound}
		holds for each triplet $(x,z,y)\in U\cap(\R^n\times\dom g\times\R^m)$.
	\end{corollary}
\end{mybox}
\begin{proof}
	We just show that the qualification condition \eqref{eq:strong_metric_subregularity} is valid.
	Then the assertion follows from \cref{lem:strong_ms_yields_error_bound}.
	
	Thus, pick $u\in\R^n$ and $\eta\in\R^m$ with
	$0=\nabla^2_{xx}\LL(\bar x,\bar y)u+c^\prime(\bar x)^\top\eta$ 
	and $\eta\in D(\partial g)(c(\bar x),\bar y)(c^\prime(\bar x)u)$.
	Taking the scalar product of the equation with $u$ gives
	$\nabla^2_{xx}\LL(\bar x,\bar y)[u,u]+\innprod{\eta}{c^\prime(\bar x)u}=0$,
	so that \eqref{eq:crucial_assumption_for_error_bound} gives
	$\nabla^2_{xx}\LL(\bar x,\bar y)[u,u]+\d^2g(c(\bar x),\bar y)(c^\prime(\bar x)u)\leq 0$.
	If $u\notin\mathcal C(\bar x)$, we have
	\begin{equation}\label{eq:not_in_critical_cone}
		\d g(c(\bar x))(c^\prime(\bar x)u)>-f^\prime(\bar x)u=\innprod{\bar y}{c^\prime(\bar x)u}
	\end{equation}
	from $\nabla_x\LL(\bar x,\bar y)=0$,
	which gives $\d^2 g(c(\bar x),\bar y)(c^\prime(\bar x)u)=\infty$,
	see \cite[formula~(5)]{BenkoMehlitz2023},
	and, thus, a contradiction.
	Hence, we have $u\in\mathcal C(\bar x)$, and \eqref{eq:SSOSC} gives $u=0$.
	Thus, $\eta\in D(\partial g)(c(\bar x),\bar y)(0)\cap\ker c^\prime(\bar x)^\top$,
	and \eqref{eq:CQ_uniqueness_of_multiplier} yields $\eta=0$.
	Consequently, \eqref{eq:strong_metric_subregularity} is valid, and the assertion follows.
\end{proof}

\begin{remark}\label{rem:error_bound_for_convex_piecewise_quadratic_functions}
	Let us note that \eqref{eq:crucial_assumption_for_error_bound} is valid
	whenever $g$ is prox-regular, subdifferentially continuous, and twice epi-differentiable
	at $c(\bar x)$ for $\bar y$, see \cref{lem:subderivative_vs_graphical_derivative}.
	Due to \cite[Ex.\ 13.30]{rockafellar1998variational},
	each proper, lower semicontinuous, convex function is prox-regular and
	subdifferentially continuous on its domain.
	Exemplary, whenever $g$ is a convex piecewise linear-quadratic function or the indicator
	function of the second-order cone, then it is twice epi-differentiable as well,
	see \cite[Prop.\ 13.9]{rockafellar1998variational} and \cite[Thm~3.1]{HangMordukhovichSarabi2020}.
\end{remark}

Observe that validity of \eqref{eq:strong_metric_subregularity} is equivalent to
validity of \eqref{eq:CQ_uniqueness_of_multiplier} and
\begin{equation}\label{eq:a_really_crucial_CQ}
	0=\nabla^2_{xx}\LL(\bar x,\bar y)u+c^\prime(\bar x)^\top \eta
	,\,
	\eta\in D(\partial g)(c(\bar x),\bar y)(c^\prime(\bar x)u)
	\implies
	u=0.
\end{equation}
The proof of \cref{lem:sufficient_condition_for_error_bound} shows that
validity of \eqref{eq:SSOSC} and \eqref{eq:crucial_assumption_for_error_bound}
implies that \eqref{eq:a_really_crucial_CQ} holds.
Let us elaborate on \eqref{eq:CQ_uniqueness_of_multiplier}.
Similar considerations in a much more specific setting can be found
in \cite[Sec.\ 8]{MohammadiMordukhovichSarabi2022} and
\cite[Sec.\ 4]{MordukhovichSarabi2019} where $g$ is assumed to be
a convex function of special type. 
Notice that the results obtained in \cite{MohammadiMordukhovichSarabi2022,MordukhovichSarabi2019}
are, expectedly, slightly stronger.

\begin{mybox}
\begin{lemma}\label{lem:uniqueness_of_multiplier_CQ}
	Let $\bar x\in\XX$ be an M-stationary point for \eqref{eq:P}, and fix $\bar y\in\Lambda(\bar x)$.
	We investigate the mapping $\ffunc{H}{\R^m}{\R^n\times\R^m}$ given by
	\begin{equation}\label{eq:definition_of_H}
		\forall y\in\R^m\colon\quad
		H(y)
		\coloneqq
		(\nabla_x\LL(\bar x,y),y)-\{0\}\times\partial g(c(\bar x)).
	\end{equation}
	Then \eqref{eq:CQ_uniqueness_of_multiplier} implies that $H$ is strongly metrically
	subregular at $(\bar y,(0,0))$, and the converse holds true whenever
	$(\partial g)^{-1}$ is metrically subregular at $(\bar y,c(\bar x))$.
\end{lemma}
\end{mybox}
\begin{proof}
	Patterning the proof of \cref{lem:strong_ms_yields_error_bound}, we find
	\[
		DH(\bar y,(0,0))(\eta)=(c^\prime(\bar x)^\top\eta,\eta)-\{0\}\times T_{\partial g(c(\bar x))}(\bar y).
	\]
	Furthermore, we have $T_{\partial g(c(\bar x))}(\bar y)\subseteq D(\partial g)(c(\bar x),\bar y)(0)$,
	and the converse holds true whenever $(\partial g)^{-1}$ is metrically subregular at
	$(\bar y,c(\bar x))$, see \cite[Thm~3.2]{BenkoMehlitz2022}.
	Thus, \eqref{eq:CQ_uniqueness_of_multiplier} implies
	\[
		\ker DH(\bar y,(0,0))=\{0\},
	\]
	and the converse holds true under the additional subregularity of $(\partial g)^{-1}$.
	Hence, the assertion follows from the Levy--Rockafellar criterion.
\end{proof}

\begin{mybox}
	\begin{corollary}\label{cor:unique_multiplier}
		Let $\bar x\in\XX$ be an M-stationary point for \eqref{eq:P}, and fix $\bar y\in\Lambda(\bar x)$.
		Then the following assertions hold.
	\begin{enumerate}
		\item \label{item:cor:unique_multiplier:cap}%
			If \eqref{eq:CQ_uniqueness_of_multiplier} holds, then there is a neighborhood
			$V\subseteq\R^m$ of $\bar y$ such that $\Lambda(\bar x)\cap V=\{\bar y\}$.
			Particularly, if $\Lambda(\bar x)$ is convex (which happens if $\partial g(c(\bar x))$ is convex),
			then $\Lambda(\bar x)=\{\bar y\}$.
		\item \label{item:cor:unique_multiplier:metricsubreg}%
			If $\Lambda(\bar x)=\{\bar y\}$, 
			if the mapping $H$ from \eqref{eq:definition_of_H} 
			is metrically subregular at $(\bar y,(0,0))$,
			and if $(\partial g)^{-1}$ is metrically subregular at $(\bar y,c(\bar x))$
			(both subregularity assumptions are satisfied if $\partial g$ is a polyhedral mapping
			as this also gives polyhedrality of $H$),
			then \eqref{eq:CQ_uniqueness_of_multiplier} holds.
	\end{enumerate}
	\end{corollary}
\end{mybox}
\begin{proof}
	Due to \cref{lem:uniqueness_of_multiplier_CQ}, the assumptions 
	in statement~\ref{item:cor:unique_multiplier:cap}
	guarantee that
	$H$ from \eqref{eq:definition_of_H} is strongly metrically subregular at $(\bar y,(0,0))$.
	Hence, we find a neighborhood $V\subseteq\R^m$ of $\bar y$ and a constant $\kappa>0$ such that
	\begin{equation}\label{eq:strong_metric_subregularity_of_H}
		\forall y\in V\colon\quad
		\norm{y-\bar y}
		\leq
		\kappa\bigl(\norm{\nabla_x\LL(\bar x,y)}+\dist(y,\partial g(c(\bar x)))\bigr).
	\end{equation}
	Particularly, this estimate shows that whenever $y\in V$ is different from $\bar y$, 
	then $y\notin\Lambda(\bar x)$. The additional statement in assertion \ref{item:cor:unique_multiplier:cap} readily follows.
			
	For assertion \ref{item:cor:unique_multiplier:metricsubreg},
	notice first that $H^{-1}(0,0)=\Lambda(\bar x)$ is valid.
	Then metric subregularity of $H$
	at $(\bar y,(0,0))$ together with $\Lambda(\bar x)=\{\bar y\}$ show
	that \eqref{eq:strong_metric_subregularity_of_H} is valid
	for some neighborhood $V\subseteq\R^m$ of $\bar y$ and some constant $\kappa>0$.
	Hence, $H$ is strongly metrically subregular at $(\bar y,(0,0))$.
	Finally, metric subregularity of $(\partial g)^{-1}$ at $(\bar y,c(\bar x))$ and
	\cref{lem:uniqueness_of_multiplier_CQ} can be used to infer validity of 
	\eqref{eq:CQ_uniqueness_of_multiplier}.
\end{proof}

In the following remark, we comment on \eqref{eq:a_really_crucial_CQ}.

\begin{remark}\label{rem:on_the_really_crucial_CQ}
	Let $\bar x\in\XX$ be an M-stationary point for \eqref{eq:P}, and fix $\bar y\in\Lambda(\bar x)$.
	Suppose that the second-order condition \eqref{eq:SSOSC} is valid, 
	and that $0\in\dom \d^2g(c(\bar x),\bar y)$.
	We first note that this means that $\bar u \coloneqq 0$ is the uniquely determined global minimizer
	of
	\begin{equation}\label{eq:second_order_minimization_problem}
		\minimize_{u}
		{}\quad{}
		\frac12\nabla^2_{xx}\LL(\bar x,\bar y)[u,u]+\frac12\d^2g(c(\bar x),\bar y)(c^\prime(\bar x)u).
	\end{equation}
	In order to see this, one has to observe two facts.
	First, for each $u\notin \mathcal C(\bar x)$, 
	we have \eqref{eq:not_in_critical_cone}	which gives 
	$\d^2g(c(\bar x),\bar y)(c^\prime(\bar x)u)=\infty$ as mentioned earlier.
	Secondly, $\d^2g(c(\bar x),\bar y)(0)=0$ follows from \cref{lem:trivial_second_subderivative}.
	
	Under the assumptions of \cref{lem:subderivative_vs_graphical_derivative}, the limiting subdifferential
	of the scaled second subderivative
	$v\mapsto\tfrac12\d^2g(c(\bar x),\bar y)(v)$ can be computed in terms of the graphical derivative
	of $\partial g$ at $(c(\bar x),\bar y)$. 
	Hence, under some suitable assumptions, the chain rule from \cite[Cor.\ 4.6]{mordukhovich2018} can be applied in order
	to derive first-order necessary optimality conditions for problem \eqref{eq:second_order_minimization_problem},
	and these conditions take the following shape:
	\[
		0=\nabla^2_{xx}\LL(\bar x,\bar y)u+c^\prime(\bar x)^\top \eta ,\quad
		\eta\in D(\partial g)(c(\bar x),\bar y)(c^\prime(\bar x)u).
	\]
	That is why \eqref{eq:a_really_crucial_CQ} demands, roughly speaking, that $\bar u \coloneqq 0$ is the uniquely determined
	stationary point of \eqref{eq:second_order_minimization_problem}.
	This is clearly different from
	postulating that this point is the uniquely determined global minimizer of this problem, i.e., \eqref{eq:SSOSC}.
	Equivalence of these conditions can only be guaranteed under some additional convexity of the
	second subderivative.
	
	Following \cite[Def.\ 3.1]{MordukhovichSarabi2018}, validity of \eqref{eq:a_really_crucial_CQ} demands that
	the multiplier $\bar y$ is so-called \emph{noncritical}.
	The concept of critical multipliers dates back to 
	\cite{Izmailov2005,IzmailovSolodov2012} where it has been introduced for standard nonlinear programs.
	In \cite{IzmailovSolodov2012}, it has been pointed out that the presence of critical multipliers
	slows down the convergence of Newton-type methods when applied for the solution of stationarity systems,
	and this observation can be extended to certain composite optimization problems as
	shown in \cite{MordukhovichSarabi2018,MordukhovichSarabi2019}.
	Let us mention that \cite[Thm~4.1]{MordukhovichSarabi2018} and \cite[Thm~5.6]{MordukhovichSarabi2019} justify
	that \eqref{eq:a_really_crucial_CQ} is equivalent to a local primal-dual (upper) error bound,
	based on a not necessarily unique Lagrange multiplier, whenever $g$ is a function whose epigraph is a convex
	polyhedral set or the indicator function of a so-called $\mathcal C^2$-cone reducible set,
	and in the latter case, further assumptions are required. 
	For brevity, we abstain here from investigating further refinements of the error bound \eqref{eq:error_bound}
	to situations where the Lagrange multiplier is not uniquely determined, 
	but indicate that this is an interesting topic for future research.
	
	The above puts our comments from \cref{rem:error_bound_for_convex_piecewise_quadratic_functions}
	into some new light.
	On the one hand, our arguments highlight that in situations where the second subderivative of $g$ is not convex, 
	\eqref{eq:SSOSC} might be too weak to yield the error bound of our interest.
	On the other hand, in the absence of any additional assumptions, \eqref{eq:strong_metric_subregularity} may not be
	sufficient for \eqref{eq:SSOSC} as uniqueness of stationary points says nothing about the existence of a global minimizer
	for \eqref{eq:second_order_minimization_problem}.
	However, the condition
	\begin{equation}\label{eq:some_new_second_order_condition_for_error_bound}
		\forall u\in\R^n\setminus\{0\},\,\forall \eta\in D(\partial g)(c(\bar x),\bar y)(c^\prime(\bar x)u)\colon\quad
		\nabla^2_{xx}\LL(\bar x,\bar y)[u,u]+\innprod{\eta}{c^\prime(\bar x)u}>0
	\end{equation}
	is clearly sufficient for \eqref{eq:a_really_crucial_CQ} and, 
	together with \eqref{eq:CQ_uniqueness_of_multiplier},
	gives \eqref{eq:strong_metric_subregularity}. 
\end{remark}

The following example shows that \eqref{eq:a_really_crucial_CQ}
does not necessarily imply \eqref{eq:SSOSC}.
Furthermore, \cref{ex:SSOSC_does_not_give_error_bound} below visualizes that \eqref{eq:SSOSC} does not
necessarily imply \eqref{eq:a_really_crucial_CQ}.
These two conditions are, thus, independent in general as indicated in \cref{rem:on_the_really_crucial_CQ}.

\begin{example}
	We consider \eqref{eq:P} for the functions $\func{f,c,g}{\R}{\R}$ given by
	$f(x)\coloneqq\tfrac12x^2$,
	$c(x)\coloneqq x$,
	$g(z)\coloneqq -z^2$,
	and choose $\bar x$ to be the origin in $\R$.
	Note that $\bar x$ is M-stationary with $\Lambda(\bar x)=\{0\}$.
	Thus, we consider the uniquely determined multiplier $\bar y \coloneqq 0$.
	Obviously, $\bar x$ is a strict local maximizer of \eqref{eq:P}
	which is why \eqref{eq:SSOSC} fails to hold at $\bar x$ for $\bar y$, 
	see \cref{cor:SSOSC}.
	Clearly, we have $\nabla_{xx}^2\LL(\bar x,\bar y)=1$, and
	twice continuous differentiability of $g$ gives
	$D(\partial g)(c(\bar x),\bar y)(c^\prime(\bar x)u)=\{-2u\}$
	for each $u\in\R$. 
	Hence, one can easily check that \eqref{eq:a_really_crucial_CQ} is valid.
\end{example}

\begin{remark}[Geometric constraints]\label{rem:CQs_for_set_indicators}
	In the special case where $g\coloneqq \indicator_D$ holds for some closed set $D\subseteq\R^m$,
	the qualification conditions 
	\eqref{eq:strong_metric_subregularity},
	\eqref{eq:CQ_uniqueness_of_multiplier},
	\eqref{eq:a_really_crucial_CQ},
	and
	\eqref{eq:some_new_second_order_condition_for_error_bound}
	involve the graphical derivative of the (limiting) normal cone mapping
	associated with $D$. For several different choices of $D$, including
	convex cones (like the semidefinite or the second-order cone) or
	convex sets given via smooth convex inequality constraints, explicit
	ready-to-use formulas for this variational object are available, see
	e.g.\ \cite{GfrererOutrata2016,Wachsmuth2017}.
	For diverse nonconvex sets $D$ of special structure, like sparsity sets of type
	$\{x\in\R^n\,|\,\norm{x}_0\leq \kappa\}$, $k\in\{1,\ldots,n-1\}$, 
	the explicit computation of this tool is possible as well.
\end{remark}

As we have seen above, the \emph{upper} error bound in \eqref{eq:error_bound} only holds in the presence of
comparatively strong assumptions. Unfortunately, due to our definition of $\Theta$ 
in \eqref{eq:ErrorBound:error} which comprises the distance to the subdifferential of $g$,
the converse \emph{lower} error bound seems to demand even more prohibitive assumptions,
as the following \cref{rem:lower_error_bound} illustrates.
Nonetheless, we will circumnavigate this potentially crucial observation later on in
\cref{sec:ALM} by the design of our algorithm.

\begin{remark}\label{rem:lower_error_bound}
	Let $\bar x\in\XX$ be an M-stationary point of \eqref{eq:P} and pick $\bar y\in\Lambda(\bar x)$.
	Then, relying on the definition of $\Theta$ in \eqref{eq:ErrorBound:error},
	it appears indispensable for estimating the distance to 
	the subdifferential to assume its \emph{inner calmness},
	see \cite[Def.\ 2.2]{BenkoGfrererOutrata2019} and \cite[Sec.\ 2]{Benko2021} for a discussion of this property.
	In particular, $\partial g$ shall be inner calm at $(c(\bar{x}),\bar{y})$,
	which entails the existence of $\kappa > 0$ and a neighborhood $V\subseteq\YY$ of $c(\bar x)$ such that
	\[
	\forall z\in V\colon\quad
	\dist(\bar y,\partial g(z))
	\leq
	\kappa\norm{z-c(\bar x)}.
	\]
	With this property at hand and exploiting \eqref{eq:Mstationary},
	for each triplet $(x,z,y)\in\XX\times\dom g\times\YY$ such that $z\in V$,
	the triangle inequality yields
	\begin{align*}
	\Theta(x,z,y)
	{}={}&
	\norm{\nabla_x\LL(x,y)} + \norm{c(x)-z}  + \dist(y,\partial g(z))
	\\
	{}\leq{}&
	\norm{\nabla_x\LL(x,y)-\nabla_x\LL(\bar x,\bar y)}
	+
	\norm{c(x)-c(\bar x)} + \norm{c(\bar x)-z}
	+
	\norm{y-\bar y}+\dist(\bar y,\partial g(z))
	\\
	{}\leq{}&
	\norm{\nabla_x\LL(x,y)-\nabla_x\LL(\bar x,\bar y)}
	+
	\norm{c(x)-c(\bar x)} + (\kappa+1)\norm{z-c(\bar x)} + \norm{y-\bar y}.
	\end{align*}
	Then, noting that $\nabla_x\LL$ and $c$ are locally Lipschitz continuous by
	\cref{ass:P}\ref{ass:f},
	there are a constant $\varrho_{\mathrm{l}}>0$ and
	a neighborhood $U$ of $(\bar x,c(\bar x),\bar y)$
	such that, for each $(x,z,y)\in U\cap(\XX\times\dom g\times\YY)$, we obtain
	the lower estimate
	\begin{equation*}
		\varrho_{\mathrm{l}}\,\Theta(x,z,y)
		\leq
		\norm{x-\bar x}+\norm{z-c(\bar x)}+\norm{y-\bar y},
	\end{equation*}
	which patterns the upper counterpart in \eqref{eq:error_bound}.
	
	However, inner calmness of the subdifferential is an impractical assumption, even for convex $g$,
	and it would restrict our considerations mainly to points where $g$ is smooth in practice.
	Exemplary, consider the absolute value function $g \coloneqq |\cdot|$.
	Since $\partial g$ is single-valued on $\R \setminus \{0\}$,
	one can easily check that inner calmness of $\partial g$ at $(0,\bar y)$ fails for every
	$\bar y\in[-1,1]=\partial g(0)$.
\end{remark}
	
In the following \cref{sec:ALM},
we will \emph{not} rely on any additional property of the subdifferential,
but leverage instead the algorithmic scheme
to derive a lower error bound \emph{along the iterates}, 
see \cref{lem:lower_error_bound:ALM} below.

\section{Augmented Lagrangian scheme and convergence}
\label{sec:ALM}

This section is devoted to describing a numerical scheme for solving \eqref{eq:P}
and to investigating its convergence properties under suitable assumptions.
In particular,
we consider the implicit AL scheme from \cite[Alg.\ 4.1]{demarchi2024implicit},
so called because it makes use of 
the AL function $\func{\LL_\mu}{\XX \times \YY}{\R}$ associated to \eqref{eq:P},
as defined in \eqref{eq:L}.
It deviates in this respect from \cite[Alg.\ 1]{demarchi2023constrained},
which builds upon \eqref{eq:ALMz:P} and treats the auxiliary variable explicitly,
see \cite{benko2021implicit} for a discussion on implicit variables
(and concealed benefits thereof)
in optimization.

\subsection{Implicit approach}\label{sec:ALM:approach}

The numerical method considered for addressing \eqref{eq:P} is stated in \cref{alg:ALM}.
Fitting into the AL framework \cite{bertsekas1996constrained,birgin2014practical,conn1991globally},
the main step of the iterative procedure involves minimizing the AL function with respect to the primal variable.
At the $k$-th iteration of the AL scheme, with some given penalty parameter $\mu_k > 0$ and multiplier estimate $\hat{y}^k \in \YY$, a subproblem involving the minimization of $\LL_{\mu_k}(x,\hat{y}^k)$ over $x \in \XX$ has to be solved approximately.
However,
the AL function may lack regularity since,
for $g$ nonconvex, the Moreau envelope $g^\mu$ is in general \emph{not} continuously differentiable.
Therefore, the concept of (approximate) $\Upsilon$-stationarity introduced in \cref{sec:VA} plays a role in characterizing adequate solutions of this AL subproblem,
see \cref{rem:Upsilon_stationarity_AL} as well.
Let us note that, in practice, such points can be computed with the aid of a
nonmonotone descent method, see \cite[Sec.\ 5]{demarchi2024implicit} for details.
Then, following classical update rules \cite{birgin2014practical},
the multiplier estimate $\hat{y}$ and penalty parameter $\mu$ are adjusted,
along with the subproblem's tolerance $\varepsilon$.

Compared to the classical AL approach for the solution of nonlinear programs, see \cite{bertsekas1996constrained,conn1991globally},
this variant uses a safeguarded update rule for the Lagrange multipliers and has stronger global convergence properties, as demonstrated in \cite{kanzow2018augmented}.
The safeguarded multiplier estimate $\hat{y}^k$ is drawn from a bounded set $\Ybounded\subseteq\YY$ at \cref{step:ALM:ysafe}.
In practice,
it is advisable to choose the safeguarded estimate $\hat{y}^k$ as the projection of the multiplier $y^{k-1}$ onto $\Ybounded$.
We refer to \cite[Sec.\ 4.1]{birgin2014practical} for a detailed discussion.

The monotonicity test at \cref{step:ALM:penalty} is adopted to monitor primal infeasibility along the iterates and update the penalty parameter accordingly.
Aimed at driving $V_k$ to zero,
the penalty parameter $\mu_k$ is reduced in case of insufficient improvement.

\begin{algorithm2e}[htb]
	\DontPrintSemicolon
	\KwData{$\mu_0 \in (0,\mu_g)$, $\theta\in(0,1)$, $\kappa \in (0,1)$, $\Ybounded\subseteq\YY$ nonempty bounded\;}
	\For{$k = 0,1,2\ldots$}{
		Select $\hat{y}^k \in \Ybounded$ and $\varepsilon_k \geq 0$\label{step:ALM:ysafe}\; 
		Compute a pair $(x^k,z^k)$ certificating $\varepsilon_k$-$\Upsilon$-stationary of $x^k$ for $\LL_{\mu_k}(\cdot,\hat{y}^k)$\label{step:ALM:subproblem}:
		\begin{center}
			$\displaystyle
			\norm{\nabla_x \LLslack_{\mu_k}(x^k,z^k,\hat{y}^k)} \leq \varepsilon_k ,
			\qquad
			z^k \in \prox_{\mu_k g}(c(x^k) + \mu_k\hat{y}^k)
			$
		\end{center}
		Set 
		\label{step:ALM:y}%
		$y^k \gets \hat{y}^k + \mu_k^{-1} [c(x^k) - z^k]$
		and $V_k \gets \norm{c(x^k) - z^k}$\;
		\lIf{$k=0$ \KwOr $ V_k \leq \theta V_{k-1}$\label{step:ALM:penalty}}{
			$\mu_{k+1} \gets \mu_k$, \textbf{else} $\mu_{k+1} \gets \kappa \mu_k$
		}
	}
	\caption{Safeguarded implicit AL method for \eqref{eq:P}}
	\label{alg:ALM}
\end{algorithm2e}

Before investigating the convergence properties of \cref{alg:ALM},
we provide some characterizations of the iterates $\{(x^k,z^k,y^k)\}$.
These are direct consequences of $z^k$ being a certificate of $\varepsilon_k$-$\Upsilon$-stationarity for $x^k$,
by \cref{step:ALM:subproblem},
and of the dual update at \cref{step:ALM:y},
see \cref{rem:Upsilon_stationarity_AL} as well.

\begin{mybox}
	\begin{proposition}\label{lem:ALM:exact_complementarity}
		Let $\{(x^k,z^k,y^k)\}$ be a sequence generated by \cref{alg:ALM}.
		Then, for each $k\in\N$,
		$z^k \in \prox_{\mu_k g}(c(x^k)+\mu_k\hat y^k) \subseteq \dom g$,
		$\norm{	\nabla f(x^k) + c^\prime(x^k)^\top y^k } \leq \varepsilon_k$,
		and
		$y^k \in \partial g(z^k)$.
	\end{proposition}
\end{mybox}

Throughout the convergence analysis, based on, and extending, that of \cite{demarchi2024implicit}, 
it is assumed that \cref{alg:ALM} is well-defined, 
thus requiring that each subproblem at \cref{step:ALM:subproblem} 
admits an approximately stationary point.
Moreover, the existence of some accumulation point $\bar x$ for a sequence $\{ x^k \}$ 
generated by \cref{alg:ALM} requires, in general, coercivity or (level) boundedness arguments.

While asymptotic M-stationarity as given in \cite[Def.\ 3.3]{demarchi2024implicit} demands the existence of a sequence
of, in a certain sense, approximately M-stationary points for \eqref{eq:P} converging to the
point of interest, no quantitative bound on this \emph{approximativity} is required.
However, for \cref{alg:ALM} to terminate, there is a need for an approximate version of
the M-stationarity concept of \cref{def:Mstationary}.
We refer to the notion delineated in \cite[Def.\ 3.2]{demarchi2024implicit},
which comes along with an explicit bound quantifying violation of M-stationarity,
while aligning with the asymptotic stipulation.

\subsection{Global convergence}\label{sec:GlobalConvergence}

In this section, we are concerned with global convergence properties related to \cref{alg:ALM},
i.e., we are going to study properties of accumulation points of sequences it generates,
regardless of how it is initialized.
For that purpose, we will assume that \cref{alg:ALM} is well-defined and produces
an infinite sequence of iterates.

Our first results pertain to a \emph{global} optimization perspective
on the subproblems at \cref{step:ALM:subproblem} of \cref{alg:ALM},
compare \cite[Ch.\ 5]{birgin2014practical} and \cite[Sec.\ 4]{kanzow2018augmented}.
Solving each subproblem up to approximate global optimality,
not necessarily with vanishing inner tolerance, 
one finds in the limit a global minimizer of the infeasibility measure,
see \cite[Lem.\ 4.2]{BoergensKanzowSteck2019} for a related result.
\begin{mybox}
	\begin{lemma}\label{lem:global_minimizer_constraint_violation}
		Let $\{(x^k,z^k,y^k)\}$ be a sequence generated by \cref{alg:ALM} with $\{\varepsilon_k\}$ bounded.
		Assume that $\dom g$ is closed and that
		\begin{equation}\label{eq:approximate_global_minimizer_subproblem}
			\forall k\in\N,\,\forall x\in\XX\colon\quad
			\LL_{\mu_k}(x^k,\hat{y}^k) \leq \LL_{\mu_k}(x,\hat{y}^k) + \varepsilon_k.
		\end{equation}
		Fix an accumulation point $\bar{x}\in\XX$ of $\{x^k\}$.
		Then $\bar{x}$ is a global minimizer of $\dist(c(\cdot),\dom g)$.
		In particular, $\bar{x}$ is feasible if the feasible set of \eqref{eq:P} is nonempty.
	\end{lemma}
\end{mybox}
\begin{proof}
	Let us consider two cases.
	If $\{\mu_k\}$ remains bounded away from zero, then \cref{step:ALM:penalty} demands 
	that $\norm{c(x^k) - z^k} \to 0$.
	As we have $\{z^k\}\subseteq\dom g$ from \cref{lem:ALM:exact_complementarity},
	it follows that
	\begin{equation*}
		0 \leq \dist(c(x^k),\dom g) \leq \norm{ c(x^k) - z^k } \to 0.
	\end{equation*}
	Owing to $\dom g$ being closed,
	we obtain $c(\bar{x})\in\dom g$,
	proving that $\bar{x}$ is feasible for \eqref{eq:P}.
	
	Consider now the case $\mu_k\downto 0$.
	Let $\{x^k\}_{k\in K}$ be a subsequence such that $x^k \to_K \bar{x}$.
	Then \eqref{eq:approximate_global_minimizer_subproblem} guarantees
	\begin{equation*}
		\forall x\in\XX\colon\quad
		f(x^k) + g^{\mu_k}(c(x^k) + \mu_k \hat{y}^k)
		\leq
		f(x) + g^{\mu_k}(c(x) + \mu_k \hat{y}^k) + \varepsilon_k
		.
	\end{equation*}
	Multiplying by $\mu_k$, taking the lower limit as $k\to_K\infty$, and using boundedness of $\{f(x^k)\}_{k\in K}$ and $\{\varepsilon_k\}$ yield
	\[
		\liminf\limits_{k\to_K\infty}
		\mu_k g^{\mu_k}(c(x^k) + \mu_k \hat{y}^k)
		\leq
		\liminf\limits_{k\to_K\infty}
		\mu_k g^{\mu_k}(c(x) + \mu_k \hat{y}^k)
	\]
	for each $x\in\XX$.
	Together with \cref{lem:dist_to_domain}, this gives
	\begin{align*}
		\frac{1}{2} \dist^2(c(\bar{x}), \dom g)
		{}={}&
		\liminf_{k\to_K\infty}
		\inf_z\left\{ \mu_k g(z) + \frac{1}{2} \norm{ z - (c(x^k) + \mu_k \hat{y}^k) }^2 \right\} \\
		{}={}&
		\liminf_{k\to_K\infty}
		\mu_k g^{\mu_k}(c(x^k) + \mu_k \hat{y}^k) 
		{}\leq{}
		\liminf_{k\to_K\infty}
		\mu_k g^{\mu_k}(c(x) + \mu_k \hat{y}^k) \\
		{}={}&
		\liminf_{k\to_K\infty}
		\inf_z\left\{ \mu_k g(z) + \frac{1}{2} \norm{ z - (c(x) + \mu_k \hat{y}^k) }^2 \right\} 
		{}={}
		\frac{1}{2} \dist^2(c(x), \dom g)
	\end{align*}
	for all $x\in\XX$, where we used boundedness of $\{\hat{y}^k\}$ and $\mu_k\downto 0$.
	This shows that $\bar{x}$ globally minimizes $\dist(c(\cdot),\dom g)$.
	Then, by closedness of $\dom g$,
	if \eqref{eq:P} is feasible, so is $\bar{x}$,
	concluding the proof.
\end{proof}

Notice that, in the proof of \cref{lem:global_minimizer_constraint_violation},
the requirement of closed $\dom g$ does not affect the case with $\mu_k \downto 0$,
and that we did not exploit the actual qualitative requirements regarding the
subproblem solver stated in \cref{step:ALM:subproblem}, apart from the necessity
of having a sequence $\{z^k\}\subset\dom g$ at hand.
Indeed, approximate global optimality of $x^k$ for the subproblem in the sense
of \eqref{eq:approximate_global_minimizer_subproblem} is all we needed.
Conversely, the upcoming two theorems rely on $z^k\in\prox_{\mu_k g}(c(x^k)+\mu_k\hat y^k)$,
additionally, but the specification
$\norm{\nabla_x \LLslack_{\mu_k}(x^k,z^k,\hat{y}^k)} \leq \varepsilon_k$
is still not required.

If, in addition to the assumptions of \cref{lem:global_minimizer_constraint_violation},
the sequence $\{\varepsilon_k\}$ is chosen so that $\varepsilon_k \downtoeq 0$, 
then (primal) accumulation points of the sequence generated by \cref{alg:ALM} correspond
to global minimizers of \eqref{eq:P}, see \cite[Thm~4.12]{BoergensKanzowSteck2019} for a related result.

\begin{mybox}
	\begin{theorem}\label{thm:global_optimality}
		Let $\{(x^k,z^k,y^k)\}$ be a sequence generated by \cref{alg:ALM} with $\varepsilon_k \downtoeq 0$.
		Let \eqref{eq:approximate_global_minimizer_subproblem} hold, and
		assume that the feasible set of \eqref{eq:P} is nonempty while $\dom g$ is closed.
		Then $\limsup_{k \to \infty} (f(x^k) + g(z^k)) \leq \varphi(x)$ holds for all feasible $x\in\XX$.
		Moreover, every accumulation point $\bar{x}\in\XX$ of $\{x^k\}$ is globally optimal for \eqref{eq:P}.
		For any index set $K \subseteq \N$ such that $x^k \to_K \bar{x}$, 
		it is also $z^k \to_K c(\bar{x})$ and $f(x^k)+g(z^k)\to_K\varphi(\bar x)$.
	\end{theorem}
\end{mybox}
\begin{proof}
	Let $x\in\XX$ be a fixed feasible point of \eqref{eq:P}.
	Then, due to \eqref{eq:approximate_global_minimizer_subproblem}, \cref{lem:upper_bound_AL}, and \cref{lem:ALM:exact_complementarity}, 
	we have
	\begin{multline*}
		f(x^k) + g(z^k) - \frac{\mu_k}{2}\norm{ \hat{y}^k }^2
		{}\leq {}
		f(x^k) + g(z^k) + \frac{1}{2\mu_k} \norm{ c(x^k)+\mu_k \hat{y}^k - z^k }^2 - \frac{\mu_k}{2}\norm{ \hat{y}^k }^2
		\\
		{}={}
		\LL_{\mu_k}(x^k,\hat{y}^k) 
		{}\leq{}
		\LL_{\mu_k}(x,\hat{y}^k) + \varepsilon_k
		{}\leq{}
		\varphi(x) + \varepsilon_k
		< \infty
	\end{multline*}
	for all $k\in\N$.
	If $\mu_k\downto 0$, then $\mu_k \norm{ \hat{y}^k }^2 \to 0$ by boundedness of $\{\hat{y}^k\}$.
	In this case, $\varepsilon_k \downtoeq 0$ and $z^k \in \dom g$ imply that $\limsup_{k \to \infty} (f(x^k) + g(z^k)) \leq \varphi(x)$.
	
	Let us now focus on the case where $\{\mu_k\}$ is bounded away from zero.
	This is possible only if $\norm{c(x^k) - z^k} \to 0$ by \cref{step:ALM:penalty}.
	Similar as above, we find
	\begin{equation*}
		f(x^k) + g(z^k) + \frac{1}{2\mu_k} \norm{c(x^k) - z^k}^2 + \innprod{\hat{y}^k}{c(x^k) - z^k}
		=
		\LL_{\mu_k}(x^k,\hat{y}^k)
		\leq
		\varphi(x)+\varepsilon_k.
	\end{equation*}
	Since $\{\hat{y}^k\}$ is bounded, $\norm{ c(x^k) - z^k } \to 0$, and $\varepsilon_k \downtoeq 0$,
	we can take the upper limit in the above estimate
	to find $\limsup_{k \to \infty} (f(x^k) + g(z^k)) \leq \varphi(x)$.
	Finally, let $\bar{x}$ be an accumulation point of $\{x^k\}$ and $K\subseteq\N$ an index set such that $x^k \to_K \bar{x}$.
	Then, as \eqref{eq:P} admits feasible points, $\bar{x}$ is feasible by \cref{lem:global_minimizer_constraint_violation}.
	
	Let us show that $\norm{c(x^k)-z^k} \to_K 0$ holds.
	If $\{\mu_k\}$ remains bounded away from zero, this is obvious by \cref{step:ALM:penalty}.
	In the case where $\mu_k\downto 0$, we can exploit feasibility of $\bar x$, boundedness of $\{\hat y^k\}$, and \cref{lem:dist_to_domain}
	to find
	\[
		\mu_kg(z^k)+\frac12\norm{z^k-c(x^k)-\mu_k\hat y^k}^2\to_K 0.
	\]
	Boundedness of $\{\hat y^k\}$ and $\mu_k\downto 0$ allow us to apply \cref{lem:norm_to_zero}
	which yields $\mu_k g(z^k)\to_K0$ as well as $\norm{z^k - c(x^k) - \mu_k\hat y^k}\to_K0$,
	and the latter gives $\norm{z^k-c(x^k)}\to_K0$.
	
	Due to $\norm{c(x^k)-z^k}\to_K0$, continuity of $c$ gives $z^k\to_K c(\bar x)$.
	Now, lower semicontinuity of $g$ yields
	\begin{multline*}
		\varphi(\bar{x}) 
		{}={}
		f(\bar x)+g(c(\bar x))
		{}\leq{}
		\liminf\limits_{k\to_K\infty} f(x^k)+\liminf\limits_{k\to_K\infty}g(z^k)\\
		{}\leq{}
		\liminf\limits_{k\to_K\infty} (f(x^k) + g(z^k) )
		{}\leq{}
		\limsup_{k \to_K \infty} (f(x^k) + g(z^k) )
		\leq 
		\varphi(x),
	\end{multline*}
	where the last inequality is due to the upper bound obtained previously in the proof.
	As $x$ is an arbitrary feasible point of \eqref{eq:P}, we have shown that $\bar{x}$ is globally optimal for \eqref{eq:P}.
	Finally, with the particular choice $x = \bar{x}$, the previous inequalities give that $f(x^k)+g(z^k)\to_K\varphi(\bar x)$,
	concluding the proof.
\end{proof}

Under an additional assumption on the multiplier estimate $\hat y^k$ in \cref{step:ALM:ysafe}, 
a stronger result can be proved that concerns the behavior of the iterates for infeasible problems.
By resetting the multiplier estimate when signs of infeasibility are detected,
the algorithm tends to minimize the objective function subject to minimal constraint violation,
see e.g.\ \cite[Thm~5.3]{birgin2014practical} for a related result.
\begin{mybox}
	\begin{theorem}\label{thm:optimality_subject_minimal_infeasibility}
		Let $\{(x^k,z^k,y^k)\}$ be a sequence generated by \cref{alg:ALM} with $\varepsilon_k \downtoeq 0$.
		Let \eqref{eq:approximate_global_minimizer_subproblem} hold,
		suppose that $\dom g$ is closed, and that,
		for all $k\in\N$, $\hat{y}^{k+1}  = 0$ if $y^k \notin \Ybounded$.
		Let $\bar{x}$ be an accumulation point of $\{x^k\}$.
		Then $\bar{x}$ is a global minimizer of $\dist(c(\cdot), \dom g)$ and,
		for all $(x,z) \in\XX \times \dom g$ such that $\| c(x) - z \| = \dist(c(\bar{x}), \dom g)$,
		it holds $\limsup_{k \to \infty} (f(x^k) + g(z^k)) \leq f(x) + g(z)$.
	\end{theorem}
\end{mybox}
\begin{proof}
	If $\dist(c(\bar{x}), \dom g) = 0$, namely $c(\bar{x}) \in \dom g$ by closedness of $\dom g$, then $\bar{x}$ is feasible and the claim follows from \cref{thm:global_optimality}.
	So, let us assume that $\dist(c(\bar{x}), \dom g) > 0$.
	Together with \cref{step:ALM:penalty} and \cref{lem:ALM:exact_complementarity}, this implies that $\mu_k\downto 0$.
	Since $\bar{x}$ is a global minimizer of $\dist(c(\cdot),\dom g)$ by \cref{lem:global_minimizer_constraint_violation}, then 
	\begin{equation}\label{eq:dist_to_domain}
		\forall k\in\N\colon\quad
		\| c(x^k) - z^k \| \geq \dist(c(x^k), \dom g) \geq \dist(c(\bar{x}), \dom g).
	\end{equation}
	Thus, by the dual update rule at \cref{step:ALM:y}, boundedness of $\Ybounded$, and $\mu_k\downto 0$, for all $k\in\N$ large enough it is $y^k \notin \Ybounded$.
	Therefore, by the estimate choice stated in the premises, it is $\hat{y}^k = 0$ for all large enough $k\in\N$.
	Then, for all $x\in\XX$, we have that
	\begin{multline*}
		f(x^k) + g(z^k) + \frac{1}{2\mu_k} \|c(x^k) - z^k\|^2
		{}={}
		f(x^k) + g^{\mu_k}(c(x^k))
		{}={}
		\LL_{\mu_k}(x^k,0) \\
		{}\leq{}
		\LL_{\mu_k}(x,0) + \varepsilon_k 
		{}={}
		f(x) + g^{\mu_k}(c(x)) + \varepsilon_k
	\end{multline*}
	for all $k\in \N$ large enough.
	This holds, in particular, for $x \in \XX$ a global minimizer of $\dist(c(\cdot),\dom g)$, namely such that $\dist(c(x), \dom g) = \dist(c(\bar{x}), \dom g)$.
	Then there is some $z \in \dom g$ such that $\|c(x) - z\| = \dist(c(\bar{x}), \dom g)$,
	and \eqref{eq:dist_to_domain} gives $\norm{c(x^k)-z^k}\geq\norm{c(x)-z}$.
	Hence, we find
	\begin{multline*}
		f(x^k) + g(z^k) + \frac{1}{2\mu_k} \|c(x^k) - z^k\|^2
		{}\leq{}
		f(x) + g^{\mu_k}(c(x)) + \varepsilon_k \\
		{}\leq{}
		f(x) + g(z) + \frac{1}{2\mu_k}\|c(x) - z\|^2 + \varepsilon_k 
		{}\leq{}
		f(x) + g(z) + \frac{1}{2\mu_k}\norm{c(x^k)-z^k}^2 + \varepsilon_k
	\end{multline*}
	for all $k\in \N$ large enough. 
	Subtracting the squared norm term on both sides and taking the upper limit yields the claim since $\varepsilon_k\downtoeq 0$.
\end{proof}

In \cref{lem:global_minimizer_constraint_violation} and \cref{thm:global_optimality,thm:optimality_subject_minimal_infeasibility},
it has been assumed that the AL subproblem can be solved up to approximate global optimality.
This, however, might be a delicate issue whenever $f$ or $g$ are nonconvex or $c$ is difficult enough.
In practice,
affordable solvers only have local scope and
return stationary points as candidate local minimizers.
Nevertheless,
(primal) accumulation points of a sequence generated by \cref{alg:ALM} can be shown to be at least \emph{asymptotically} M-, or AM-, stationary points for \eqref{eq:P},
see \cite[Def.\ 3.3, Thm~4.1]{demarchi2024implicit} for a detailed discussion.
Notice that the (approximate) global optimality is not relaxed to \emph{local} optimality, but to mere ($\Upsilon$-)stationarity,
while the subsequential $g$-attentive convergence of certain iterates is required.
We refer to \cite[Ex.\ 3.4]{demarchi2023constrained} for an illustration of the importance of attentive convergence.
It should be noted that a mild asymptotic regularity condition is enough to guarantee that AM-stationary point of
\eqref{eq:P} are indeed M-stationary, see \cite[Cor.\ 3.1]{demarchi2024implicit} or \cite[Cor.\ 2.7]{demarchi2023constrained}
for related results.
Thus, \cref{alg:ALM} is likely to compute M-stationary points of \eqref{eq:P}.

A result analogous to \cref{lem:global_minimizer_constraint_violation}
is available in this affordable setting, too, see
\cite[Prop.\ 4.5]{demarchi2024implicit},
whereas \cite[Prop.\ 4.3]{demarchi2024implicit}
provides some sufficient conditions for the
feasibility of accumulation points.

\subsection{Local convergence}\label{sec:LocalConvergence}

In this section, we investigate the behavior of \cref{alg:ALM} in the vicinity of stationary points of \eqref{eq:P} under various assumptions.
Particularly, we are interested in the existence of strict local minimizers of the AL subproblems in a neighborhood
of a strict local minimizer to \eqref{eq:P} and convergence rates associated with \cref{alg:ALM} in such situations.

\subsubsection{Existence of local minimizers}

Let us consider
the existence of local minimizers of the AL function $\LL_\mu(\cdot,\hat y)$
for $\mu>0$ sufficiently small and an approximate multiplier $\hat y\in Y$, where $Y\subseteq\YY$ is a bounded set, see \cref{alg:ALM}.
Note that, by construction, any such local minimizer would be an $\Upsilon$-stationary
point of $\LL_\mu(\cdot,\hat y)$, and \cref{step:ALM:subproblem} would be meaningful,
see \cref{rem:Upsilon_stationarity_AL} as well.

We fix some strict local minimizer $\bar x\in\XX$ of \eqref{eq:P} and
proceed as suggested in \cite[Sec.\ 7]{BoergensKanzowSteck2019}. For sufficiently small $r>0$, 
consider the localized AL subproblem
\begin{equation}\label{eq:subproblem_existence_local_minimizers}
	\minimize_{x\in\XX} \quad \LL_{\mu_k}(x,\hat{y}^k)
	\qquad
	\stt\quad x \in \closedball_r(\bar{x}),
\end{equation}
where $\hat y^k\in Y$ is the chosen multiplier estimate and $\mu_k\in(0,\mu_g)$.
Clearly, by continuity of the Moreau envelope, see e.g.\ \cite[Thm~1.25]{rockafellar1998variational},
\eqref{eq:subproblem_existence_local_minimizers} possesses a global minimizer $x^k\in\XX$.
Under suitable assumptions it is possible to show
$\norm{x^k-\bar x}<r$ for sufficiently small $\mu_k$, and the localization in
\eqref{eq:subproblem_existence_local_minimizers} becomes superfluous.
In fact, if $\mu_k\downto 0$, we are in position to verify $x^k\to\bar x$ as desired.

As we will see in the subsequent lemma, strict local minimality of
some feasible point $\bar x\in\XX$ of \eqref{eq:P} serves as a sufficient condition 
for a sequence of asymptotically feasible points to converge to $\bar{x}$,
see \cite[Cor.\ 6.2]{BoergensKanzowSteck2019} for a related result under stronger assumptions.

\begin{mybox}
	\begin{lemma}\label{corollary327}
		Let $\bar{x} \in \XX$ be a strict local minimizer of \eqref{eq:P}.
		Then there exists $r > 0$ such that, 
		whenever $\{x^k\} \subseteq \closedball_r(\bar{x})$ and $\{z^k\}\subseteq \dom g$ are sequences 
		with $\norm{ c(x^k) - z^k } \to 0$ and $\limsup_{k \to \infty} (f(x^k) + g(z^k)) \leq \varphi(\bar{x})$, 
		then $x^k \to \bar{x}$ and $z^k \to c(\bar{x})$.
	\end{lemma}
\end{mybox}
\begin{proof}
	The stated assumptions guarantee the existence of $r>0$ such that
	\begin{equation}\label{eq:strict_local_minimality}
		\forall x\in\closedball_r(\bar x)\setminus\{\bar x\}\colon\quad
		\varphi(x)>\varphi(\bar x).
	\end{equation}
	Let us pick arbitrary sequences $\{x^k\}$ and $\{z^k\}$ satisfying the requirements.
	Suppose now that $x^k\nto\bar x$.
	Hence, as $\{x^k\}$ belongs to the compact set $\closedball_r(\bar x)$, 
	it possesses an accumulation point $\tilde x\in\closedball_r(\bar x)$ such that $\tilde x\neq\bar x$.
	Let $K\subseteq\N$ be an index set such that $x^k\to_K\tilde x$.
	Continuity of $c$ yields $c(x^k)\to_Kc(\tilde x)$, and $\norm{c(x^k)-z^k}\to 0$ gives $z^k\to_Kc(\tilde x)$.
	Furthermore, continuity of $f$ and lower semicontinuity of $g$ can be used to infer
	\begin{align*}
		\varphi(\bar x)
		&\geq
		\limsup\limits_{k\to\infty}(f(x^k)+g(z^k))
		\geq
		\liminf\limits_{k\to_K\infty}(f(x^k)+g(z^k))
		\geq
		\varphi(\tilde x),
	\end{align*}
	but this contradicts \eqref{eq:strict_local_minimality}. Thus, it must be $x^k\to\bar x$.
	Then, continuity of $c$ yields $c(x^k)\to c(\bar x)$, and $z^k\to c(\bar x)$
	follows from $\norm{c(x^k)-z^k}\to 0$.
\end{proof}

Observe that SOSC from \cref{def:SOSC} provides a sufficient
condition for strict local minimality
which can be checked in terms of initial problem data,
see \cref{prop:SOSC_P}.
It is remarkable that,
in analogous ways, one can show that
\cite[Cor.\ 6.2]{BoergensKanzowSteck2019} 
(stated in an infinite-dimensional setting) remains true whenever the 
considered point of interest therein is supposed to be a strict local minimizer.
One does not need to assume validity of a second-order sufficient condition
for that purpose.

The following is an analog of \cite[Lem.\ 7.1]{BoergensKanzowSteck2019}.

\begin{mybox}
	\begin{lemma}\label{lem:uniform_convergence_minimizers_subproblems}
		Let $\bar{x}\in\XX$ be a strict local minimizer of \eqref{eq:P}.
		Furthermore, let $\Ybounded\subseteq\YY$ be bounded.
		Then there is a radius $r > 0$ such that whenever
		$\{\hat{y}^k\} \subseteq \Ybounded$,
		$\mu_k\downto 0$, $\varepsilon_k \downtoeq 0$, and, 
		for all $k\in\N$, $x^k$ is an $\varepsilon_k$-minimizer 
		of \eqref{eq:subproblem_existence_local_minimizers} in the sense that
		\begin{equation}\label{eq:approx_minimality_AL}
			\forall x\in\closedball_r(\bar x)\colon\quad
			\LL_{\mu_k}(x^k,\hat y^k)\leq\LL_{\mu_k}(x,\hat y^k)+\varepsilon_k,
		\end{equation}
		then $x^k \to \bar{x}$.
	\end{lemma}
\end{mybox}
\begin{proof}
	Let $r > 0$ be as in \cref{corollary327}.
	For large enough $k\in\N$, $\mu_k\in(0,\mu_g)$ holds, and we can
	fix $z^k \in \prox_{\mu_k g}(c(x^k) + \mu_k \hat{y}^k)$.
	For any such $k\in\N$, \eqref{eq:approx_minimality_AL} and \cref{lem:upper_bound_AL} yield
	\begin{align*}		
		f(x^k)
		&
		+ 
		g(z^k) + \frac{1}{2 \mu_k} \| c(x^k) + \mu_k \hat{y}^k - z^k \|^2
		- \frac{\mu_k}{2} \| \hat{y}^k \|^2 
		\\
		&
		=
		f(x^k)
		+ g^{\mu_k}\left( c(x^k) + \mu_k \hat{y}^k \right)
		- \frac{\mu_k}{2} \| \hat{y}^k \|^2 \\
		&
		=
		\LL_{\mu_k}(x^k,\hat{y}^k)
		\leq 
		\LL_{\mu_k}(\bar{x},\hat{y}^k) + \varepsilon_k
		\leq
		\varphi(\bar{x}) + \varepsilon_k
		<
		\infty.
	\end{align*}
	Multiplying by $\mu_k$, by the boundedness of $\{\hat{y}^k\}$ and $\{f(x^k)\}$, $\mu_k\downto 0$, and $\varepsilon_k\downtoeq 0$,
	we obtain
	\[
		\limsup\limits_{k\to\infty}
		\left(
			\mu_k g(z^k) + \frac12\norm{z^k-c(x^k)-\mu_k\hat y^k}^2
		\right)
		\leq 
		0.
	\]
	We apply \cref{lem:norm_to_zero} to find $\norm{z^k-c(x^k)-\mu_k\hat y^k}\to 0$ and, thus,
	$\norm{ c(x^k) - z^k } \to 0$.
	Moreover, the above estimate also guarantees $\limsup_{k \to \infty} ( f(x^k) + g(z^k) ) \leq \varphi(\bar{x})$,
	again by boundedness of $\{\hat y^k\}$ and $\mu_k\downto 0$.
	Hence, \cref{corollary327} is applicable and yields the desired convergence.
\end{proof}

As a consequence of \cref{lem:uniform_convergence_minimizers_subproblems}, 
we find the following result which parallels \cite[Thm~7.2]{BoergensKanzowSteck2019}.

\begin{mybox}
\begin{theorem}\label{lem:localMin:convergence_to_strict_local_min}
	Let $\bar{x}\in\XX$ be a strict local minimizer of \eqref{eq:P}.
	Furthermore, let $\Ybounded\subseteq\YY$ be bounded.
	Then there is a radius $ r > 0$ such that, for every $\hat{y} \in \Ybounded$ and $\mu \in (0,\mu_g)$, the function $\LL_\mu(\cdot,\hat{y})$ 
	admits a local minimizer $x(\mu,\hat{y})$ which lies in $\closedball_r(\bar{x})$.
	Moreover, $x(\mu,\hat y) \to \bar{x}$ uniformly on $\Ybounded$ as $\mu \downto 0$. 
\end{theorem}
\end{mybox}
\begin{proof}
	Let $r>0$ be as in \cref{lem:uniform_convergence_minimizers_subproblems}.
	For $\mu\in(0,\mu_g)$,
	the Moreau envelope $g^{\mu}$ is a continuous function, 
	and this extends to $\LL_\mu(\cdot,\hat y)$ for arbitrary $\hat y\in\Ybounded$.
	Hence, this function possesses a global minimizer over $\closedball_r(\bar x)$ which we denote by $x(\mu,\hat y)$.
	As $\mu\downto 0$, we find $x(\mu,\hat y)\to\bar x$ from \cref{lem:uniform_convergence_minimizers_subproblems},
	and this convergence is uniform for $\hat y\in\Ybounded$.
\end{proof}

Let us note that \cref{lem:uniform_convergence_minimizers_subproblems} and \cref{lem:localMin:convergence_to_strict_local_min}
merely assume strict local minimality of the reference point.
As mentioned before, this holds true whenever SOSC in valid at the point of interest.
Note that SOSC does not demand uniqueness of the underlying Lagrange multiplier.

\subsubsection{Rates of convergence}

Our rates-of-convergence analysis of \cref{alg:ALM} is based on a primal-dual pair $(\bar x,\bar y)\in\XX\times\YY$
which solves the M-stationarity system \eqref{eq:Mstationary} associated with \eqref{eq:P} such that
an (upper) error bound condition is valid. 

For brevity of notation, we will partially make use of the following assumptions.
\begin{mybox}
	\begin{assum}[Rates of convergence]\label{ass:rate_convergence}
		\phantom{The following hold.}
		\begin{enumerate}[label=(\roman*)]
			\item\label{item:error_bound}%
			Let $\bar x\in \XX$ be an M-stationary point of \eqref{eq:P}, and let $\bar y\in\Lambda(\bar x)$
			be chosen such that there are a constant $\varrho_{\mathrm{u}}>0$ and
			a neighborhood $U$ of $(\bar x,c(\bar x),\bar y)$ such that
			the upper error bound condition \eqref{eq:error_bound} holds
			for each triplet $(x,z,y)\in U\cap(\XX\times\dom g\times\YY)$.
			\item\label{item:ALMsequence}%
			Let $\{(x^k,z^k,y^k)\}$ be a sequence generated by \cref{alg:ALM} with $\varepsilon_k\downtoeq 0$.
			\item\label{item:pure_convergence}%
			The primal-dual sequence $\{(x^k,y^k)\}$ converges to $(\bar{x},\bar{y})$.
			\item\label{item:no_safeguarding}%
			For each $k\in\N$ large enough, $\hat{y}^k = y^{k-1}$ is valid.
		\end{enumerate}
	\end{assum}
\end{mybox}

Note that we already know, by \cref{prop:SOSC_P} and \cref{lem:localMin:convergence_to_strict_local_min}, 
that the AL admits approximate local minimizers and stationary points in a neighborhood of
some M-stationary point $\bar{x}\in\XX$ which satisfies SOSC.
We shall now see that,
under the error bound conditions from \cref{sec:ErrorBounds} involving a fixed multiplier $\bar y\in\Lambda(\bar x)$,
if the algorithm chooses these local minimizers (or any other points sufficiently close to $\bar{x}$), 
then we automatically obtain the convergence $(x^k,z^k,y^k) \to (\bar{x},c(\bar x),\bar{y})$.
In this case, the sequence $\{y^k\}$ is necessarily bounded, 
so it is reasonable to assume that the safeguarded multipliers are eventually chosen as $\hat{y}^k = y^{k-1}$.

Let us recall that even validity of 
the second-order condition \eqref{eq:SSOSC}, which is more restrictive than SOSC,
may not be sufficient for the error bound condition, see \cref{rem:on_the_really_crucial_CQ}.
However, \cref{sec:ErrorBounds} provides a number of sufficient conditions which still can be checked
in terms of initial problem data, so we will abstain here from postulating any more specific
assumptions on the upper error bound.
Furthermore, we do not stipulate any lower error bound conditions,
deviating from all other related papers,
where the lower estimate was never problematic,
see \cref{rem:lower_error_bound}.

The following result, which is motivated by \cite[Prop.\ 4.29]{steck2018dissertation}, 
can be considered as (retrospective) justification for
\cref{ass:rate_convergence}\ref{item:pure_convergence}--\ref{item:no_safeguarding}
in the presence of \cref{ass:rate_convergence}\ref{item:error_bound}--\ref{item:ALMsequence}.
Besides the error bound condition, a CQ is needed.
As we require an M-stationary point of \eqref{eq:P}, this is not too restrictive.

\begin{mybox}
	\begin{proposition}\label{lem:ErrorBound:convergence_vanishing_error}
		Let \cref{ass:rate_convergence}\ref{item:error_bound}--\ref{item:ALMsequence} hold
		and suppose that (at least) one of the following conditions is valid:
		\begin{enumerate}
			\item\label{item:ErrorBound:convergence_vanishing_error:fullrank}%
			$c^\prime(\bar x)$ possesses full row rank $m$;
			\item\label{item:ErrorBound:convergence_vanishing_error:some_CQ}%
			condition \eqref{eq:metric_regularity_CQ} is valid, $\dom g$ is closed, and $g$ is continuous relative to its domain.
		\end{enumerate}
		Then there exists a radius $r > 0$ such that,
		if $x^k \in \closedball_r(\bar{x})$ for all sufficiently large $k\in\N$, 
		then we have the convergences $\Theta(x^k,z^k,y^k) \to 0$ and $(x^k,z^k,y^k) \to (\bar{x},c(\bar{x}),\bar{y})$ as $k\to\infty$.
	\end{proposition}
\end{mybox}
\begin{proof}
	Let $r > 0$ be small enough so that \eqref{eq:error_bound} holds 
	for all $(x,z,y) \in \XX \times \dom g \times \YY$ with $x \in \closedball_r(\bar{x})$.
	Assume now that $x^k \in \closedball_r(\bar{x})$ for all $k\in\N$ sufficiently large.
	The proof is divided into multiple steps.
	
	We first show that $c(x^k) - z^k \to 0$ as $k\to \infty$.
	Consider two cases.
	If $\{\mu_k\}$ remains bounded away from zero, this assertion readily follows from the penalty updating scheme at \cref{step:ALM:penalty}.
	Instead, if $\mu_k\downto 0$, then we can argue from
	\cref{lem:ALM:exact_complementarity} that
	\begin{equation}\label{eq:some_intermediate_convergence}
		c^\prime(x^k)^\top \bigl[ c(x^k) + \mu_k\hat y^k - z^k \bigr] \to 0
	\end{equation}
	as $\mu_k\downto 0$, by boundedness of $\{\nabla f(x^k)\}$, $\{\hat{y}^k\}$, and $\{\varepsilon_k\}$.
	Let us now show $c(x^k) + \mu_k\hat y^k-z^k \to 0$,
	which readily yields $c(x^k)-z^k\to 0$ since $\{\hat y^k\}$ is bounded and $\mu_k\downto 0$.
	In case \ref{item:ErrorBound:convergence_vanishing_error:fullrank} where $c^\prime(\bar x)$ has full row rank,
	the matrices $c^\prime(x) c^\prime(x)^\top$ are uniformly invertible
	on $\closedball_r(\bar x)$, potentially after shrinking $r$, and \eqref{eq:some_intermediate_convergence}
	gives $c(x^k) + \mu_k\hat y^k-z^k \to 0$.
	Next, for case \ref{item:ErrorBound:convergence_vanishing_error:some_CQ},
	assume that \eqref{eq:metric_regularity_CQ} holds while $\dom g$ is closed and $g$ is continuous on its domain.
	Note that, for each $k\in\N$, we even have $\mu_k^{-1}(c(x^k)+\mu_k\hat y^k-z^k)\in\regsub g(z^k)$ by
	definition of the prox-operator and compatibility of the regular subdifferential with smooth additions, and
	this also gives $(c(x^k)+\mu_k\hat y^k-z^k,-\mu_k)\in\widehat{N}_{\epi g}(z^k,g(z^k))$.
	Recall that \eqref{eq:metric_regularity_CQ} is equivalent to the metric regularity of
	$\Xi$ from \eqref{eq:svm_composition} at $((\bar x,g(c(\bar x))),(0,0))$.
	\Cref{lem:some_stability_of_MR} now yields the existence of $s>0$ such that, for
	all sufficiently large $k\in\N$, we have
	\begin{equation}\label{eq:uniform_condition}
			\closedball_s(0,0)
			\subseteq
			\begin{bmatrix}
				c^\prime(x^k) 	&	0
				\\
				0			&	1
			\end{bmatrix}
			\closedball_1(0,0)
			-
			\bigl(T_{\epi g}(z^k,g(z^k))
			\cap
			\closedball_1(0,0)\bigr)
	\end{equation}
	as $g$ is continuous on $\dom g$.
	In order to see this, we need to make sure that
	$(z^k,g(z^k))$ is sufficiently close to $(c(\bar x),g(c(\bar x)))$
	for large enough $k\in\N$, and due to the continuity of $g$,
	this boils down to showing
	that $z^k$ is sufficiently close to $c(\bar x)$ for large enough $k\in\N$.
	
	Along the tail of the sequence (without relabeling), 
	we have that $\{ x^k \}$ is close to $\bar{x}$.
	For every $k$, the optimality of $z^k$ in the proximal minimization subproblem reads
	\begin{equation*}
		\forall z \in \YY \colon\qquad
		\mu_k g(z^k) + \frac{1}{2} \norm{ c(x^k) + \mu_k \hat{y}^k - z^k }^2
		\leq
		\mu_k g(z) + \frac{1}{2} \norm{ c(x^k) + \mu_k \hat{y}^k - z }^2
		.
	\end{equation*}
	Taking the specific choice $z \coloneqq c(\bar{x}) \in \dom g$ and dividing both sides by $\mu_k > 0$ results in
	\begin{equation*}
		g(z^k) + \frac{1}{2 \mu_k} \norm{ c(x^k) + \mu_k \hat{y}^k - z^k }^2
		\leq
		g(c(\bar{x})) + \frac{1}{2 \mu_k} \norm{ c(x^k) + \mu_k \hat{y}^k - c(\bar{x}) }^2
		<
		\infty
		.
	\end{equation*}
	By invoking the triangle, Cauchy--Schwarz, and Young's inequalities, this implies that
	\begin{align*}
		\norm{ z^k -  c(\bar{x}) }^2
		{}={}&
		\norm{ z^k - [c(x^k) + \mu_k \hat{y}^k] - c(\bar{x}) + [c(x^k) + \mu_k \hat{y}^k] }^2 \\
		{}\leq{}&
		\norm{ z^k - [c(x^k) + \mu_k \hat{y}^k] }^2
		+ \norm{ c(\bar{x}) -  [c(x^k) + \mu_k \hat{y}^k] }^2 \\
		{}{}&\qquad
		+ 2 \norm{ z^k - [c(x^k) + \mu_k \hat{y}^k] } \norm{ c(\bar{x}) - [c(x^k) + \mu_k \hat{y}^k] } \\
		{}\leq{}&
		2 \norm{ z^k - [c(x^k) + \mu_k \hat{y}^k] }^2 + 2 \norm{ c(\bar{x}) - [c(x^k) + \mu_k \hat{y}^k] }^2 \\
		{}\leq{}&
		4 \left[ \mu_k g( c(\bar{x}) ) - \mu_k g(z^k) + \norm{ c(x^k) + \mu_k \hat{y}^k - c(\bar{x}) }^2 \right]
		.
	\end{align*}
	Rearranging gives
	\begin{equation*}
		\mu_k g(z^k) + \frac{1}{4} \norm{ z^k -  c(\bar{x}) }^2
		\leq
		\mu_k g( c(\bar{x}) ) + \norm{ c(x^k) + \mu_k \hat{y}^k - c(\bar{x}) }^2.
	\end{equation*}
	Since $c(\bar{x}) \in \dom g$, the term $\mu_k g( c(\bar{x}) )$ vanishes as $\mu_k\downto 0$,
	and so does $\mu_k \hat{y}^k$.
	Therefore,
	possibly shrinking the neighborhood considered around $\bar{x}$,
	the right-hand side remains bounded by some arbitrarily small $C > 0$ for all large $k\in\N$,
	by continuous differentiability of $c$, i.e.,
	\begin{equation*}
		\mu_k g(z^k) + \frac{1}{4} \norm{ z^k -  c(\bar{x}) }^2
		\leq
		C
	\end{equation*}
	holds for all $k\in\N$ large enough.
	By virtue of the prox-boundedness of $g$, 
	\cite[Lem.\ 4.1]{demarchi2022proximal} yields boundedness of $\{z^k\}$
	and, thus, of $\{g(z^k)\}$ by continuity of $g$ on its domain which is assumed to be closed 
	(Heine's theorem yields uniform continuity of $g$ on closed, bounded subsets of $\dom g$).
	As $\mu_k\downto 0$ and $C>0$ can be made arbitrarily small if only $r>0$ is chosen small enough,
	it follows that $\{z^k\}$ is arbitrarily close to $c(\bar x)$ for all large enough $k\in\N$.
	
	Pick $w\in \closedball_s(0)$ arbitrary. 
	Then, for each sufficiently large $k\in\N$, we can rely on \eqref{eq:uniform_condition}
	to find $(u^k,\alpha_k)\in \closedball_1(0,0)$ and 
	$(v^k,\beta_k)\in T_{\epi g}(z^k,g(z^k))\cap\closedball_1(0,0)$ such that
	$(w,0)=(c^\prime(x^k)u^k-v^k, \alpha_k-\beta_k)$.
	From $(v^k,\beta_k)\in T_{\epi g}(z^k,g(z^k))$, we find
	$\innprod{(v^k,\beta_k)}{(c(x^k)+\mu_k\hat y^k-z^k,-\mu_k)}\leq 0$
	due to \eqref{eq:polarization_rules}.
	Thus
	\begin{align*}
		\innprod{w}{c(x^k)+\mu_k\hat y^k-z^k}
		={}&
		\innprod{c^\prime(x^k)u^k-v^k}{c(x^k)+\mu_k\hat y^k-z^k}
		\\
		={}&
		\innprod{u^k}{c^\prime(x^k)^\top(c(x^k)+\mu_k\hat y^k-z^k)}
		-
		\mu_k\beta_k
		-
		\innprod{(v^k,\beta_k)}{(c(x^k)+\mu_k\hat y^k-z^k,-\mu_k)}
		\\
		\geq{}&
		\innprod{u^k}{c^\prime(x^k)^\top(c(x^k)+\mu_k\hat y^k-z^k)}
		-
		\mu_k\beta_k
		\to 0,
	\end{align*}
	where we used boundedness of $\{u^k\}$ and $\{\beta_k\}$ 
	as well as $\mu_k\downto 0$ and \eqref{eq:some_intermediate_convergence}.
	Testing this expression with
	$w \coloneqq \pm s(c(x^k)+\mu_k\hat y^k-z^k)/\norm{c(x^k)+\mu_k\hat y^k-z^k}$
	gives
	$c(x^k)+\mu_k\hat y^k-z^k\to 0$.
	
	We now demonstrate that $\Theta(x^k,z^k,y^k) \to 0$ as $k\to\infty$.
	Observe that 
	\[
		\nabla_x\LL(x^k,y^k)
		=
		\nabla_x \LLslack(x^k,z^k,y^k) 
		= 
		\nabla_x \LLslack_{\mu_k}(x^k,z^k,\hat{y}^k)
	\] 
	holds for all $k\in\N$ by construction of the dual update rule in \cref{step:ALM:y}.
	Then the first summand in $\Theta(x^k,z^k,y^k)$ satisfies
	\begin{equation*}
		\norm{ \nabla_x \LL(x^k,y^k) }
		=
		\norm{ \nabla_x \LLslack_{\mu_k}(x^k,z^k,\hat{y}^k) }
		\leq
		\varepsilon_k ,
	\end{equation*}
	which converges to zero by \cref{ass:rate_convergence}\ref{item:ALMsequence}.
	Hence, as $\norm{ c(x^k) - z^k } \to 0$ was obtained previously, the second term in \eqref{eq:ErrorBound:error} vanishes, too.
	For the third and last term, it remains to show that $\dist(y^k, \partial g(z^k))\to 0$.
	This, however, readily follows from \cref{lem:ALM:exact_complementarity}.
	
	Finally, recall that $x^k \in \closedball_r(\bar{x})$ for all $k\in\N$ and that $\Theta(x^k,z^k,y^k) \to 0$.
	Hence, the convergence $(x^k,z^k,y^k) \to (\bar{x},c(\bar x),\bar{y})$ 
	is an immediate consequence of \eqref{eq:error_bound}.
\end{proof}

Subsequently, we will prove convergence rates for the sequence $\{(x^k,z^k,y^k)\}$ in the presence of \cref{ass:rate_convergence}.
Since the distance of $(x^k,z^k,y^k)$ to $(\bar{x},c(\bar x),\bar{y})$ admits an estimate 
relative to the residual terms $\Theta_k \coloneqq \Theta(x^k,z^k,y^k)$ by \eqref{eq:error_bound}, 
we will largely base our analysis on the sequence $\{\Theta_k\}$,
and the results on the sequence $\{(x^k,z^k,y^k)\}$ will follow directly.
However, this correspondence heavily relies on a \emph{two-sided} error bound,
see the proof of \cref{thm:rate_convergence} below.
In stark contrast to \cref{rem:lower_error_bound},
the following \cref{lem:lower_error_bound:ALM} shows that,
along a sequence generated by \cref{alg:ALM},
a lower error bound holds.
This exploits the fact that,
as a consequence of \cref{lem:ALM:exact_complementarity},
the distance-to-subdifferential in $\Theta$ does not play a role
for the error bound \emph{at the iterates}.
Therefore,
complementing the upper estimate of \cref{ass:rate_convergence}\ref{item:error_bound},
a two-sided error bound becomes \emph{algorithmically} available,
enabling the derivation of convergence rates.
\begin{mybox}
	\begin{lemma}\label{lem:lower_error_bound:ALM}
		Let $\bar x\in\XX$ be an M-stationary point of \eqref{eq:P} and $\bar y\in\Lambda(\bar x)$ be arbitrary.
		Suppose \cref{ass:rate_convergence}\ref{item:ALMsequence}--\ref{item:no_safeguarding} hold.
		Then there are a constant $\varrho_{\mathrm{l}}>0$ and
		a neighborhood $U$ of $(\bar x,c(\bar x),\bar y)$
		such that,
		for each triplet $(x^k,z^k,y^k)\in U\cap(\XX\times\dom g\times\YY)$,
		we have
		\begin{equation}\label{eq:error_bound_lower:ALM}
		\varrho_{\mathrm{l}}\,\Theta(x^k,z^k,y^k)
		\leq
		\norm{x^k-\bar x}+\norm{z^k-c(\bar x)}+\norm{y^k-\bar y}.
		\end{equation}
	\end{lemma}
\end{mybox}
\begin{proof}
	For each triplet $(x^k,z^k,y^k)\in\XX\times\dom g\times\YY$, 
	we can exploit \eqref{eq:Mstationary} and the triangle inequality to obtain
	\begin{align*}
	\Theta(x^k,z^k,y^k)
	{}={}&
	\norm{\nabla_x\LL(x^k,y^k)} + \norm{c(x^k)-z^k}
	\\
	{}\leq{}&
	\norm{\nabla_x\LL(x^k,y^k)-\nabla_x\LL(\bar x,\bar y)}
	+
	\norm{c(x^k)-c(\bar x)} + \norm{z^k-c(\bar x)}
	,
	\end{align*}
	where the equality is due to
	\cref{ass:rate_convergence}\ref{item:ALMsequence}
	and \cref{lem:ALM:exact_complementarity},
	which imply $\dist(y^k,\partial g(z^k)) = 0$ for all $k\in\N$.
	Then, noting that $\nabla_x\LL$ and $c$ are
	locally Lipschitz continuous,
	the claim follows.
\end{proof}

Our next result,
preparatory for \cref{thm:rate_convergence} below,
has been inspired by \cite[Lem.\ 4.30]{steck2018dissertation}.

\begin{mybox}
	\begin{lemma}\label{lem:ErrorBound:convergence_rate_error}
		Let \cref{ass:rate_convergence} hold and set $\Theta_k \coloneqq \Theta(x^k,z^k,y^k)$ for each $k\in\N$.
		Then
		\begin{equation*}
			( 1 - \varrho_{\mathrm{u}} \mu_k ) \Theta_k \leq \varepsilon_k + \varrho_{\mathrm{u}} \mu_k \Theta_{k-1} 
		\end{equation*}
		for all $k\in\N$ large enough,
		where $\varrho_{\mathrm{u}}>0$ is the constant from \eqref{eq:error_bound}.
	\end{lemma}
\end{mybox}
\begin{proof}
	Due to \cref{lem:ALM:exact_complementarity}, $y^k \in \partial g(z^k)$ holds for all $k\in\N$.
	Then, by \eqref{eq:ErrorBound:error} and \cref{step:ALM:subproblem}, we have
	\begin{equation*}
		\Theta_k
		\leq
		\varepsilon_k + \norm{ c(x^k) - z^k }
		=
		\varepsilon_k + \mu_k \norm{ y^k - y^{k-1} }
		\leq
		\varepsilon_k + \mu_k \norm{ y^k - \bar{y} } + \mu_k \norm{ y^{k-1} - \bar{y} }
		,
	\end{equation*}
	where the equality is due to the update rule at \cref{step:ALM:y} and \cref{ass:rate_convergence}\ref{item:no_safeguarding}.
	By \cref{ass:rate_convergence}\ref{item:error_bound}, 
	since $x^k \to \bar{x}$ and, due to \cref{lem:ErrorBound:convergence_vanishing_error}, $z^k \to c(\bar{x})$, 
	we find $\norm{y^k - \bar{y}} \leq \varrho_{\mathrm{u}} \Theta_k$ for all $k\in\N$ large enough.
	Hence,
	$\Theta_k\leq\varepsilon_k + \varrho_{\mathrm{u}} \mu_k (\Theta_k + \Theta_{k-1})$
	holds for all $k\in\N$ large enough,
	and reordering gives the assertion.
\end{proof}

With the above lemma and the two-sided error bound enabled by
\cref{ass:rate_convergence}\ref{item:error_bound} and \cref{lem:lower_error_bound:ALM},
one can deduce convergence rates for the sequence $\{(x^k,z^k,y^k)\}$,
see \cite[Thm~4.31]{steck2018dissertation} as well.
Notice that the condition $\varepsilon_k \in \oo(\Theta_{k-1})$ can be easily guaranteed in practice.
For instance, one could compute
the next iterate $(x^k,z^k,y^k)$ with a precision $\varepsilon_k \leq \nu_k \Theta_{k-1}$ 
where $\{\nu_k\}$ is a given null sequence.
It should be mentioned also that the value $\Theta_k$ from \eqref{eq:ErrorBound:error} 
becomes algorithmically available thanks to \cref{lem:ALM:exact_complementarity}, 
and can readily be obtained by the dual update rule at \cref{step:ALM:y}.

\begin{mybox}
	\begin{theorem}\label{thm:rate_convergence}
		Let \cref{ass:rate_convergence} hold and
		assume that $\varepsilon_k \in \oo(\Theta_{k-1})$,
		where $\Theta_k\coloneqq\Theta(x^k,z^k,y^k)$ for each $k\in\N$.
		Then the following assertions hold.
		\begin{enumerate}
			\item For every $q\in(0,1)$, there exists $\bar{\mu}(q)$ such that, if $\mu_k \leq \bar{\mu}(q)$ for sufficiently large $k\in\N$, 
				then $(x^k,z^k,y^k) \to (\bar{x},c(\bar x),\bar{y})$ Q-linearly with rate $q$.
			\item\label{thm:rate_convergence:super}%
			If $\mu_k\downto 0$, then $(x^k,z^k,y^k) \to (\bar{x},c(\bar x),\bar{y})$ Q-superlinearly.
		\end{enumerate}
	\end{theorem}
\end{mybox}
\begin{proof}
	Let $k\in\N$ be large enough so that $\hat{y}^k = y^{k-1}$.
	By \cref{lem:ErrorBound:convergence_rate_error}, if $\mu_k$ is small enough so that $1 - \varrho_{\mathrm{u}} \mu_k > 0$, then
	\begin{equation*}
		\frac{\Theta_k}{\Theta_{k-1}}
		\leq
		\frac{\varrho_{\mathrm{u}} \mu_k}{1 - \varrho_{\mathrm{u}} \mu_k} + \oo(1).
	\end{equation*}
	The desired rates for $\{(x^k,z^k,y^k)\}$ are an easy consequence of
	the upper and lower estimates in \eqref{eq:error_bound}
	and \cref{lem:lower_error_bound:ALM},
	as these give
	\begin{align*}
		\frac{\norm{x^k-\bar x}+\norm{z^k-c(\bar x)}+\norm{y^k-\bar y}}{\norm{x^{k-1}-\bar x}+\norm{z^{k-1}-c(\bar x)}+\norm{y^{k-1}-\bar y}}
		\leq
		\frac{\varrho_{\mathrm{u}}}{\varrho_{\mathrm{l}}}
		\frac{\varrho_{\mathrm{u}} \mu_k}{1 - \varrho_{\mathrm{u}} \mu_k} + \oo(1)
	\end{align*}
	for all $k\in\N$ large enough.
\end{proof}

The following result, analogous to
\cite[Cor.\ 4.32]{steck2018dissertation},
establishes the boundedness of $\{\mu_k\}$ away from zero
in the case of exact subproblem solutions,
thus preventing the fast local convergence of \cref{thm:rate_convergence}\cref{thm:rate_convergence:super}.
\begin{mybox}
	\begin{corollary}
		Let \cref{ass:rate_convergence} hold and assume that the subproblems occurring at \cref{step:ALM:subproblem} of \cref{alg:ALM} are solved exactly, i.e., that $\varepsilon_k = 0$ for all $k\in\N$.
		Then $\{\mu_k\}$ remains bounded away from zero.
	\end{corollary}
\end{mybox}
\begin{proof}
	For each $k\in\N$, we make use of $V_k\coloneqq\norm{c(x^k)-z^k}$ and $\Theta_k\coloneqq\Theta(x^k,z^k,y^k)$.
	Let $k\in\N$ be large enough so that $\hat{y}^k = y^{k-1}$.
	Arguing as in the proof of \cref{lem:ErrorBound:convergence_rate_error}, we have for all $k\in\N$ that $\Theta_k \leq V_k$ since $\varepsilon_k = 0$.
	Furthermore, using the triangle inequality, the convergences, $x^k\to\bar{x}$, $z^k \to c(\bar{x})$, and
	\cref{ass:rate_convergence}\ref{item:no_safeguarding}, we obtain 
	$V_k=\mu_k\norm{y^k-y^{k-1}}\leq\varrho_{\mathrm{u}} \mu_k (\Theta_k + \Theta_{k-1})$
	from \eqref{eq:error_bound}.
	Combining these inequalities yields
	\begin{equation*}
		\frac{V_k}{V_{k-1}}
		\leq
		\varrho_{\mathrm{u}} \mu_k \left(\frac{\Theta_k}{\Theta_{k-1}} + 1 \right)
		.
	\end{equation*}
	Finally, assuming that $\mu_k\downto 0$, we deduce from the proof of \cref{thm:rate_convergence} 
	that $\Theta_k/\Theta_{k-1} \to 0$, and then $V_k/V_{k-1} \to 0$ follows.
	Hence, $V_k / V_{k-1} \leq \theta$ for all $k\in\N$ sufficiently large, 
	where $\theta\in(0,1)$ is a fixed parameter of \cref{alg:ALM},
	so that \cref{step:ALM:penalty} gives a contradiction, thus proving the assertion.
\end{proof}

In summary,
local fast convergence of \cref{alg:ALM}, even for nonconvex functions $g$,
can be obtained in the presence of suitable second-order conditions
(one ensuring the existence of minimizers of the subproblems and
another one to guarantee validity of an upper error bound)
and a first-order CQ which, in principle, gives us the
full convergence of the primal-dual sequence.

In comparison with the noteworthy results from \cite{FernandezSolodov2012,HangSarabi2021}, these assumptions may seem quite strong.
However, let us mention that in the settings discussed in these papers, the (convex)
function $g$ under consideration is chosen in such a way that the aforementioned two
second-order conditions can already be merged into one, 
see \cref{rem:error_bound_for_convex_piecewise_quadratic_functions}.
Furthermore, it is likely that the additional postulation of a first-order CQ
could be avoided in these papers, too,
since $g$ (or at least its derivative) is convex and/or polyhedral \emph{enough}
while, for convex functions, the proximal operator is well-behaved.
It remains a question for future research whether, for example,
a generalized polyhedral structure of $g$
(where its domain and epigraph are unions of finitely many convex polyhedral sets)
makes the additional assumption of a first-order CQ
superfluous.

\section{Some exemplary settings}\label{sec:examples}

In light of our theoretical findings for the general problem \eqref{eq:P},
this section examines two notable illustrative settings:
sparsity-promoting and complementarity-constrained optimization.

\subsection{Sparsity-promoting optimization}

Here, we take a closer look at the sparsity-promoting optimization problem
\eqref{eq:sparse_programming} which has been already discussed in \cref{ex:conjugates_for_nonconvex_g_messy}.

Let us fix some point $\bar x\in\R^n$.
For $\bar y\in\partial\norm{\cdot}_0(c(\bar x))$, we make use of
\[
	I^{00}(\bar x,\bar y)\coloneqq \{i\in I^0(\bar x)\,|\,\bar y_i=0\},
	\quad
	I^{0\pm}(\bar x,\bar y)\coloneqq \{i\in I^0(\bar x)\,|\,\bar y_i\neq 0\}
\]
where $I^0(\bar x)$ has been defined in \cref{ex:conjugates_for_nonconvex_g_messy}.
With the definition of $I^\pm(\bar x)$ therein,
one obtains 
\begin{equation}\label{eq:tangent_cone_ell0_quasi_norm}
	T_{\gph\partial\norm{\cdot}_0}(c(\bar x),\bar y)
	=
	\left\{
		(v,\eta)\in\R^m\times\R^m\,\middle|\,
		\begin{aligned}
			&\forall i\in I^\pm(\bar x)\colon&&\eta_i=0\\
			&\forall i\in I^{0\pm}(\bar x,\bar y)\colon&&v_i=0\\
			&\forall i\in I^{00}(\bar x,\bar y)\colon&& v_i\eta_i=0
		\end{aligned}
	\right\}.
\end{equation}
This can be used to see that \eqref{eq:CQ_uniqueness_of_multiplier} reduces to the linear independence 
of the family $(\nabla c_i(\bar x))_{i\in I^0(\bar x)}$.

For $u\in\R^n\setminus\{0\}$ and $\eta\in D(\partial \norm{\cdot}_0)(c(\bar x),\bar y)(c^\prime(\bar x)u)$,
we easily find 
\[
	\innprod{\eta}{c^\prime(\bar x)u}
	=
	\sum_{i\in I^{\pm}(\bar x)}\underbrace{\eta_i}_{=0}c^\prime_i(\bar x)u
	+
	\sum_{i\in I^0(\bar x)}\underbrace{\eta_i c^\prime_i(\bar x)u}_{=0}
	=
	0,
\]
so that \eqref{eq:some_new_second_order_condition_for_error_bound} reduces to
\begin{equation}\label{eq:sufficient_condition_error_bound_sparse_programming}
	\forall u\in\{u^\prime\in\R^n\,|\,\forall i\in I^{0\pm}(\bar x,\bar y)\colon\,c^\prime_i(\bar x)u^\prime=0\}\setminus\{0\}\colon\quad
	\nabla^2_{xx}\LL(\bar x,\bar y)[u,u]>0,
\end{equation}
and this corresponds to a classical second-order sufficient condition for the nonlinear program
\[
	\minimize_x
	{}\quad{}
	f(x)
	{}\quad{}
	\stt
	{}\quad{}
	c_i(x)=0\quad i\in I^{0\pm}(\bar x,\bar y).
\]

From \cite[Ex.\ 6.3]{BenkoMehlitz2023}, we find
\[
	\mathcal C(\bar x)
	=
	\{u\in\R^n\,|\,f^\prime(\bar x)u\leq 0,\,\forall i\in I^{0}(\bar x)\colon\,c^\prime_i(\bar x)u=0\}.
\]
However, for $u\in\R^n$ satisfying $c^\prime_i(\bar x)u=0$ for all $i\in I^0(\bar x)$, we already have
\[
	f^\prime(\bar x)u
	=
	-\innprod{\bar y}{c^\prime(\bar x)u}
	=
	-\sum_{i\in I^{\pm}(\bar x)}\underbrace{\bar y_i}_{=0}c^\prime_i(\bar x)u
	-\sum_{i\in I^{0}(\bar x)}\bar y_i\underbrace{c^\prime_i(\bar x)u}_{=0}
	=
	0,	
\]
i.e., $u\in\mathcal C(\bar x)$ due to $\bar y\in\partial\norm{\cdot}(c(\bar x))$,
and this particularly holds for $\bar y\in\Lambda(\bar x)$ which exists whenever $\bar x$ is M-stationary.
In the latter case, we thus obtain the simplified representation
\begin{equation}\label{eq:critical_cone_sparse_programming_simplified}
	\mathcal C(\bar x)=\{u\in\R^n\,|\,\forall i\in I^0(\bar x)\colon\,c^\prime_i(\bar x)u=0\}.
\end{equation}
Noting that $\bar y_i=0$ holds for each $i\in I^{\pm}(\bar x)$, \cite[Ex.\ 6.3]{BenkoMehlitz2023}
shows that SOSC is implied by
\[
	\forall u\in\mathcal C(\bar x)\setminus\{0\},\,\exists y\in\Lambda(\bar x)\colon\quad
	\nabla^2_{xx}\LL(\bar x,y)[u,u]>0,
\]
while 
\[
	\exists\bar y\in\Lambda(\bar x),\,\forall u\in\mathcal C(\bar x)\setminus\{0\}\colon\quad
	\nabla^2_{xx}\LL(\bar x,\bar y)[u,u]>0
\]
is sufficient for \eqref{eq:SSOSC}, and these correspond to certain second-order sufficient
optimality conditions for the optimization problem
\[
	\minimize_x
	{}\quad{}
	f(x)
	{}\quad{}
	\stt
	{}\quad{}
	c_i(x)=0\quad i\in I^{0}(\bar x).
\]
Clearly, due to \eqref{eq:critical_cone_sparse_programming_simplified}, 
both conditions are implied by \eqref{eq:sufficient_condition_error_bound_sparse_programming}.
The following example shows that \eqref{eq:sufficient_condition_error_bound_sparse_programming}
can, indeed, be stronger than \eqref{eq:SSOSC}.

\begin{example}\label{ex:SSOSC_does_not_give_error_bound}
	We consider \eqref{eq:sparse_programming}
	for the functions $\func{f}{\R^2}{\R}$ and $\func{c}{\R^2}{\R^2}$ given by
	\[
		f(x)\coloneqq \frac12(x_1-x_2)^2+x_1-x_2,
		\quad
		c(x)\coloneqq \begin{pmatrix}x_1-x_2\\x_1+x_2\end{pmatrix},
	\]
	and choose $\bar x$ to be the origin in $\R^2$.
	Note that $I^0(\bar x)=\{1,2\}$ and $\Lambda(\bar x)=\{(-1,0)\}$,
	i.e., $\bar y \coloneqq (-1,0)$ is the uniquely determined multiplier in this situation.
	As the critical cone $\mathcal C(\bar x)$ reduces to the origin,
	\eqref{eq:SSOSC} is trivially satisfied.
	
	Observe that we have
	\[
		\nabla_{xx}^2\LL(\bar x,\bar y)
		=
		\begin{pmatrix}
			1	&	-1	\\	-1	&	1
		\end{pmatrix},
	\]
	and for $\bar u\coloneqq(1,1)$ and $\bar\eta\coloneqq(0,0)$, 
	we find $\bar \eta\in D(\partial\norm{\cdot}_0)(c(\bar x),\bar y)(c^\prime(\bar x)\bar u)$
	from \eqref{eq:tangent_cone_ell0_quasi_norm}.
	Furthermore, $\nabla_{xx}^2\LL(\bar x,\bar y)\bar u+c^\prime(\bar x)^\top\bar\eta=0$ is valid.
	Hence, \eqref{eq:a_really_crucial_CQ} does not hold,
	and this also shows that the stronger condition \eqref{eq:some_new_second_order_condition_for_error_bound}
	fails---the latter being equivalent to \eqref{eq:sufficient_condition_error_bound_sparse_programming}
	in the present setting.
\end{example}

\subsection{Complementarity-constrained optimization}

Let $m \coloneqq 2p$ for some $p\in\N$ and consider the special situation $g \coloneqq \indicator_{C_\textup{cc}}$ 
where $C_\textup{cc}\subseteq\R^{2p}$ is given by
\[
	C_\textup{cc} \coloneqq \{z\in\R^{2p}\,|\,\forall i\in\{1,\ldots,p\}\colon\,0\leq z_i\perp z_{p+i}\geq 0\},
\]
i.e., $C_\textup{cc}$ is the standard complementarity set. Problem \eqref{eq:P}, thus, reduces to
\begin{equation}
	\tag{MPCC}\label{eq:MPCC}
		\minimize_x
		{}\quad{}
		f(x)
		{}\quad{}
		\stt
		{}\quad{}
		c(x) \in C_\textup{cc},
\end{equation}
a \emph{mathematical problem with complementarity constraints} (MPCC),
see the classical monographs \cite{LuoPangRalph1996,OutrataKocvaraZowe1998}.
Note that standard inequality and equality constraints can be added without any
difficulty and are omitted here for brevity of presentation.

For a feasible point $\bar x\in\XX$ of \eqref{eq:MPCC}, we make use of the index sets
\begin{align*}
	I^{+0}(\bar x)& \coloneqq \{i\in\{1,\ldots,p\}\,|\,c_i(\bar x)>0,\,c_{p+1}(\bar x)=0\},\\
	I^{0+}(\bar x)& \coloneqq \{i\in\{1,\ldots,p\}\,|\,c_i(\bar x)=0,\,c_{p+i}(\bar x)>0\},\\
	I^{00}(\bar x)& \coloneqq \{i\in\{1,\ldots,p\}\,|\,c_i(\bar x)=0,\,c_{p+i}(\bar x)=0\},
\end{align*}
which provide a disjoint partition of $\{1,\ldots,p\}$.
As we have 
\begin{align*}
	\partial g(c(\bar x))
	=
	\partial^\infty g(c(\bar x))
	&=
	\limnormalcone_{C_\textup{cc}}(c(\bar x))
	=
	\left\{
		y\in\R^{2p}\,\middle|\,
		\begin{aligned}
			&\forall i\in I^{+0}(\bar x)\colon&&y_i=0\\
			&\forall i\in I^{0+}(\bar x)\colon&&y_{p+i}=0\\
			&\forall i\in I^{00}(\bar x)\colon&&(y_i\leq 0\,\land\,y_{p+i}\leq 0)\,\lor\,y_iy_{p+i}=0
		\end{aligned}
	\right\},	
\end{align*}
we can specify the precise meaning of the CQ \eqref{eq:metric_regularity_CQ}.
Note that, as $\indicator_{C_\textup{cc}}$ is continuous on its closed domain $C_\textup{cc}$, 
\eqref{eq:metric_regularity_CQ} can be used in \cref{lem:ErrorBound:convergence_vanishing_error}.
We also note that
\begin{align*}
	\widehat N_{C_\textup{cc}}(c(\bar x))
	=
	\left\{
		y\in\R^{2p}\,\middle|\,
		\begin{aligned}
			&\forall i\in I^{+0}(\bar x)\colon&&y_i=0\\
			&\forall i\in I^{0+}(\bar x)\colon&&y_{p+i}=0\\
			&\forall i\in I^{00}(\bar x)\colon&& y_i\leq 0,\,y_{p+i}\leq 0
		\end{aligned}
	\right\}.
\end{align*}

Let us now assume that $\bar x$ is an M-stationary point of \eqref{eq:MPCC}.
Some calculations show that the associated critical cone is given by
\begin{align*}
	\mathcal C(\bar x)
	=
	\left\{
		u\in\XX\,\middle|\,
		\begin{aligned}
			&&&f^\prime(\bar x)u\leq 0\\
			&\forall i\in I^{+0}(\bar x)\colon&&c^\prime_{p+i}(\bar x)u=0\\
			&\forall i\in I^{0+}(\bar x)\colon&&c^\prime_i(\bar x)u=0\\
			&\forall i\in I^{00}(\bar x)\colon&&0\leq c^\prime_i(\bar x)u\perp c^\prime_{p+i}(\bar x)u\geq 0
		\end{aligned}
	\right\},
\end{align*}
see \cite[Sec.\ 5.1]{BenkoMehlitz2023}.
If $\Lambda(\bar x)\cap\widehat{N}_{C_{\textup{cc}}}(c(\bar x))$ is nonempty,
i.e., if $\bar x$ is so-called strongly stationary,
a simplified representation of the critical cone is available which does not involve the
$\nabla f(\bar x)$ anymore but depends on a multiplier 
$y\in\Lambda(\bar x)\cap\widehat{N}_{C_{\textup{cc}}}(c(\bar x))$ and is given by
\begin{align*}
	\mathcal C(\bar x)
	=
	\left\{
		u\in\XX\,\middle|\,
		\begin{aligned}
			&&&f^\prime(\bar x)u\leq 0\\
			&\forall i\in I^{+0}(\bar x)\colon&&c^\prime_{p+i}(\bar x)u=0\\
			&\forall i\in I^{0+}(\bar x)\colon&&c^\prime_i(\bar x)u=0\\
			&\forall i\in I^{00}_{--}(\bar x,y)\colon&&c^\prime_i(\bar x)u=0,\,c^\prime_{p+i}(\bar x)= 0\\
			&\forall i\in I^{00}_{-0}(\bar x,y)\colon&&c^\prime_i(\bar x)u= 0,\,c^\prime_{p+i}(\bar x)u\geq 0\\
			&\forall i\in I^{00}_{0-}(\bar x,y)\colon&&c^\prime_i(\bar x)u\geq 0,\,c^\prime_{p+i}(\bar x)u=0\\
			&\forall i\in I^{00}_{*}(\bar x,y)\colon&&0\leq c^\prime_i(\bar x)u\perp c^\prime_{p+i}(\bar x)u\geq 0
		\end{aligned}
	\right\},
\end{align*}
see \cite[Lem.\ 4.1]{mehlitz2020}. Here, we used
\begin{align*}
	I^{00}_{--}(\bar x,y)& \coloneqq \{i\in I^{00}(\bar x)\,|\,y_i<0,\,y_{p+i}<0\},&
	I^{00}_{-0}(\bar x,y)& \coloneqq \{i\in I^{00}(\bar x)\,|\,y_i<0,\,y_{p+i}=0\},&\\
	I^{00}_{0-}(\bar x,y)& \coloneqq \{i\in I^{00}(\bar x)\,|\,y_i=0,\,y_{p+i}<0\},&
	I^{00}_{*}(\bar x,y)& \coloneqq \{i\in I^{00}(\bar x)\,|\,y_i=0,\,y_{p+i}=0\},&
\end{align*}
which provide a disjoint partition of $I^{00}(\bar x)$.

Following the arguments provided at the end of \cite[Sec.\ 3.1]{BenkoMehlitz2023}, we find
\begin{align*}
	\Lambda(\bar x,u)
	=
	\left\{y\in\Lambda(\bar x)\,\middle|\,
		\begin{aligned}
			&\forall i\in I^{00}_{+0}(\bar x,u)\colon&&y_i=0\\
			&\forall i\in  I^{00}_{0+}(\bar x,u)\colon&&y_{p+i}=0\\
			&\forall i\in I^{00}_{00}(\bar x,u)\colon&& y_i\leq 0,\,y_{p+1}\leq 0
		\end{aligned}
	\right\}.
\end{align*}
for each $u\in\mathcal C(\bar x)$,
where we made use of a disjoint partition of $I^{00}(\bar x)$ given by
\begin{align*}
	I^{00}_{+0}(\bar x,u)& \coloneqq \{i\in I^{00}(\bar x)\,|\,c^\prime_i(\bar x)u>0,\,c^\prime_{p+i}(\bar x)u=0\},\\
	I^{00}_{0+}(\bar x,u)& \coloneqq \{i\in I^{00}(\bar x)\,|\,c^\prime_i(\bar x)u=0,\,c^\prime_{p+i}(\bar x)u>0\},\\
	I^{00}_{00}(\bar x,u)& \coloneqq \{i\in I^{00}(\bar x)\,|\,c^\prime_i(\bar x)u=0,\,c^\prime_{p+i}(\bar x)u=0\}.
\end{align*}
Thus, due to \cite[Thm~5.4]{BenkoMehlitz2023}, SOSC can be stated in the form
\[
	\forall u\in\mathcal C(\bar x)\setminus\{0\},\,\exists y\in\Lambda(\bar x,u)\colon\quad
	\nabla^2_{xx}\LL(\bar x,y)[u,u]>0.
\]
As shown in the proof of \cref{cor:SSOSC}, any multiplier $\bar y\in\Lambda(\bar x)$ suitable to appear
the second-order condition \eqref{eq:SSOSC} necessarily belongs to $\bigcap_{u\in\mathcal C(\bar x)\setminus\{0\}}\Lambda(\bar x,u)\subset\Lambda(\bar x)$,
and for any such multiplier $\bar y$, $\mathrm d^2\indicator_{C_\textup{cc}}(c(\bar x),\bar y)(c^\prime(\bar x)u)$ vanishes,
see \cite[Lem.\ 3.2, Prop.\ 3.6]{BenkoMehlitz2023}. Hence, \eqref{eq:SSOSC} takes the form
\[
	\exists \bar y\in \bigcap_{u\in\mathcal C(\bar x)\setminus\{0\}}\Lambda(\bar x,u),\,
	\forall u\in\mathcal C(\bar x)\setminus\{0\}\colon\quad
	\nabla^2_{xx}\LL(\bar x,\bar y)[u,u]>0.
\]
It follows from \cite[Lem.\ 3.2]{Gfrerer2014} that this is a less restrictive assumption than
the standard second-order sufficient condition for \eqref{eq:MPCC} which takes the form
\[
	\exists \bar y\in\Lambda(\bar x)\cap\widehat{N}_{C_{\textup{cc}}}(c(\bar x)),\,
	\forall u\in\mathcal C(\bar x)\setminus\{0\}\colon\quad
	\nabla^2_{xx}\LL(\bar x,\bar y)[u,u]>0
\]
and is based on a strongly stationary point.
A detailed study on the relationship between SOSC and \eqref{eq:SSOSC}
as well as other MPCC-tailored second-order optimality conditions
is beyond the scope of this paper,
see e.g.\ \cite{Gfrerer2014,GuoLinYe2013} for an overview.

The graphical derivative of the limiting normal cone mappings associated with $C_\textup{cc}$ has been computed
recently in \cite[Sec.\ 4.4.1]{BenkoMehlitz2024}, and the obtained formulas can be used to specify the CQs
\eqref{eq:strong_metric_subregularity},
\eqref{eq:CQ_uniqueness_of_multiplier},
\eqref{eq:a_really_crucial_CQ},
and
\eqref{eq:some_new_second_order_condition_for_error_bound}
in the recent setting.

\section{Concluding remarks}\label{sec:conclusions}

The results in this paper could be extended to cover the
extra feature in \eqref{eq:P} of a geometric \emph{convex} constraint $x \in X$,
which was not included here for reasons of exposition.
It remains unclear, instead, how to address such additional constraint with nonconvex $X$,
if not reformulating into \eqref{eq:P} and accepting $x\in X$ as a soft constraint,
see \cite[Rem.\ 5.1]{demarchi2024implicit}.

Another challenging question is whether
it is possible to dispose the additional constraint qualification in the nonconvex polyhedral case
(i.e., $\epi g$ being the union of finitely many convex polyhedra) in the analysis
of \cref{sec:ALM}.
Such a result would yield convergence rates merely via some second-order sufficient conditions and the (upper) error bound,
generalizing \cite{FernandezSolodov2012}.
In specific situations, this should be possible even in the nonpolyhedral setting,
as \cite{HangSarabi2021} has shown for the case of convex linear-quadratic $g$.
As already pointed out in \cref{sec:introduction}, however, such a generalization does not seem to be available
in nonpolyhedral settings.

Future research may also
focus on the relationship between the proximal point algorithm and the augmented Lagrangian method
in the fully nonconvex setting, in the vein of \cite{rockafellar1976augmented,rockafellar2022convergence},
and investigate saddle-point properties of the augmented Lagrangian function in primal-dual terms
as in \cite{steck2018dissertation}.

\section*{Acknowledgments}

The authors thank the anonymous reviewer whose valuable 
comments and suggestions helped to strengthen some results,
particularly \cref{corollary327},
and to improve the presentation of this paper.

{%
\phantomsection
\addcontentsline{toc}{section}{References}%
\small
\bibliographystyle{jnsao}
\bibliography{biblio}

\begin{thebibliography}{10}

\bibitem{BankGuddatKlatteKummerTammer1983}
B{.\nobreak\kern 0.33333em}Bank, J{.\nobreak\kern 0.33333em}Guddat,
  D{.\nobreak\kern 0.33333em}Klatte, B{.\nobreak\kern 0.33333em}Kummer, and
  K{.\nobreak\kern 0.33333em}Tammer, \emph{Non-Linear Parametric Optimization},
  Birkh{\"a}user, Basel, 1983,
  \href{https://dx.doi.org/10.1007/978-3-0348-6328-5}{\nolinkurl{doi:10.1007/978-3-0348-6328-5}}.

\bibitem{BauschkeCombettes2011}
H.\,H{.\nobreak\kern 0.33333em}Bauschke and P.\,L{.\nobreak\kern
  0.33333em}Combettes, \emph{Convex Analysis and Monotone Operator Theory in
  {H}ilbert Spaces}, Springer, New York, 2011,
  \href{https://dx.doi.org/10.1007/978-1-4419-9467-7}{\nolinkurl{doi:10.1007/978-1-4419-9467-7}}.

\bibitem{Benko2021}
M{.\nobreak\kern 0.33333em}Benko, On inner calmness*, generalized calculus, and
  derivatives of the normal-cone mapping, \emph{Journal of Nonsmooth Analysis
  and Optimization} 2 (2021),  5881,
  \href{https://dx.doi.org/10.46298/jnsao-2021-5881}{\nolinkurl{doi:10.46298/jnsao-2021-5881}}.

\bibitem{BenkoGfrererOutrata2019}
M{.\nobreak\kern 0.33333em}Benko, H{.\nobreak\kern 0.33333em}Gfrerer, and
  J.\,V{.\nobreak\kern 0.33333em}Outrata, Calculus for directional limiting
  normal cones and subdifferentials, \emph{Set-Valued and Variational Analysis}
  27 (2019),  713--745,
  \href{https://dx.doi.org/10.1007/s11228-018-0492-5}{\nolinkurl{doi:10.1007/s11228-018-0492-5}}.

\bibitem{BenkoGfrererYeZhangZhou2022}
M{.\nobreak\kern 0.33333em}Benko, H{.\nobreak\kern 0.33333em}Gfrerer,
  J.\,J{.\nobreak\kern 0.33333em}Ye, J{.\nobreak\kern 0.33333em}Zhang, and
  J{.\nobreak\kern 0.33333em}Zhou, Second-order optimality conditions for
  general nonconvex optimization problems and variational analysis of
  disjunctive systems, \emph{SIAM Journal on Optimization} 33 (2023),
  2625--2653,
  \href{https://dx.doi.org/10.1137/22M1484742}{\nolinkurl{doi:10.1137/22m1484742}}.

\bibitem{benko2021implicit}
M{.\nobreak\kern 0.33333em}Benko and P{.\nobreak\kern 0.33333em}Mehlitz, On
  implicit variables in optimization theory, \emph{Journal of Nonsmooth
  Analysis and Optimization} 2 (2021),  7215,
  \href{https://dx.doi.org/10.46298/jnsao-2021-7215}{\nolinkurl{doi:10.46298/jnsao-2021-7215}}.

\bibitem{BenkoMehlitz2022}
M{.\nobreak\kern 0.33333em}Benko and P{.\nobreak\kern 0.33333em}Mehlitz,
  Calmness and calculus: two basic patterns, \emph{Set-Valued and Variational
  Analysis} 30 (2022),  81--117,
  \href{https://dx.doi.org/10.1007/s11228-021-00589-x}{\nolinkurl{doi:10.1007/s11228-021-00589-x}}.

\bibitem{BenkoMehlitz2023}
M{.\nobreak\kern 0.33333em}Benko and P{.\nobreak\kern 0.33333em}Mehlitz, Why
  second-order sufficient conditions are, in a way, easy --- or --- revisiting
  calculus for second subderivatives, \emph{Journal of Convex Analysis} 30
  (2023),  541--589,
  \url{https://www.heldermann.de/JCA/JCA30/JCA302/jca30031.htm}.

\bibitem{BenkoMehlitz2024}
M{.\nobreak\kern 0.33333em}Benko and P{.\nobreak\kern 0.33333em}Mehlitz, On the
  directional asymptotic approach in optimization theory, \emph{arXiv}  (2024),
  \href{https://arxiv.org/abs/2402.16530v1}{\nolinkurl{arXiv:2402.16530v1}}.

\bibitem{bertsekas1996constrained}
D.\,P{.\nobreak\kern 0.33333em}Bertsekas, \emph{Constrained Optimization and
  {L}agrange Multiplier Methods}, Athena Scientific, 1996.

\bibitem{birgin2014practical}
E.\,G{.\nobreak\kern 0.33333em}Birgin and J.\,M{.\nobreak\kern
  0.33333em}Mart\'inez, \emph{Practical Augmented {L}agrangian Methods for
  Constrained Optimization}, SIAM, Philadelphia, 2014,
  \href{https://dx.doi.org/10.1137/1.9781611973365}{\nolinkurl{doi:10.1137/1.9781611973365}}.

\bibitem{bolte2018nonconvex}
J{.\nobreak\kern 0.33333em}Bolte, S{.\nobreak\kern 0.33333em}Sabach, and
  M{.\nobreak\kern 0.33333em}Teboulle, Nonconvex {L}agrangian-based
  optimization: monitoring schemes and global convergence, \emph{Mathematics of
  Operations Research} 43 (2018),  1210--1232,
  \href{https://dx.doi.org/10.1287/moor.2017.0900}{\nolinkurl{doi:10.1287/moor.2017.0900}}.

\bibitem{BoergensKanzowSteck2019}
E{.\nobreak\kern 0.33333em}B\"{o}rgens, C{.\nobreak\kern 0.33333em}Kanzow, and
  D{.\nobreak\kern 0.33333em}Steck, Local and global analysis of multiplier
  methods for constrained optimization in {B}anach spaces, \emph{SIAM Journal
  on Control and Optimization} 57 (2019),  3694--3722,
  \href{https://dx.doi.org/10.1137/19M1240186}{\nolinkurl{doi:10.1137/19m1240186}}.

\bibitem{chen2017augmented}
X{.\nobreak\kern 0.33333em}Chen, L{.\nobreak\kern 0.33333em}Guo,
  Z{.\nobreak\kern 0.33333em}Lu, and J.\,J{.\nobreak\kern 0.33333em}Ye, An
  augmented {L}agrangian method for non-{L}ipschitz nonconvex programming,
  \emph{SIAM Journal on Numerical Analysis} 55 (2017),  168--193,
  \href{https://dx.doi.org/10.1137/15M1052834}{\nolinkurl{doi:10.1137/15m1052834}}.

\bibitem{conn1991globally}
A.\,R{.\nobreak\kern 0.33333em}Conn, N.\,I.\,M{.\nobreak\kern 0.33333em}Gould,
  and P.\,L{.\nobreak\kern 0.33333em}Toint, A globally convergent augmented
  {L}agrangian algorithm for optimization with general constraints and simple
  bounds, \emph{SIAM Journal on Numerical Analysis} 28 (1991),  545--572,
  \href{https://dx.doi.org/10.1137/0728030}{\nolinkurl{doi:10.1137/0728030}}.

\bibitem{demarchi2024implicit}
A{.\nobreak\kern 0.33333em}De~Marchi, Implicit augmented {L}agrangian and
  generalized optimization, \emph{Journal of Applied and Numerical
  Optimization} 6 (2024),  291--320,
  \href{https://dx.doi.org/10.23952/jano.6.2024.2.08}{\nolinkurl{doi:10.23952/jano.6.2024.2.08}}.

\bibitem{demarchi2023constrained}
A{.\nobreak\kern 0.33333em}De~Marchi, X{.\nobreak\kern 0.33333em}Jia,
  C{.\nobreak\kern 0.33333em}Kanzow, and P{.\nobreak\kern 0.33333em}Mehlitz,
  Constrained composite optimization and augmented {L}agrangian methods,
  \emph{Mathematical Programming} 201 (2023),  863--896,
  \href{https://dx.doi.org/10.1007/s10107-022-01922-4}{\nolinkurl{doi:10.1007/s10107-022-01922-4}}.

\bibitem{demarchi2022proximal}
A{.\nobreak\kern 0.33333em}De~Marchi and A{.\nobreak\kern 0.33333em}Themelis,
  Proximal gradient algorithms under local {L}ipschitz gradient continuity,
  \emph{Journal of Optimization Theory and Applications} 194 (2022),  771--794,
  \href{https://dx.doi.org/10.1007/s10957-022-02048-5}{\nolinkurl{doi:10.1007/s10957-022-02048-5}}.

\bibitem{DempeDuttaMordukhovich2007}
S{.\nobreak\kern 0.33333em}Dempe, J{.\nobreak\kern 0.33333em}Dutta, and
  B.\,S{.\nobreak\kern 0.33333em}Mordukhovich, New necessary optimality
  conditions in optimistic bilevel programming, \emph{Optimization} 56 (2007),
  577--604,
  \href{https://dx.doi.org/10.1080/02331930701617551}{\nolinkurl{doi:10.1080/02331930701617551}}.

\bibitem{DhingraKhongJavanovic2019}
N.\,K{.\nobreak\kern 0.33333em}Dhingra, S.\,Z{.\nobreak\kern 0.33333em}Khong,
  and M.\,R{.\nobreak\kern 0.33333em}Jovanovi{\'{c}}, The proximal augmented
  {L}agrangian method for nonsmooth composite optimization, \emph{IEEE
  Transactions on Automatic Control} 64 (2019),  2861--2868,
  \href{https://dx.doi.org/10.1109/TAC.2018.2867589}{\nolinkurl{doi:10.1109/tac.2018.2867589}}.

\bibitem{DontchevRockafellar2014}
A.\,L{.\nobreak\kern 0.33333em}Dontchev and R.\,T{.\nobreak\kern
  0.33333em}Rockafellar, \emph{Implicit Functions and Solution Mappings},
  Springer, Heidelberg, 2014,
  \href{https://dx.doi.org/10.1007/978-0-387-87821-8}{\nolinkurl{doi:10.1007/978-0-387-87821-8}}.

\bibitem{FernandezSolodov2012}
D{.\nobreak\kern 0.33333em}Fern{\'a}ndez and M{.\nobreak\kern
  0.33333em}Solodov, Local convergence of exact and inexact augmented
  {L}agrangian methods under the second-order sufficient optimality condition,
  \emph{SIAM Journal on Optimization} 22 (2012),  384--407,
  \href{https://dx.doi.org/10.1137/10081085X}{\nolinkurl{doi:10.1137/10081085x}}.

\bibitem{Gfrerer2014}
H{.\nobreak\kern 0.33333em}Gfrerer, Optimality conditions for disjunctive
  programs based on generalized differentiation with application to
  mathematical programs with equilibrium constraints, \emph{SIAM Journal on
  Optimization} 24 (2014),  898--931,
  \href{https://dx.doi.org/10.1137/130914449}{\nolinkurl{doi:10.1137/130914449}}.

\bibitem{GfrererOutrata2016}
H{.\nobreak\kern 0.33333em}Gfrerer and J.\,V{.\nobreak\kern 0.33333em}Outrata,
  On computation of generalized derivatives of the normal-cone mapping and
  their applications, \emph{Mathematics of Operations Research} 41 (2016),
  1535--1556,
  \href{https://dx.doi.org/10.1287/moor.2016.0789}{\nolinkurl{doi:10.1287/moor.2016.0789}}.

\bibitem{GuoLinYe2013}
L{.\nobreak\kern 0.33333em}Guo, G.\,H{.\nobreak\kern 0.33333em}Lin, and
  J.\,J{.\nobreak\kern 0.33333em}Ye, Second-order optimality conditions for
  mathematical programs with equilibrium constraints, \emph{Journal of
  Optimization Theory and Applications} 158 (2013),  33--64,
  \href{https://dx.doi.org/10.1007/s10957-012-0228-x}{\nolinkurl{doi:10.1007/s10957-012-0228-x}}.

\bibitem{HangMordukhovichSarabi2020}
N.\,T.\,V{.\nobreak\kern 0.33333em}Hang, B.\,S{.\nobreak\kern
  0.33333em}Mordukhovich, and M.\,E{.\nobreak\kern 0.33333em}Sarabi,
  Second-order variational analysis in second-order cone programming,
  \emph{Mathematical Programming} 180 (2020),  75--116,
  \href{https://dx.doi.org/10.1007/s10107-018-1345-6}{\nolinkurl{doi:10.1007/s10107-018-1345-6}}.

\bibitem{HangMordukhovichSarabi2022}
N.\,T.\,V{.\nobreak\kern 0.33333em}Hang, B.\,S{.\nobreak\kern
  0.33333em}Mordukhovich, and M.\,E{.\nobreak\kern 0.33333em}Sarabi, Augmented
  {L}agrangian method for second-order cone programs under second-order
  sufficiency, \emph{Journal of Global Optimization} 82 (2022),  51--81,
  \href{https://dx.doi.org/10.1007/s10898-021-01068-1}{\nolinkurl{doi:10.1007/s10898-021-01068-1}}.

\bibitem{HangSarabi2021}
N.\,T.\,V{.\nobreak\kern 0.33333em}Hang and M.\,E{.\nobreak\kern
  0.33333em}Sarabi, Local convergence analysis of augmented {L}agrangian
  methods for piecewise linear-quadratic composite optimization problems,
  \emph{SIAM Journal on Optimization} 31 (2021),  2665--2694,
  \href{https://dx.doi.org/10.1137/20M1375188}{\nolinkurl{doi:10.1137/20m1375188}}.

\bibitem{hang2023convergence}
N.\,T.\,V{.\nobreak\kern 0.33333em}Hang and M.\,E{.\nobreak\kern
  0.33333em}Sarabi, Convergence of augmented {L}agrangian methods for composite
  optimization problems, \emph{arXiv}  (2023),
  \href{https://arxiv.org/abs/2307.15627v2}{\nolinkurl{arXiv:2307.15627v2}}.

\bibitem{hestenes1969multiplier}
M.\,R{.\nobreak\kern 0.33333em}Hestenes, Multiplier and gradient methods,
  \emph{Journal of Optimization Theory and Applications} 4 (1969),  303--320,
  \href{https://dx.doi.org/10.1007/BF00927673}{\nolinkurl{doi:10.1007/bf00927673}}.

\bibitem{IoffeOutrata2008}
A.\,D{.\nobreak\kern 0.33333em}Ioffe and J.\,V{.\nobreak\kern
  0.33333em}Outrata, On metric and calmness qualification conditions in
  subdifferential calculus, \emph{Set-Valued Analysis} 16 (2008),  199--227,
  \href{https://dx.doi.org/10.1007/s11228-008-0076-x}{\nolinkurl{doi:10.1007/s11228-008-0076-x}}.

\bibitem{Izmailov2005}
A.\,F{.\nobreak\kern 0.33333em}Izmailov, On the analytical and numerical
  stability of critical {L}agrange multiplier, \emph{Computational Mathematics
  and Mathematical Physics} 45 (2005),  930--946,
  \url{https://www.mathnet.ru/php/archive.phtml?wshow=paper&jrnid=zvmmf&paperid=636&option\_lang=eng}.

\bibitem{IzmailovKurennoySolodov2015}
A.\,F{.\nobreak\kern 0.33333em}Izmailov, A.\,S{.\nobreak\kern
  0.33333em}Kurennoy, and M.\,V{.\nobreak\kern 0.33333em}Solodov, Local
  convergence of the method of multipliers for variational and optimization
  problems under the noncriticality assumption, \emph{Computational
  Optimization and Applications} 60 (2015),  111--140,
  \href{https://dx.doi.org/10.1007/s10589-014-9658-8}{\nolinkurl{doi:10.1007/s10589-014-9658-8}}.

\bibitem{IzmailovSolodov2012}
A.\,F{.\nobreak\kern 0.33333em}Izmailov and M.\,V{.\nobreak\kern
  0.33333em}Solodov, Stabilized {S}{Q}{P} revisited, \emph{Mathematical
  Programming} 133 (2012),  93--120,
  \href{https://dx.doi.org/10.1007/s10107-010-0413-3}{\nolinkurl{doi:10.1007/s10107-010-0413-3}}.

\bibitem{JiaKanzowMehlitzWachsmuth2023}
X{.\nobreak\kern 0.33333em}Jia, C{.\nobreak\kern 0.33333em}Kanzow,
  P{.\nobreak\kern 0.33333em}Mehlitz, and G{.\nobreak\kern 0.33333em}Wachsmuth,
  An augmented {L}agrangian method for optimization problems with structured
  geometric constraints, \emph{Mathematical Programming} 199 (2023),
  1365--1415,
  \href{https://dx.doi.org/10.1007/s10107-022-01870-z}{\nolinkurl{doi:10.1007/s10107-022-01870-z}}.

\bibitem{KanzowSteck2018}
C{.\nobreak\kern 0.33333em}Kanzow and D{.\nobreak\kern 0.33333em}Steck, On
  error bounds and multiplier methods for variational problems in {B}anach
  spaces, \emph{SIAM Journal on Control and Optimization} 56 (2018),
  1716--1738,
  \href{https://dx.doi.org/10.1137/17M1146518}{\nolinkurl{doi:10.1137/17m1146518}}.

\bibitem{KanzowSteck2019}
C{.\nobreak\kern 0.33333em}Kanzow and D{.\nobreak\kern 0.33333em}Steck,
  Improved local convergence results for augmented {L}agrangian methods in
  {$C^2$}-reducible constrained optimization, \emph{Mathematical Programming}
  177 (2019),  425--438,
  \href{https://dx.doi.org/10.1007/s10107-018-1261-9}{\nolinkurl{doi:10.1007/s10107-018-1261-9}}.

\bibitem{kanzow2018augmented}
C{.\nobreak\kern 0.33333em}Kanzow, D{.\nobreak\kern 0.33333em}Steck, and
  D{.\nobreak\kern 0.33333em}Wachsmuth, An augmented {L}agrangian method for
  optimization problems in {B}anach spaces, \emph{SIAM Journal on Control and
  Optimization} 56 (2018),  272--291,
  \href{https://dx.doi.org/10.1137/16M1107103}{\nolinkurl{doi:10.1137/16m1107103}}.

\bibitem{Levy1996}
A.\,B{.\nobreak\kern 0.33333em}Levy, Implicit multifunction theorems for the
  sensitivity analysis of variational conditions, \emph{Mathematical
  Programming} 74 (1996),  333--350,
  \href{https://dx.doi.org/10.1007/BF02592203}{\nolinkurl{doi:10.1007/bf02592203}}.

\bibitem{LuoPangRalph1996}
Z.\,Q{.\nobreak\kern 0.33333em}Luo, J.\,S{.\nobreak\kern 0.33333em}Pang, and
  D{.\nobreak\kern 0.33333em}Ralph, \emph{Mathematical Programs with
  Equilibrium Constraints}, Cambridge University Press, Cambridge, 1996,
  \href{https://dx.doi.org/10.1017/CBO9780511983658}{\nolinkurl{doi:10.1017/cbo9780511983658}}.

\bibitem{mehlitz2020}
P{.\nobreak\kern 0.33333em}Mehlitz, On the linear independence constraint
  qualification in disjunctive programming, \emph{Optimization} 69 (2020),
  2241--2277,
  \href{https://dx.doi.org/10.1080/02331934.2019.1679811}{\nolinkurl{doi:10.1080/02331934.2019.1679811}}.

\bibitem{MohammadiMordukhovichSarabi2022}
A{.\nobreak\kern 0.33333em}Mohammadi, B.\,S{.\nobreak\kern
  0.33333em}Mordukhovich, and M.\,E{.\nobreak\kern 0.33333em}Sarabi,
  Variational analysis of composite models with applications to continuous
  optimization, \emph{Mathematics of Operations Research} 47 (2022),  397--426,
  \href{https://dx.doi.org/10.1287/moor.2020.1074}{\nolinkurl{doi:10.1287/moor.2020.1074}}.

\bibitem{mordukhovich2018}
B.\,S{.\nobreak\kern 0.33333em}Mordukhovich, \emph{Variational Analysis and
  Applications}, Springer, Cham, 2018,
  \href{https://dx.doi.org/10.1007/978-3-319-92775-6}{\nolinkurl{doi:10.1007/978-3-319-92775-6}}.

\bibitem{MordukhovichSarabi2018}
B.\,S{.\nobreak\kern 0.33333em}Mordukhovich and M.\,E{.\nobreak\kern
  0.33333em}Sarabi, Critical multipliers in variational systems via
  second-order generalized differentiation, \emph{Mathematical Programming} 169
  (2018),  605--648,
  \href{https://dx.doi.org/10.1007/s10107-017-1155-2}{\nolinkurl{doi:10.1007/s10107-017-1155-2}}.

\bibitem{MordukhovichSarabi2019}
B.\,S{.\nobreak\kern 0.33333em}Mordukhovich and M.\,E{.\nobreak\kern
  0.33333em}Sarabi, Criticality of {L}agrange multipliers in variational
  systems, \emph{SIAM Journal on Optimization} 29 (2019),  1524--1557,
  \href{https://dx.doi.org/10.1137/18M1206862}{\nolinkurl{doi:10.1137/18m1206862}}.

\bibitem{OutrataKocvaraZowe1998}
J.\,V{.\nobreak\kern 0.33333em}Outrata, M{.\nobreak\kern
  0.33333em}Ko{\v{c}}vara, and J{.\nobreak\kern 0.33333em}Zowe, \emph{Nonsmooth
  Approach to Optimization Problems with Equilibrium Constraints}, Kluwer
  Academic, Dordrecht, 1998,
  \href{https://dx.doi.org/10.1007/978-1-4757-2825-5}{\nolinkurl{doi:10.1007/978-1-4757-2825-5}}.

\bibitem{powell1969method}
M.\,J.\,D{.\nobreak\kern 0.33333em}Powell, A method for nonlinear constraints
  in minimization problems, in \emph{Optimization}, Academic Press, 1969,
  283--298.

\bibitem{Robinson1981}
S.\,M{.\nobreak\kern 0.33333em}Robinson, Some continuity properties of
  polyhedral multifunctions, in \emph{Mathematical Programming at
  {O}berwolfach}, H{.\nobreak\kern 0.33333em}K{\"o}nig, B{.\nobreak\kern
  0.33333em}Korte, and K{.\nobreak\kern 0.33333em}Ritter (eds.), Springer,
  Berlin, 1981,  206--214,
  \href{https://dx.doi.org/10.1007/BFb0120929}{\nolinkurl{doi:10.1007/bfb0120929}}.

\bibitem{rockafellar1973dual}
R.\,T{.\nobreak\kern 0.33333em}Rockafellar, A dual approach to solving
  nonlinear programming problems by unconstrained optimization,
  \emph{Mathematical Programming} 5 (1973),  354--373,
  \href{https://dx.doi.org/10.1007/BF01580138}{\nolinkurl{doi:10.1007/bf01580138}}.

\bibitem{rockafellar1976augmented}
R.\,T{.\nobreak\kern 0.33333em}Rockafellar, Augmented {L}agrangians and
  applications of the proximal point algorithm in convex programming,
  \emph{Mathematics of Operations Research} 1 (1976),  97--116,
  \href{https://dx.doi.org/10.1287/moor.1.2.97}{\nolinkurl{doi:10.1287/moor.1.2.97}}.

\bibitem{rockafellar1993lagrange}
R.\,T{.\nobreak\kern 0.33333em}Rockafellar, Lagrange multipliers and
  optimality, \emph{SIAM Review} 35 (1993),  183--238,
  \href{https://dx.doi.org/10.1137/1035044}{\nolinkurl{doi:10.1137/1035044}}.

\bibitem{rockafellar2000extended}
R.\,T{.\nobreak\kern 0.33333em}Rockafellar, Extended nonlinear programming, in
  \emph{Nonlinear Optimization and Related Topics}, G{.\nobreak\kern
  0.33333em}Di~Pillo and F{.\nobreak\kern 0.33333em}Giannessi (eds.), volume~36
  of Applied Optimization, Kluwer, 2000,  381--399,
  \href{https://dx.doi.org/10.1007/978-1-4757-3226-9_20}{\nolinkurl{doi:10.1007/978-1-4757-3226-9_20}}.

\bibitem{rockafellar2022augmented}
R.\,T{.\nobreak\kern 0.33333em}Rockafellar, Augmented {L}agrangians and hidden
  convexity in sufficient conditions for local optimality, \emph{Mathematical
  Programming} 198 (2023),  159--194,
  \href{https://dx.doi.org/10.1007/s10107-022-01768-w}{\nolinkurl{doi:10.1007/s10107-022-01768-w}}.

\bibitem{rockafellar2022convergence}
R.\,T{.\nobreak\kern 0.33333em}Rockafellar, Convergence of augmented
  {L}agrangian methods in extensions beyond nonlinear programming,
  \emph{Mathematical Programming} 199 (2023),  375--420,
  \href{https://dx.doi.org/10.1007/s10107-022-01832-5}{\nolinkurl{doi:10.1007/s10107-022-01832-5}}.

\bibitem{rockafellar1998variational}
R.\,T{.\nobreak\kern 0.33333em}Rockafellar and R.\,J.\,B{.\nobreak\kern
  0.33333em}Wets, \emph{Variational Analysis}, Springer, 1998,
  \href{https://dx.doi.org/10.1007/978-3-642-02431-3}{\nolinkurl{doi:10.1007/978-3-642-02431-3}}.

\bibitem{sabach2019lagrangian}
S{.\nobreak\kern 0.33333em}Sabach and M{.\nobreak\kern 0.33333em}Teboulle,
  Lagrangian methods for composite optimization, in \emph{Processing, Analyzing
  and Learning of Images, Shapes, and Forms: Part 2}, R{.\nobreak\kern
  0.33333em}Kimmel and X.\,C{.\nobreak\kern 0.33333em}Tai (eds.), volume~20 of
  Handbook of Numerical Analysis, Elsevier, 2019,  401--436,
  \href{https://dx.doi.org/10.1016/bs.hna.2019.04.002}{\nolinkurl{doi:10.1016/bs.hna.2019.04.002}}.

\bibitem{steck2018dissertation}
D{.\nobreak\kern 0.33333em}Steck, \emph{{L}agrange Multiplier Methods for
  Constrained Optimization and Variational Problems in {B}anach Spaces}, PhD
  thesis, University of W\"urzburg, 2018,
  \url{https://nbn-resolving.org/urn:nbn:de:bvb:20-opus-174444}.

\bibitem{Wachsmuth2017}
G{.\nobreak\kern 0.33333em}Wachsmuth, Strong stationarity for optimization
  problems with complementarity constraints in absence of polyhedricity,
  \emph{Set-Valued and Variational Analysis} 25 (2017),  133--175,
  \href{https://dx.doi.org/10.1007/s11228-016-0370-y}{\nolinkurl{doi:10.1007/s11228-016-0370-y}}.

\bibitem{WalkupWets1969}
D.\,W{.\nobreak\kern 0.33333em}Walkup and R.\,J.\,B{.\nobreak\kern
  0.33333em}Wets, A {L}ipschitz characterization of convex polyhedra,
  \emph{Proceedings of the American Mathematical Society} 23 (1969),  167--173,
  \href{https://dx.doi.org/10.2307/2037511}{\nolinkurl{doi:10.2307/2037511}}.

\end{thebibliography}
}%

\end{document}